\documentclass[english, hidelinks, 12pt]{article}
\usepackage{geometry}
\usepackage{setspace}
\usepackage{authblk}
\usepackage[utf8]{inputenc}
\usepackage{color,soul,yhmath}
\usepackage{xcolor}
\usepackage{bm}
\usepackage{graphicx}
\graphicspath{ {.} }
\usepackage{hyperref}
\usepackage{cleveref}
\usepackage{natbib}
\usepackage{bbm}
\usepackage{amsmath}
\usepackage{amssymb}
\usepackage{amsthm}
\usepackage{booktabs}
\usepackage{multirow}

\usepackage{enumitem}

\usepackage[makeroom]{cancel}

\hypersetup{
    colorlinks=true,
    linkcolor=blue,
    citecolor=blue,
    filecolor=magenta,
    urlcolor=blue
}

\newtheorem{theorem}{Theorem}

\newtheorem{corollary}{Corollary}

\newtheorem{lemma}{Lemma}

\newtheorem{proposition}{Proposition}
\newtheorem{remark}{Remark}

\numberwithin{equation}{section}
\def\c#1{\mathcal{#1}}
\def\b#1{\mathbbm{#1}}
\def\bf#1{\mathbf{#1}}

\def\wt#1{\widetilde{#1}}
\def\wh#1{\widehat{#1}}
\def\vec#1{\textnormal{\textbf{vec}}(#1)}

\def\top{\intercal}

\def\cA{\mathcal{A}}
\def\cB{\mathcal{B}}

\def\cI{\mathcal{I}}

\def\cL{\mathcal{L}}

\def\cN{\mathcal{N}}
\def\cQ{\mathcal{Q}}
\def\cS{\mathcal{S}}

\def\aa{\mathbf{a}}
\def\bb{\mathbf{b}}
\def\g{\mathbf{g}}
\def\h{\mathbf{h}}

\def\v{\mathbf{v}}
\def\x{\mathbf{x}}
\def\y{\mathbf{y}}

\def\A{\mathbf{A}}
\def\B{\mathbf{B}}
\def\C{\mathbf{C}}

\def\Q{\mathbf{Q}}

\def\D{\mathbf{D}}

\def\M{\mathbf{M}}
\def\F{\mathbf{F}}
\def\E{\mathbf{E}}
\def\H{\mathbf{H}}
\def\U{\mathbf{U}}

\def\R{\mathbf{R}}
\def\I{\mathbf{I}}
\def\L{\mathbf{L}}

\def\X{\mathbf{X}}
\def\Z{\mathbf{Z}}
\def\0{\mathbf{0}}
\def\1{\mathbf{1}}
\def\diag{\textnormal{\text{diag}}}

\def\bSigma{\boldsymbol{\Sigma}}

\def\balpha{\boldsymbol{\alpha}}
\def\bbeta{\boldsymbol{\beta}}

\def\btheta{\boldsymbol{\theta}}

\def\bepsilon{\boldsymbol{\epsilon}}
\def\var{\textnormal{Var}}

\newcommand{\norm}[1]{\big\| #1 \big\|}

\def\MSE{\operatorname{MSE}}
\def\med{\operatorname{med}}

\allowdisplaybreaks[1]

\usepackage{parskip}

\begin{document}
\def\spacingset#1{\renewcommand{\baselinestretch}%
{#1}\small\normalsize} \spacingset{1}

	\begin{titlepage}
		
		\title{Main Effect Factor Models in High-Dimensional Matrix Time Series: Identification and Sparsity}
                \author{Zetai Cen\thanks{Equal contribution; corresponding author. School of Mathematics, University of Bristol.
                Email: \url{zetai.cen@bristol.ac.uk}.
                Supported by Engineering and Physical Sciences Research Council (EP/Z531327/1).}}
                \author{Kaixin Liu\thanks{Equal contribution. Department of Statistics, London School of Economics. 
                Email: \url{K.Liu31@lse.ac.uk}.}}
                \author{Clifford Lam\thanks{Department of Statistics, London School of Economics. 
                Email: \url{C.Lam2@lse.ac.uk}.}}

		\affil{}
		
		\date{}
		
		\maketitle

\begin{abstract}
We propose a general identification framework for main effect factor models for matrix-valued time series. The classical sum‑to‑zero restriction on the row and column main effects is replaced by a broad class of shift‑varying functions, which includes weighted averages, quantiles, and reference‑unit anchors as special cases. We prove that the shift‑varying property is necessary and sufficient for parameter identification, and we derive closed‑form estimators under minimal assumptions. Asymptotic convergence rates of all estimated components are derived under weak, heterogeneous factor strengths. Building on this flexible identification, we address sparsity in the main effects by selecting a minimum‑based shift‑varying function that sets the smallest main effect to zero, and we introduce a doubly adaptive fused Lasso estimator that consistently recovers the true sparse and dense blocks. A modified Mallows's statistic is developed for tuning parameter selection, and the entire procedure is computationally practical via an equivalent generalized lasso formulation. Simulation experiments are performed under a variety of settings, showing our proposed methods work well. Two real applications on employment growth and economic indices are demonstrated to illustrate the empirical value of our approach.
\end{abstract}
		
		\bigskip
		\bigskip


\noindent

	\end{titlepage}
	
	\setcounter{page}{2}

\maketitle


\renewcommand{\baselinestretch}{2}

\newpage
\section{Introduction}\label{sec: introduction}

With the rapid development of computational hardware and technology for big data, researchers obtain and analyze datasets that are ever larger in size and complexity. This leads to significant advances in the area of factor analysis, which is a useful tool in multivariate analysis with a wide range of applications in psychology \citep{McCraeJohn1992}, biology \citep{Hirzeletal2002, Hochreiteretal2006}, economics and finance \citep{FamaFrench1993, StockWatson2002, StockWatson2002b}, to name but a few areas. In particular, since the early work of \cite{ChamberlainRothchild1983}, approximate factor models have been well studied over the past few decades; see \text{e.g.} \cite{BaiNg2002}, \cite{PanYao2008}, \cite{Lametal2011}, and the references therein.

\paragraph{Related literature and motivation.}
To facilitate interpretation of the estimated factor structure, one of the major solutions is through sparsity, which is not new in time series analysis. Sparsity can be imposed on the data covariance matrix \citep{BickelLevina2008} and the noise covariance matrix \citep{Fanetal2013}, among many other methods. See also Section~1.4 in \cite{UematsuYamagata2023} for a discussion on sparse principal components. More recently, researchers are interested in studying the sparsity in factor loadings \citep{Freyaldenhoven2022, UematsuYamagata2023}, which is a concept closely related to factor strength as discussed in Section~7 in \cite{BarigozziHallin2024}. In fact, various forms of sparsity in factor loadings are discussed in the literature. An example is a multilevel/group factor model \citep{Wang2008, Huetal2025}, where each observed vector $\x_t^s\in \b{R}^{p_s}$ ($t=1,\dots,T$, $s=1,\dots,S$) can be represented by $\x_t^s = \A^s \g_t + \B^s \mathbf{f}_t^s +\mathbf{e}_t^s$,  with $\g_t$ and $\mathbf{f}_t^s$ called the global and group-specific factors respectively, $\A^s$, $\B^s$ their corresponding factor loading matrices, and $\mathbf{e}_t^s$ the noise. The model can be rewritten as
\[
\begin{pmatrix}
    \x_t^1 \\ \x_t^2 \\ \vdots \\ \x_t^S
\end{pmatrix} = \begin{pmatrix}
    \A^1 & \B^1 & \0 & \ldots & \0 \\
    \A^2 & \0 & \B^2 & \ldots & \0 \\
    \vdots & \vdots & \vdots & \ddots & \vdots \\
    \A^S & \0 & \0 & \ldots & \B^S
\end{pmatrix} \begin{pmatrix}
    \g_t \\ \mathbf{f}_t^1 \\ \vdots \\ \mathbf{f}_t^S
\end{pmatrix} + \begin{pmatrix}
    \mathbf{e}_t^1 \\ \mathbf{e}_t^2 \\ \vdots \\ \mathbf{e}_t^S
\end{pmatrix} ,
\]
so that the factor loading matrix has a specific sparse structure, which is based on a priori data grouping. Another similar example of sparse loadings is models of the above kind, except that the membership of the data is unknown; see \cite{AndoBai2017} and \cite{Zhangetal2024}. In addition to block sparsity due to grouping, \cite{UematsuYamagata2023, UematsuYamagata2023_infer} studies a sparsity-induced weak factor model under rather restricted conditions to address the identification issue that sparse factor loadings are generally not rotation-invariant. \cite{WeiZhang2024} further investigates the near-sparsity preservation property of the estimated loadings. Other related examples include \cite{Fanetal2023} which relates sparsity and factor models through factor-augmented regression, and \cite{Mosleyetal2024} which focuses on sparsity in the loading matrix of dynamic factor models. For Bayesian methods on sparsity in factor structures, see \cite{Zhangetal2025} and the references therein.

Despite the growing literature from the above discussion, all of them focus only on vector-valued time series. As first proposed by \cite{Wangetal2019} and later studied in broader literature \citep[e.g.][]{Yuetal2022, Heetal2024}, matrix factor models have become powerful tools for studying economic and financial data, leading to more insightful interpretations and improved estimations. In particular, a series of matrices $\X_t$, $t\in \{1,\ldots,T\} =: [T]$, is observed and admits the decomposition $\X_t = \A_r \F_t \A_c^\top +\E_t$, where $\A_r,\A_c$ are the row and column factor loading matrices respectively, and $\F_t$ is the core factor matrix. Although one may write $\vec{\X_t} = (\A_c\otimes \A_r) \vec{\F_t} + \vec{\E_t}$, where $\vec{\X_t}$ is the stacking of the columns of $\X_t$ into a vector, and thus estimate a factor model for the vectorized data with sparsity, such an approach is neither efficient nor appropriate. For instance, the global inference step in \cite{UematsuYamagata2023_infer} is now inapplicable as it neglects the Kronecker product structure in the factor loading matrix.

This paper provides an alternative solution to address sparsity in matrix factor models. It is worth pointing out that it remains a challenge to address the interaction between rows and columns in factor modeling for matrix-valued time series, for which very limited work has been done. An example is \cite{Heetal2023_onewaytwoway} where a specification test is proposed. Another attempt is by \cite{lam2024matrix}, which generalizes the Tucker-decomposition factor model to a Main Effect Factor Model (MEFM), significantly improving model interpretability. In particular, given a sequence of mean-zero matrix observations $\X_t\in \b{R}^{p\times q}$ for $t=1,\ldots,T$, MEFM decomposes each $\X_t$ as
\begin{equation}
\label{eqn: MEFM}
\X_t = \mu_t\1_p \1_q^\top + \balpha_t^\ast \1_q^\top + \1_p \bbeta_t^{\ast \top} +\C_t + \E_t,
\end{equation}
where $\mu_t \in\b{R}$ is the grand mean,
$\balpha_t^\ast \in\b{R}^p$ and $\bbeta_t^\ast \in\b{R}^q$ are respectively the row and column main effects with potentially many zero entries, $\C_t =\A_r \F_t \A_c^\top$ is the common component, $\A_r\in \b{R}^{p\times k_r}$ and $\A_c\in \b{R}^{q\times k_c}$ are the row and column factor loading matrices, $\F_t$ is the core factor matrix, and $\E_t$ is the noise.
In particular, the identification of parameters, \text{e.g.} main effects and factor loadings, relies heavily on the condition
\begin{equation}
\label{eqn: MEFM_identification}
\1_p^\top \balpha_t^\ast =0, \;\;\;
\1_q^\top \bbeta_t^\ast =0, \;\;\;
\1_p^\top \A_r =\0, \;\;\;
\1_q^\top \A_c =\0 .
\end{equation}
\cite{lam2024matrix} shows that any Tucker-decomposition matrix factor model can be rewritten into MEFM, whereas the converse generally requires far more factors and can still be empirically inferior. Note that for any fixed row (\text{resp.} column), the contribution of $\balpha_t^\ast$ (\text{resp.} $\bbeta_t^\ast$) towards the observed $\X_t$ is constant across columns (\text{resp.} rows). Hence row or column structures are solely featured by the main effects, and the interaction between rows and columns by the common component.

However, the identification condition \eqref{eqn: MEFM_identification} gives rise to complications when sparsity of the main effects is desired. 
For instance, suppose $\X_t$ records the values of economic indicators with countries indexed by rows and indicators indexed by columns. If a small group of countries form an economic entity with pervasive effects on each indicator of the member countries, then the row main effects vector $\balpha_t^\ast$ is sparse with non-zero entries corresponding to those member countries.
Yet \eqref{eqn: MEFM_identification} implies that either some member countries have opposite effects which are unnatural, or some non-member countries have non-zero effects which are not realistic under this scenario.
To resolve this, we propose a flexible class of identification conditions based on \emph{shift-varying} functions. Such a function satisfies that adding a constant to its argument strictly changes its value; typical examples include weighted sums (yielding the sum‑to‑zero condition as a special case) and quantile functions (e.g., anchoring the main effects at their minimum or median). 
By choosing a shift‑varying function that sets a particular quantile to zero, one can naturally accommodate sparse and interpretable main effects, while the sum‑to‑zero condition becomes only one of many possibilities. This new identification significantly enhances the interpretability of the model and is the foundation of our sparsity recovery methodology.

\begin{figure}[t!]
    \centering
    \includegraphics[width=\textwidth]{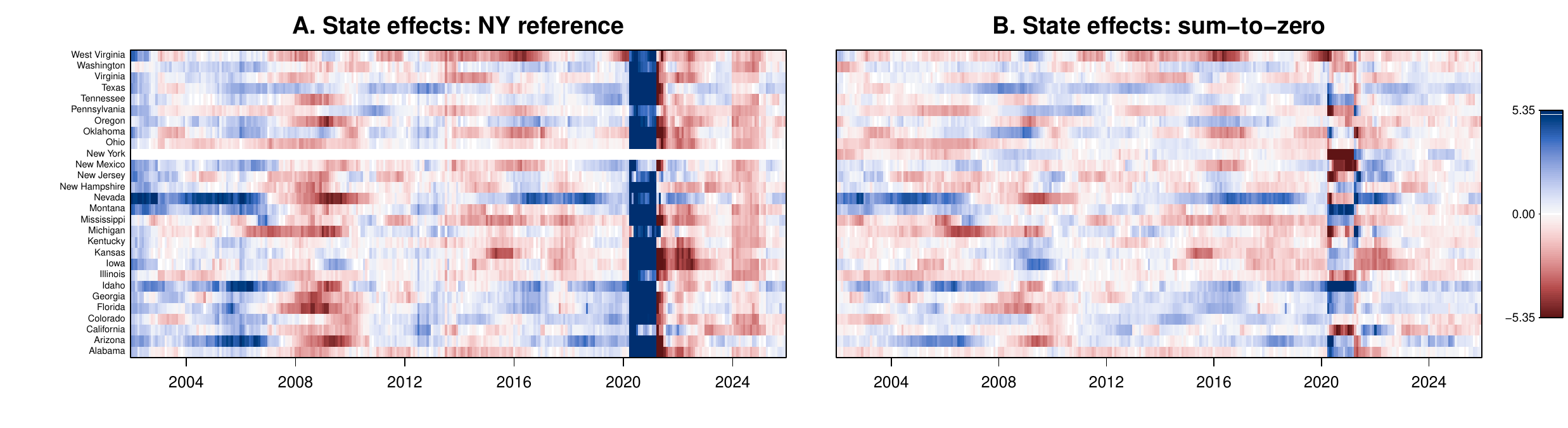}
    \caption{QCEW state effects under identifications of \textit{NY reference} (panel~A) and \textit{sum-to-zero} (panel~B).
    Color scales are symmetrically clipped at the 98th percentile of absolute effect magnitudes for readability.}
    \label{fig:qcew_state_effect}
\end{figure}

As an empirical example, consider the monthly employment data from the U.S.\ Bureau of Labor Statistics Quarterly Census of Employment and Wages (QCEW).
The dataset contains $T=288$ months from 2002-M01 to 2025-M12, and
at each month $t$, $\X_t\in\b{R}^{28\times 17}$ denotes the employment growth for 28 states and 17 non-overlapping private-sector industries;
see Section~\ref{subsec: application_QCEW} for detailed investigations and we only highlight the value of flexible identifications here.
The row main effects $\balpha_t^\ast \in 
\b{R}^{28}$, namely the state effects, are estimated under two identifications for all $t\in[288]$: 
\textit{NY reference} where the element of $\balpha_t^\ast$ corresponding to New York is zero;
\textit{sum-to-zero} such that $\1_{28}^\top \balpha_t^\ast=0$ as in \cite{lam2024matrix}.
While the interpretation of \textit{sum-to-zero} is not straightforward, the choice of \textit{NY reference} is already meaningful such that the state effects represent the \textit{excess} employment growth of each state with respect to the growth in New York.
Later, we will show that \textit{NY reference} is a valid identification condition for the main effects.
We also show in Figure~\ref{fig:qcew_state_effect} the estimated effects under these two settings.
The comparison is sheer around the pandemic period from 2020: 
\textit{NY reference} clearly shows excess employment growth in all other states, reflecting the significant employment crash in New York, which is vague under \textit{sum-to-zero}.
From \textit{NY reference},
the labor market in New York also seems to experience a post-pandemic rebound during 2024,
which again shows the benefits of considering such identification.

\paragraph{Contributions and organization.}
Built upon the existing MEFM framework, we develop a general theory for its identification that allows researchers to anchor the row and column main effects via any shift‑varying function. This class includes the sum‑to‑zero restriction of \cite{lam2024matrix} as a special case, as well as quantile‑based and reference‑unit (viz.\ a particular cross‑sectional unit’s effect is fixed to be zero) identifications that are far better flexible and interpretable in different scenarios, and are suited to possibly sparse or tail‑heavy data.
We show that the shift‑varying condition is necessary and sufficient for parameter identification, derive closed‑form estimators, and establish their asymptotic properties under weak, heterogeneous factor strengths.
As a first in the literature, we consider MEFM with sparsity in the main effects and develop methods to recover the structure for more natural and better interpretation of the main effects over time. A special (and sparsest) case is the traditional matrix factor model where all main effects are zero. 
A period of non-zero main effects in relation to other periods of zero main effects can indicate a significant change in circumstances for the rows or columns in question. 
By employing a quantile‑based identification (e.g., zero‑minimum constraint), we can directly interpret zero main effects as the absence of row‑ or column‑specific deviations, and the sparsity pattern itself becomes a quantity of interest.

\iftrue
The rest of this paper is organized as follows. The new model identification is detailed in Section~\ref{sec: model_estimation}, followed by parameter estimation. 
Section~\ref{sec: assumption_theories} presents the main assumptions and theoretical results of the estimators. 
Section~\ref{sec: sparsity_ME} develops regularized estimation for sparse main effects and studies the asymptotic properties of the proposed method for sparsity recovery. 
We also discuss the implementation of the algorithm and hyperparameter tuning based on a modified $C_p$ statistic. 
Lastly, numerical results are shown in Section~\ref{sec: numerical_results}, where two real data applications are also analyzed. 
Section~\ref{sec: conclusion} concludes the paper.
All theoretical proofs and additional information are relegated to the supplement.
\fi

\paragraph{Notations.}
We use the lower-case or capital letter, bold lower-case letter, and bold capital letter, i.e., $a$ or $A$, $\bf{a}$, $\bf{A}$, to denote a scalar, a vector, and a matrix respectively. We also use $a_i, A_{i,j}, \bf{A}_{i\cdot}, \bf{A}_{\cdot i}$ to denote, respectively, the $i$-th element of $\bf{a}$, the $(i,j)$-th element of $\bf{A}$, the $i$-th row vector (as a column vector) of $\bf{A}$, and the $i$-th column vector of $\bf{A}$. We use $\circ$ to denote the Hadamard product and $\otimes$ the Kronecker product; $a\asymp b$ represents $a=O(b)$ and $b=O(a)$;
$a\lor b:=\max(a,b)$,
and $a \wedge b := \min(a,b)$.
{For a vector $\aa$, we write $\min\{\aa\}$, $\med\{\aa\}$, and $\max\{\aa\}$ for its minimum, median, and maximum, respectively.}
A random variable $X$ is sub-Gaussian with variance proxy $\sigma^2$, denoted as $X \sim \text{subG}(\sigma^2)$, if $\b{E}[\exp\{s(X-\b{E}[X])\}] \leq \exp(s^2\sigma^2/2)$ for all $s\in\b{R}$. A random variable $X$ is sub-exponential
with parameter $\lambda$, denoted as $X \sim \text{subE}(\lambda)$, if $\b{E}[\exp\{s(X-\b{E}[X])\}] \leq \exp(s^2\lambda^2/2)$ for all $|s|\leq 1/\lambda$. Given a positive integer $a$, define $[a]:=\{1,\dots,a\}$. We define $\1_a$ as a vector of ones with length $a$, and $\M_a = \I_{a} - a^{-1}\1_{a}\1_{a}^\top$. 
Given a matrix $\A$, 
$\lambda_i(\bf{A})$ denotes its $i$-th largest eigenvalue,
$\bf{A}^\top$ denotes its transpose, 
and without loss of generality, the eigenvalues of a matrix are arranged by descending orders, and so are their corresponding eigenvectors.
We use $\diag(\{\A_1, \dots, \A_n\})$ to represent the block diagonal matrix with $\{\A_1, \dots, \A_n\}$ on the diagonal.
We denote by 
$|\cdot|$ the cardinality of a set,
$\|\cdot\|$ the spectral norm of a matrix or the $L_2$ norm of a vector, 
and $\|\cdot\|_F$ the Frobenius norm of a matrix.
Lastly, $\|\cdot\|_1$ and $\|\cdot\|_{\infty}$ denote the $L_1$ and $L_{\infty}$ norms of a matrix or a vector, 
respectively.

\section{Model and Estimation}\label{sec: model_estimation}

\subsection{Identification of main effect factor models}
\label{subsec: SMEFM}

For a matrix time series admitting the MEFM representation depicted in \eqref{eqn: MEFM}, we investigate thoroughly in this subsection the identification of main effects.
To start with, we say a function $g(\cdot): \b{R}^d\mapsto \b{R}$ is \textit{shift-varying} at $c_0$ if for any vector $\bf{a}\in \b{R}^d$ with $g(\bf{a}) =c_0$, we have $g(\bf{a} + c_1\1_d) \neq c_0$ for any constant $c_1\neq 0$.
We write in short $g(\cdot)$ is shift-varying if $g(\cdot)$ is shift-varying at any real number.
The following condition serves to identify parameters in model~\eqref{eqn: MEFM}.

\begin{itemize}
\item[(IC1)] (Identification).
{\em
We assume:
(i) $\A_r^\top \1_p=\0$ and $\A_c^\top \1_q=\0$;
(ii) Given functions $g_\alpha(\cdot): \b{R}^p\mapsto \b{R}$ and $g_\beta(\cdot): \b{R}^q\mapsto \b{R}$, we have $g_\alpha(\balpha_t^\ast) = 0$ and $g_\beta(\bbeta_t^\ast) = 0$ for any $t\in[T]$;
and (iii) The functions $g_\alpha(\cdot)$ and $g_\beta(\cdot)$ are shift-varying at zero.
}
\end{itemize}

Note immediately that for any constant $c$, the function $f(\cdot) = g(\cdot) + c$ is shift-varying at $c_0+c$ if and only if $g(\cdot)$ is shift-varying at $c_0$. Hence in parts~(ii)--(iii) of (IC1), it is equivalent to assume $g_\alpha(\balpha_t^\ast) = c_\alpha$ and $g_\beta(\bbeta_t^\ast) = c_\beta$, with $g_\alpha(\cdot)$ and $g_\beta(\cdot)$ being shift-varying at some constants $c_\alpha$ and $c_\beta$, respectively.
Hence Assumption~(IC1) is stated with $c_\alpha = c_\beta=0$, without loss of generality.
In fact, even the shift-varying functions and the corresponding constants can be time-varying at the cost of heavier notations, and we do not pursue this line of analysis.

Notably, the class of shift-varying functions encapsulates a large, useful set of functions, thus making the identified MEFM framework sufficiently general.
The functions $g_\alpha$ and $g_\beta$ are pre-specified by researchers to reflect any prior knowledge on the matrix time series. In what follows, we discuss two important types of shift-varying functions in terms of $g_\alpha$, together with examples of the practical interpretation of the identified row main effects by Assumption~(IC1).

\begin{enumerate}
    \item (Most linear functions).
    Given $\bf{x}\in \b{R}^p$, define the set of functions
    \begin{align}
        \cL:= \Big\{ g(\x) = \v^\top \x - c: \v \in \b{R}^p \text{\ with\ } \1_p^\top \v \neq 0, \; c\in \b{R} \Big\} .
        \label{eqn: set_linear_function}
    \end{align}
    Then all functions in $\cL$ are shift-varying.
    With non-negative entries in $\v$, the condition in Assumption~(IC1) can be viewed as imposing a level $c$ on the weighted average of the main effects.
    The identification strategy in \cite{lam2024matrix} is a special case of this with $\v=\1_p$ and $c=0$.
    Note that an extreme case can be $g_\alpha(\balpha_t^\ast) = \alpha_{t,i}^\ast$, that is, the identification of the row main effects is anchored at the $i$-th cross-sectional unit. This can be of particular interest when researchers wish to fix a reference unit and measure all other units relative to this baseline,
    and the main effects can be interpreted as deviations from the chosen benchmark.
    
    \item (All quantile functions). 
    With $q_u(\x)$ representing the $u$-th quantile of elements in a vector $\x$,
    define the set of functions
    \begin{align}
        \cQ:= \Big\{ g(\x) = 
        q_u(\x) - c:  
        u\in[0,1] , c\in \b{R} \Big\} .
        \label{eqn: set_quantile_function}
    \end{align}
    Then all functions in $\cQ$ are shift-varying.
    For $u=0$, $0.5$, and $1$, the identification on the row main effects can be read as $\min\{\balpha_t^\ast\}=c$, $\med\{\balpha_{t}^\ast\} = c$, and $\max\{\balpha_t^\ast\}=c$, respectively.
    The two extrema setups naturally accommodate the empirical range of the observed time series, such as non-negative values for metric reading in monitoring river flows in ecological application.
    On the other hand, main effects identified based on median function could serve as a tail-robust framework allowing for outliers, which are typical for financial return data.
\end{enumerate}

From the above discussion, our model design can be significantly flexible compared to MEFM in \cite{lam2024matrix}. While the above examples hold similarly for the column main effects, it is worth pointing out that the identification design for the column main effects can be different from that for the row main effects, depending upon the data sets in practice and the purpose of analysis.
Moreover, regardless of how the main effects are identified in (IC1), they always imply that, e.g.\ from $\balpha_t^\ast \1_q^\top$, the contribution of $\balpha_t^\ast$ towards $\X_t$ only varies over different rows and remains the same over columns. Hence the main effects specialize in row-wise and column-wise contribution, while the common component $\C_t$ picks up the interaction between rows and columns.
It is worth pointing out that the grand mean $\mu_t$ should be understood as the remainder from row and column main effects under the specified identification.

In all previous works on MEFM, the grand mean $\mu_t$ is first identified via the identification condition \eqref{eqn: MEFM_identification}, so that the row and column main effects can be subsequently identified. The same argument is inadmissible under our Assumption~(IC1), so we take another route by identifying the main effects first. In fact, by our new identification strategy, we can show that part~(iii) of Assumption~(IC1) is a necessary and sufficient condition to identify the main effects. Additionally, for a certain class of functions specified in Assumption~(IC1), model~\eqref{eqn: MEFM} is a more general framework than the Tucker-decomposition matrix factor model. We formalize all such results in the following theorem.

\begin{theorem}\label{thm: identification}
For model~\eqref{eqn: MEFM}, we have the following.
\begin{itemize}
    \item [(i)] Under Assumption~(IC1)~(i)--(ii), the parameters $\mu_t$, $\balpha_t^\ast$, $\bbeta_t^\ast$ and $\C_t$ can be identified if and only if (IC1)~(iii) holds.
    \item [(ii)] Let $g_\alpha(\cdot)$ and $g_\beta(\cdot)$ be the functions from Assumption~(IC1).
    Further assume $g_\alpha(\cdot)$ satisfies that for any $\bf{a}\in \b{R}^p$, there exists some scalar $c$ such that $g_\alpha(\bf{a}-c\1_p)= 0$; and assume similarly for $g_\beta(\cdot)$.
    If $\X_t$ follows a Tucker-decomposition matrix factor model such that $\X_t =\Acute\A_r \F_t\Acute\A_c^\top +\E_t$
    with arbitrary factor loading matrices $\Acute\A_r$ and $\Acute\A_c$,
    then $\X_t$ also follows model \eqref{eqn: MEFM} with the resulting parameters satisfying Assumption~(IC1).
\end{itemize}
\end{theorem}

Theorem~\ref{thm: identification}~(i) shows that
shift-varying functions fully characterize the identifiability of the particular structure in MEFM. 
In particular, the grand mean is identified on top of the main effects, which is different from \cite{lam2024matrix} where the grand mean can be directly identified.
Part~(ii) of the above theorem generalizes their Remark~1. 
We emphasize that the additional requirement on functions $g_\alpha$ and $g_\beta$ in Theorem~\ref{thm: identification}~(ii) is very mild, e.g.\ it can be fulfilled by all the functions in \eqref{eqn: set_linear_function} and \eqref{eqn: set_quantile_function}.

\subsection{Parameter estimation}\label{subsec: estimation}
In this subsection, we discuss the procedure of estimating parameters in model~\eqref{eqn: MEFM}. 
First, we state another condition to facilitate estimation.

\begin{itemize}
\item[(IC2)] (Closed-form solution).
{\em
Let $g_\alpha(\cdot)$ and $g_\beta(\cdot)$ be the shift-varying functions from Assumption~(IC1).
For any $\bf{x}\in \b{R}^p$, the equation $g_\alpha(\bf{x}-c\1_p)= 0$ has a closed-form solution $c = \wt{g}_\alpha(\bf{x})$ for some known function $\wt{g}_\alpha(\cdot)$; and assume similarly for $g_\beta(\cdot)$, with the solution function denoted as $\wt{g}_\beta(\cdot)$.
}
\end{itemize}

Assumption~(IC2) requires the existence of solution for the equation in Theorem~\ref{thm: identification}~(ii), and that the solution can be computed explicitly.
While not required, note that the solution must be unique by the definition of shift-varying functions and Assumption~(IC1)~(iii).
It is also straightforward to check that all functions in $\cL$ and $\cQ$ fulfill Assumption~(IC2);
we present their solution functions at the end of this section.

\paragraph{Estimation of grand mean and main effects.}
Fix any $t\in[T]$.
We first detail the estimation procedure for row main effects.
Right-multiplying $\1_q$ on $\X_t$ in \eqref{eqn: MEFM}, by Assumption~(IC1)~(i) we have
\begin{align*}
    \X_t \1_q = \1_p (q\mu_t +\1_q^\top \bbeta_t^\ast) + q\balpha_t^\ast + \E_t \1_q ,
\end{align*}
which, together with Assumption~(IC1)~(ii), gives
\begin{align*}
    0 = g_\alpha\left( \balpha_t^\ast \right)
    = g_\alpha\left( q^{-1} \X_t \1_q 
    -q^{-1} \E_t \1_q 
    - \1_p (\mu_t + q^{-1} \1_q^\top \bbeta_t^\ast) \right) .
\end{align*}
Then by Assumption~(IC2), we have $(\mu_t + q^{-1} \1_q^\top \bbeta_t^\ast) = \wt{g}_\alpha( q^{-1} \X_t \1_q  -q^{-1} \E_t \1_q )$, and hence
\begin{align*}
    \balpha_t^\ast 
    &=
    q^{-1} \X_t \1_q - \1_p \cdot \wt{g}_\alpha( q^{-1} \X_t \1_q  -q^{-1} \E_t \1_q ) - q^{-1} \E_t \1_q .
\end{align*}
Therefore, a feasible estimator for the row main effect at $t$ can be
\begin{align}
    \wt\balpha_t &:= 
    q^{-1} \X_t \1_q - \1_p \cdot \wt{g}_\alpha( q^{-1} \X_t \1_q) .
    \label{eqn: general_tilde_balpha_t}
\end{align}
Analogous arguments yield the estimator for the column main effect
\begin{align}
    \wt\bbeta_t &:= 
    p^{-1} \X_t^\top \1_p - \1_q \cdot \wt{g}_\beta( p^{-1} \X_t^\top \1_p) .
    \label{eqn: general_tilde_bbeta_t}
\end{align}
Next, we estimate the grand mean at time $t\in[T]$. This is straightforward using the row and column main effect estimators.
To this end, by left-multiplying $\1_p^\top$ and right-multiplying $\1_q$ on $\X_t$, we have
\[
\1_p^\top \X_t \1_q = pq\mu_t + q\1_p^\top \balpha_t^\ast + p\1_q^\top \bbeta_t^\ast + \1_p^\top \E_t \1_q .
\]
This leads to the grand mean estimator
\begin{align}
    \wt\mu_t &:= (pq)^{-1}\1_p^\top \X_t \1_q - p^{-1}\1_p^\top \wt\balpha_t - q^{-1}\1_q^\top \wt\bbeta_t.
    \label{eqn: general_tilde_mu_t}
\end{align}

\paragraph{Estimation of factor structure.}
Note that due to rotational indeterminacy, we cannot identify the loading matrices $\A_r$ and $\A_c$ exactly, but only their column spaces. To take into account potentially heterogeneous weak factors \citep[e.g.][]{LamYao2012, Zhangetal2024}, we normalize the loadings and the core factor as $\Q_r= \A_r\Z_r^{-1/2}$, $\Q_c= \A_c\Z_c^{-1/2}$, and $\F_{Z,t}= \Z_r^{1/2} \F_t \Z_c^{1/2}$, with $\Z_r$ and $\Z_c$ from Assumption~(L1); see more details in Section~\ref{subsec: assumption}. Since $\C_t=\Q_r \F_{Z,t} \Q_c^\top$, we may equivalently estimate the normalized parameters. To this end, define the matrix
\begin{align}
    \wt\L_t &:= \X_t - \wt\mu_t \1_p \1_q^\top - \wt\balpha_t \1_q^\top - \1_p \wt\bbeta_t^\top.
    \label{eqn: tilde_L_t}
\end{align}
Then the estimator for the normalized row loading matrix, denoted by $\wh\Q_r$, is defined as the eigenvector matrix corresponding to the $k_r$ largest eigenvalues of the matrix $T^{-1} \sum_{t=1}^T \wt\L_t \wt\L_t^\top$. Similarly, the normalized column loading matrix estimator $\wh\Q_c$ is the eigenvector matrix corresponding to the $k_c$ largest eigenvalues of $T^{-1} \sum_{t=1}^T \wt\L_t^\top \wt\L_t$. 
Lastly, the core factor and the common component can be estimated by
\begin{align*}
    \wh\F_{Z,t} &:= \wh\Q_r^\top \wt\L_t \wh\Q_c, \quad
    \wh\C_t := \wh\Q_r \wh\F_{Z,t} \wh\Q_c^\top 
    = \wh\Q_r \wh\Q_r^\top \X_t \wh\Q_c\wh\Q_c^\top ,
\end{align*}
where the last equality uses the fact that $\wh\Q_r^\top \wt\L_t \wh\Q_c = \wh\Q_r^\top \X_t \wh\Q_c$, since both $\wh\Q_r$ and $\wh\Q_c$ fulfill Assumption~(IC1)~(i) with $\A_r^\top \1_p=\0$ and $\A_c^\top \1_q=\0$ (see Lemma~\ref{lemma: loading_est_fulfill_iden} in the supplement).

\paragraph{Main effect estimators under identification of $\cL$ and $\cQ$, respectively.}
Lastly, we give explicit forms of the main effect estimators in \eqref{eqn: general_tilde_balpha_t} and \eqref{eqn: general_tilde_bbeta_t}, when the shift-varying functions in Assumption~(IC1) are specified from $\cL$ and $\cQ$, respectively.
To this end, first note that for any shift-varying functions in $\cL$, using the notations in \eqref{eqn: set_linear_function}, their solution functions and thus the estimator in \eqref{eqn: general_tilde_balpha_t} take the form
\begin{equation}
\label{eqn: linear_tilde_balpha_t}
\begin{split}
    \wt{g}(\x) &= (\v^\top \1_p)^{-1} \left(\v^\top \x - c \right) , \\
    \wt\balpha_t &= 
    q^{-1} \X_t \1_q - \1_p \cdot (\v^\top \1_p)^{-1} \left( q^{-1} \v^\top \X_t \1_q - c \right) .
\end{split}
\end{equation}
For instance, when $\v=\1_p$ and $c=0$, the above main effect estimator simplifies to 
\begin{align*}
    \wt\balpha_t = q^{-1} \X_t \1_q - (pq)^{-1} \1_p \1_p^\top \X_t \1_q 
    = 
    q^{-1} (\I_p - p^{-1} \1_p \1_p^\top) \X_t \1_q 
    = q^{-1} \M_p \X_t \1_q ,
\end{align*}
which matches Equations~(3.2)--(3.3) in \cite{lam2024matrix}.

On the other hand, consider shift-varying functions from $\cQ$. With the notations in \eqref{eqn: set_quantile_function}, we can write the solution function and the corresponding estimator respectively as
\begin{equation}
\label{eqn: quantile_tilde_balpha_t}
\begin{split}
    \wt{g}(\x) &= q_u(\x) -c , \\
    \wt\balpha_t &= 
    q^{-1} \X_t \1_q - q^{-1} \1_p \cdot    q_u(\X_t \1_q) + c \1_p  .
\end{split}
\end{equation}
Suppose $c=0$.
When $u=0.5$, i.e., $q_u(\cdot)$ takes the median of a vector, the main effect estimator reads as $\wt\balpha_t = q^{-1} \X_t \1_q - q^{-1} \1_p \cdot \med\{\X_t \1_q\}$.
Similarly, if we set $u=0$, then the main effect estimator is $\wt\balpha_t = q^{-1} \X_t \1_q - q^{-1} \1_p \cdot \min\{\X_t \1_q\}$.
Markedly, given $c=0$, model identification based on any shift-varying functions from $\cQ$ indicates that the main effects vanish at a certain quantile. This allows practitioners to study the sparsity pattern on main effects, which can be of particular interest; we defer the detailed discussion in Section~\ref{sec: sparsity_ME}.

\section{Assumptions and Theoretical Results}\label{sec: assumption_theories}

\subsection{Assumptions}\label{subsec: assumption}

We discuss the main theoretical results under general identification based on Assumptions~(IC1) and (IC2).
For simplicity, throughout this paper, we assume the shift-varying functions for the main effects
have the same functional form which is generically denoted by $g(\cdot)$.
To begin with, we present additional assumptions in the following to characterize model~\eqref{eqn: MEFM}. The explanations to each assumption are deferred to the end of this subsection.

\begin{itemize}
\item[(L1)] (Factor strength).
{\em
We assume that $\A_r$ and $\A_c$ are of full rank and independent of $\{\F_t\}$ and $\{\E_t\}$. Furthermore, as $p,q\to\infty$,
\begin{equation}
\label{eqn: L1}
    \Z_r^{-1/2}\A_r^\top \A_r\Z_r^{-1/2}
    \to \bSigma_{A,r}, \quad
    \Z_c^{-1/2}\A_c^\top \A_c\Z_c^{-1/2}
    \to \bSigma_{A,c},
\end{equation}
where $\Z_r=\textnormal{diag}(\A_r^\top \A_r)$, $\Z_c=\textnormal{diag}(\A_c^\top \A_c)$, and both $\bSigma_{A,r}$ and $\bSigma_{A,c}$ are positive definite with all eigenvalues bounded away from 0 and infinity. We assume $(\Z_r)_{jj}\asymp p^{\delta_{r,j}}$ for $j\in[k_r]$ and $1/2<\delta_{r,k_r}\leq \dots\leq \delta_{r,2}\leq \delta_{r,1}\leq 1$. Similarly, we assume $(\Z_c)_{jj}\asymp q^{\delta_{c,j}}$ for $j\in[k_c]$, with $1/2< \delta_{c,k_c}\leq \dots\leq \delta_{c,2}\leq \delta_{c,1}\leq 1$.
}
\end{itemize}

\begin{itemize}
    \item[(F1)] (Time series in $\F_t$).
{\em
There is $\X_{f,t}$ the same dimension as $\F_t$, such that $\F_t = \sum_{w\geq 0}a_{f,w}\X_{f,t-w}$. The time series $\{\X_{f,t}\}$ has \text{i.i.d.} elements with mean $0$ and variance $1$, with uniformly bounded fourth order moments. The coefficients $a_{f,w}$ are such that $\sum_{w\geq 0}a_{f,w}^2=1$ and $\sum_{w\geq 0}|a_{f,w}|\leq c$ for some constant $c$.
}
\end{itemize}

\begin{itemize}
  \item[(E1)] (Decomposition of $\E_t$).
{\em We assume that
\begin{equation}
    \label{eqn: E1}
    \E_t = \A_{e,r}\F_{e,t}\A_{e,c}^\top + \bSigma_{\epsilon}\circ \bepsilon_t ,
\end{equation}
where $\F_{e,t}$ is a matrix of size $k_{e,r}\times k_{e,c}$, containing independent elements with mean $0$ and variance $1$. The matrix $\bepsilon_t\in\b{R}^{p\times q}$ contains independent elements with mean 0 and variance 1, with $\{\bepsilon_t\}$ independent of $\{\F_{e,t}\}$. The matrix $\bSigma_{\epsilon}$ contains the standard deviations of the corresponding elements in $\bepsilon_t$, and has elements uniformly bounded away from 0 and infinity.
Moreover, $\A_{e,r}$ and $\A_{e,c}$ are (approximately) sparse matrices with sizes $p\times k_{e,r}$ and $q\times k_{e,c}$ respectively, such that $\|\A_{e,r}\|_1, \|\A_{e,c}\|_1=O(1)$, with $k_{e,r}, k_{e,c} = O(1)$.
}
\end{itemize}
\begin{itemize}
  \item[(E2)] (Time series in $\E_t$).
{\em
There is $\X_{e,t}$ the same dimension as $\F_{e,t}$, and $\X_{\epsilon,t}$ the same dimension as $\bepsilon_t$, such that $\F_{e,t} = \sum_{w\geq 0}a_{e,w}\X_{e,t-w}$ and $\bepsilon_t = \sum_{w\geq 0}a_{\epsilon,w}\X_{\epsilon,t-w}$, with $\{\X_{e,t}\}$ and $\{\X_{\epsilon,t}\}$ independent of each other, and each time series has independent elements with mean $0$ and variance $1$ with uniformly bounded fourth order moments. Both $\{\X_{e,t}\}$ and $\{\X_{\epsilon,t}\}$ are independent of $\{\X_{f,t}\}$ from Assumption~(F1).
The coefficients $a_{e,w}$ and $a_{\epsilon,w}$ satisfy $\sum_{w\geq 0}a_{e,w}^2 = \sum_{w\geq 0}a_{\epsilon,w}^2 = 1$ and $\sum_{w\geq 0}|a_{e,w}|, \sum_{w\geq 0}|a_{\epsilon,w}|\leq c$ for some constant $c$.
}
\end{itemize}

\begin{itemize}
  \item[(E3)] (Tail condition in $\F_t$ and $\E_t$).
{\em
Each element in the time series $\{\X_{f,t}\}$ from Assumption (F1), $\{\X_{e,t}\}$ and $\{\X_{\epsilon,t}\}$ from Assumption~(E2) is sub-Gaussian.
}
\end{itemize}

\begin{itemize}
  \item[(R1)] (Rate assumptions).
{\em
We assume that,
\begin{align*}
    & T^{-1} p^{2(1-\delta_{r,k_r})}q^{1-2\delta_{c,1}} = o(1) , \quad
    p^{1-2\delta_{r,k_r}}q^{2(1-\delta_{c,1})} = o(1), \\
    & T^{-1} q^{2(1-\delta_{c,k_c})}p^{1-2\delta_{r,1}} = o(1) , \quad
    q^{1-2\delta_{c,k_c}}p^{2(1-\delta_{r,1})} = o(1) .
\end{align*}
}
\end{itemize}

\begin{itemize}
  \item[(S1)] (Lipschitz continuity).
{\em
The solution functions $\wt{g}_\alpha(\cdot)$ and $\wt{g}_\beta(\cdot)$ in Assumption~(IC2) are Lipschitz continuous with respect to some vector norm $\|\cdot\|_\ast$, such that there exist constants $L_\alpha$ and $L_\beta$, that may depend on $T,p,q$, such that for all $\x_1, \x_2\in \b{R}^p$ and $\y_1, \y_2 \in \b{R}^q$:
\begin{align*}
    |\wt{g}_\alpha( \x_1 ) - \wt{g}_\alpha( \x_2 ) |  \leq L_\alpha \|\x_1 - \x_2\|_\ast ,
    \quad
    |\wt{g}_\beta( \y_1 ) - \wt{g}_\beta( \y_2 ) |  \leq L_\beta \|\y_1 - \y_2\|_\ast .
\end{align*}
}
\end{itemize}

With Assumption~(L1), we can define the normalized row and column loading matrices $\Q_r=\A_r\Z_r^{-1/2}$ and $\Q_c=\A_c\Z_c^{-1/2}$. Hence $\Q_r^\top \Q_r\to\bSigma_{A,r}$ and $\Q_c^\top \Q_c\to\bSigma_{A,c}$. Assumption~(L1) allows for weak factors with heterogeneous factor strengths, which is also similarly seen in Remark~1 in \cite{LamYao2012}, differing from traditional approximate factor models that consider either pervasive factors only \citep{Heetal2024} or weak factors with the same factor strength \citep{Wangetal2019}.

Assumptions~(F1), (E1) and (E2) characterize the dynamics in the core factor and noise series, which are necessary in matrix factor models. Rather than directly presenting, \text{e.g.} the weak dependence in noise as in Assumptions~(D) in \cite{Yuetal2022}, our assumptions naturally specify the data generating process, allowing the noise and core factors to be general linear processes. In particular, the core factor can be serially correlated under (F1) which also implies that for each $t\in[T]$, $\b{E}[\F_t \F_t^\top] = k_c \I_{k_r}$ and $\b{E}[\F_t^\top \F_t] = k_r \I_{k_c}$. This together with Assumption~(L1) thus serve as an alternative set of identification conditions where the dependence among the latent dynamics driven by the core factors is featured by the loading matrices; see the discussion in Section~3 in \cite{BaiNg2002} for instance. The decomposition by Assumption~(E1) together with the general linear process in (E2) allow the noise to have both serial and cross-sectional dependence. Hence our Assumptions~(E1) and (E2) are comparable to, for example, conditional independence as in Assumption~3.2 in \cite{Heetal2024}.

Assumption~(E3) controls the tail of the random variables in model \eqref{eqn: MEFM}. It is arguably very mild since stronger assumptions are often needed in the regularized regression literature, such as Condition~(a) in \cite{Zou2006}, Assumption~(A1) in \cite{huang2008adaptive}, and Assumption~(E) in \cite{Rinaldo2009}, to name but a few. In particular, as pointed out in Remark~3 in \cite{fan2024latent}, such sub-Gaussianity would not hold under highly correlated covariates in the regression problem, but the covariate matrix in our formulation is an identity matrix (see Section~\ref{subsec: practical_implementation}) and hence this issue is circumvented.

Assumption~(R1) spells out the rates required on the dimensions and factor strengths, which directly hold if all factors are pervasive.
Assumption~(S1) imposes a smoothness condition on the solution functions, so as to study the convergence rate of the main effect estimators.
Markedly, using different vector norms controls the deviation of the grand mean and main effect estimators from the true parameters to various extents, since the rates of convergence depend on some forms of the noise measured in $\|\cdot\|_\ast$.
Note that the vector norms in Assumption~(S1) can be specified differently over row and column main effects and timestamps, if the shift-varying functions in Assumption~(IC1) are set differently.  
Throughout this paper, we assume all these are identical given the observed time series for simplicity.
To further elaborate this assumption, we discuss some examples using the row main effects.

Suppose $\|\cdot\|_\ast$ is set as the Euclidean norm $\|\cdot\|$.
Then for any shift-varying function in the set $\cL$, its solution function takes the form \eqref{eqn: linear_tilde_balpha_t} and hence
\[
|\wt{g}_\alpha(\x_1) - \wt{g}_\alpha(\x_2)|
= \left| (\v^\top\1_p)^{-1} \v^\top(\x_1-\x_2) \right|
\leq 
\frac{\| \v\|}{|\v^\top\1_p|} \cdot \|\x_1-\x_2\| ,
\]
which indicates its Lipschitz constant as $L_\alpha = \|\v\| /|\v^\top\1_p|$. For instance, $L_\alpha = p^{-1/2}$ when $\v=\1_p$.
The solution function of any shift-varying function in $\cQ$ admits the form in \eqref{eqn: quantile_tilde_balpha_t}, which has the corresponding Lipschitz constant $L_\alpha= 1$, since
\begin{align*}
    |\wt{g}_\alpha(\x_1) - \wt{g}_\alpha(\x_2)|
    = |q_u(\x_1) - q_u(\x_2)|
    \leq 
    \|\x_1-\x_2\|_{\infty} 
    \leq
    \|\x_1-\x_2\| .
\end{align*}
On the other hand, suppose $\|\cdot\|_{\infty}$ is the $L_{\infty}$ norm.
Then for functions in $\cL$, similar steps as above give 
$|\wt{g}_\alpha(\x_1) - \wt{g}_\alpha(\x_2)| \leq (\| \v\|_1/ |\v^\top\1_p|) \cdot \|\x_1-\x_2\|_{\infty}$,
thus implying $L_\alpha = \|\v\|_1 /|\v^\top\1_p|$.
For functions in $\cQ$, their Lipschitz constants are again $L_\alpha = 1$ from the above display.

\subsection{Theoretical results}\label{subsec: theorem}

We now study the statistical properties for the estimators in Section~\ref{subsec: estimation}.
First, we present the theoretical guarantee on estimating the main effects and grand mean, under identification using arbitrary shift-varying functions in Assumption~(IC1).

\begin{theorem}[Properties of $\wt\balpha_t$, $\wt\bbeta_t$, and $\wt\mu_t$, under arbitrary shift-varying functions]
\label{thm: general_initial_consistency_maineffect}
Let Assumptions (IC1), (IC2), (E1), (E2), and (S1) hold.
Then the estimators for the main effects and the grand mean in Section~\ref{subsec: estimation} satisfy that for any $t\in[T]$, as $\min(T,p,q) \to \infty$,
    \begin{align*}
        \frac{1}{p}\|\wt\balpha_t - \balpha_t^\ast\|^2 &=
        O_P\bigg( \frac{L_\alpha^2}{q^2} \left\| \E_t \1_q \right\|_\ast^2 
        + \frac{1}{q} \bigg) , \\
        \frac{1}{q}\|\wt\bbeta_t - \bbeta_t^\ast\|^2 &=
        O_P\bigg( \frac{L_\beta^2}{p^2} \left\| \E_t^\top \1_p \right\|_\ast^2 
        + \frac{1}{p} \bigg) , \\
        (\wt\mu_t - \mu_t)^2 &= 
        O_P\bigg( \frac{L_\alpha^2}{q^2} \left\| \E_t \1_q \right\|_\ast^2 
        + \frac{L_\beta^2}{p^2} \left\| \E_t^\top \1_p \right\|_\ast^2 
        + \frac{1}{pq} \bigg) .
    \end{align*}
\end{theorem}

Theorem~\ref{thm: general_initial_consistency_maineffect} spells out the rates of convergence for the estimators defined in \eqref{eqn: general_tilde_balpha_t}, \eqref{eqn: general_tilde_bbeta_t}, and \eqref{eqn: general_tilde_mu_t}. The theorem indicates that, for instance, the row (resp.\ column) main effect estimators at least carry the rate of $1/q$ (resp.\ $1/p$).
Also notably, the convergence rate for the grand mean estimator has some identical terms as in the main effect estimators. While it is unsurprising from its construction in \eqref{eqn: general_tilde_mu_t},
the presented result only serves as a feasible bound for general identification in Assumption~(IC1).
The following theorem shows an exception of $\v=\1$, where the grand mean estimator has better rate than the main effect estimators;
it also gives explicit results when the shift-varying functions in Assumption~(IC1) are specified in $\cL$ or $\cQ$.

\begin{theorem}[Properties of $\wt\balpha_t$, $\wt\bbeta_t$, and $\wt\mu_t$, under shift-varying functions in $\cL$ or $\cQ$]
\label{thm: particular_initial_consistency_maineffect}
In addition to all assumptions in Theorem~\ref{thm: general_initial_consistency_maineffect}, further let Assumption~(E3) hold.
Then for any $t\in[T]$, as $\min(T,p,q) \to \infty$, the following hold.
\begin{itemize}
    \item [(i)] When the shift-varying functions are in $\cL$, i.e., $g(\x) = \v^\top \x - c$, we have
    \begin{align*}
        \frac{1}{p}\|\wt\balpha_t - \balpha_t^\ast\|^2 &=
        O_P\Bigg\{ \frac{ \big( p \|\v\|^2 \big) 
        \bigwedge \big\{ \log(p) \|\v\|_1^2 \big\} }{q (\v^\top\1_p)^2}
        + \frac{1}{q} \Bigg\}  , \\
        \frac{1}{q}\|\wt\bbeta_t - \bbeta_t^\ast\|^2 &=
        O_P\Bigg\{ \frac{ \big( q \|\v\|^2 \big) 
        \bigwedge \big\{ \log(q) \|\v\|_1^2 \big\} }{p (\v^\top\1_q)^2}
        + \frac{1}{p} \Bigg\}  , \\
        (\wt\mu_t - \mu_t)^2 &= 
        O_P\Bigg\{ 
        \frac{ \big( p \|\v\|^2 \big) 
        \bigwedge \big\{ \log(p) \|\v\|_1^2 \big\} }{q (\v^\top\1_p)^2} 
        + \frac{ \big( q \|\v\|^2 \big) 
        \bigwedge \big\{ \log(q) \|\v\|_1^2 \big\} }{p (\v^\top\1_q)^2}
        + \frac{1}{pq} \Bigg\} .
    \end{align*}
    In particular, when $\v=\1$ (with dimensions $p$ and $q$ for the row and column main effects, respectively), we have $p^{-1} \|\wt\balpha_t - \balpha_t^\ast\|^2 = O_P(1/q)$, 
    $q^{-1} \|\wt\bbeta_t - \bbeta_t^\ast\|^2 = O_P(1/p)$,
    and markedly, $(\wt\mu_t - \mu_t)^2 = O_P\{ 1/(pq) \}$.

    \item [(ii)] When the shift-varying functions are in $\cQ$, i.e., $g(\x) = q_u(\x) - c$, we have
    \begin{align*}
        \frac{1}{p}\|\wt\balpha_t - \balpha_t^\ast\|^2 &=
        O_P\bigg\{ \frac{\log(p)}{q} \bigg\} , \\
        \frac{1}{q}\|\wt\bbeta_t - \bbeta_t^\ast\|^2 &=
        O_P\bigg\{ \frac{\log(q)}{p} \bigg\} , \\
        (\wt\mu_t - \mu_t)^2 &= 
        O_P\bigg\{ \frac{\log(p)}{q} + \frac{\log(q)}{p} \bigg\} .
    \end{align*}
\end{itemize}
\end{theorem}

Notably, both parts of the theorem are agnostic to the value of $c$.
From Theorem~\ref{thm: particular_initial_consistency_maineffect}~(i), 
the choice of $\v$ guides the convergence rates. For instance, if $\v$ is sparse, then $p \|\v\|^2$ tends to dominate $\log(p) \|\v\|_1^2$ and the latter would then be the resulting term in the final presented rate. 
One such instance is when $\v$ is some standard basis vector, in which case we would have $p^{-1} \|\wt\balpha_t - \balpha_t^\ast\|^2 = O_P\{\log(p)/q\}$, etc. 
Setting $\v=\1$ illustrates a particular case where the grand mean estimator has better rate of convergence than the main effect estimators. This is essentially due to the cancellation of main effect estimators when contributed towards the estimation error of $\wt\mu_t$.
We mark that part~(i) also implies any non-zero scaling of $\v$ would not change the rates of convergence.

Theorem~\ref{thm: particular_initial_consistency_maineffect}~(ii) is a consequence of Theorem~\ref{thm: general_initial_consistency_maineffect} with $L_{\infty}$ norm in Assumption~(S1).
The resulting rates carry some logarithmic terms of the dimensions,
which is intuitive since the extrema of the noise series are inevitably involved due to the worst-case identification with $u$ close to 0 or 1.
However, when the shift-varying functions are central quantiles (e.g., median with $u=0.5$),
the Lipschitz constant in Assumption~(S1) can be of smaller order if the main effects are sparse and the noise distribution has a smooth density, thus leading to an improved rate;
we do not pursue this direction here.

Next, define the following square matrices
\begin{equation*}
\begin{split}
        \H_r &:= \frac{1}{T} \wh\D_r^{-1} \wh\Q_r^\top \Q_r \sum_{t=1}^T (\F_{Z,t} \Q_c^\top \Q_c \F_{Z,t}^\top) , \\
        \H_c &:= \frac{1}{T} \wh\D_c^{-1} \wh\Q_c^\top \Q_c \sum_{t=1}^T (\F_{Z,t}^\top \Q_r^\top \Q_r \F_{Z,t}) ,
\end{split}
\end{equation*}
where $\wh\D_r := \wh\Q_r^\top(T^{-1} \sum_{t=1}^T \wt\L_t\wt\L_t^\top)\wh\Q_r$ is the $k_r \times k_r$ diagonal matrix consisting of eigenvalues of $T^{-1} \sum_{t=1}^T \wt\L_t\wt\L_t^\top$, and $\wh\D_c := \wh\Q_c^\top (T^{-1} \sum_{t=1}^T \wt\L_t^\top\wt\L_t)\wh\Q_c$ is the $k_c \times k_c$ diagonal matrix of eigenvalues of $T^{-1} \sum_{t=1}^T \wt\L_t^\top\wt\L_t$.
Our next theorem then presents the asymptotic properties of $\H_r$, $\H_c$, and our estimators on the factor structure in Section~\ref{subsec: estimation}.

\begin{theorem}[Asymptotic rates for estimators on factor structure]
\label{thm: general_initial_consistency}
Let Assumptions (IC1), (IC2), (L1), (F1), (E1), (E2), and (R1) hold.
Then as $\min(T,p,q) \to \infty$,
$\H_r$ and $\H_c$ are asymptotically invertible,
and both the row and column factor loading estimators are consistent such that
\begin{align*}
    \frac{1}{p} \big\| \wh\Q_r-\Q_r\H_r^\top \big\|_F^2 &= O_P\Big( T^{-1}p^{1-2\delta_{r,k_r}}q^{1-2\delta_{c,1}} + p^{-3\delta_{r,k_r}} q^{2(1-\delta_{c,1})} \Big),  \\
    \frac{1}{q} \big\| \wh\Q_c-\Q_c\H_c^\top \big\|_F^2 &= O_P\Big( T^{-1}q^{1-2\delta_{c,k_c}}p^{1-2\delta_{r,1}} + q^{-3\delta_{c,k_c}} p^{2(1-\delta_{r,1})} \Big).
\end{align*}
Furthermore, the estimated core factor series and common components are consistent such that for any $t\in[T]$, $i\in[p]$, $j\in[q]$, we have
{\small
\begin{align*}
        &
        \big\| \wh\F_{Z,t} - (\H_r^{-1})^\top \F_{Z,t} \H_c^{-1} \big\|_F^2 
        =
        O_P\big( p^{1-\delta_{r,k_r}} q^{1-\delta_{c,k_c}}
        + T^{-1} p^{1+ 2\delta_{r,1} - 2\delta_{r,k_r}} q^{1-\delta_{c,1}} +  p^{1 +\delta_{r,1} -3\delta_{r,k_r}} q^{2-\delta_{c,1}} \\
        &\quad\quad\quad\quad\quad\quad\quad\quad\quad\quad\quad\quad\quad
        + T^{-1} q^{1+ 2\delta_{c,1} - 2\delta_{c,k_c}} p^{1-\delta_{r,1}} +  q^{1 +\delta_{c,1} -3\delta_{c,k_c}} p^{2-\delta_{r,1}}\big) ,\\
        &
        (\wh{C}_{t,ij} - C_{t,ij})^2 
        =
        O_P\big(p^{1-2\delta_{r,k_r}} q^{1-2\delta_{c,k_c}}
        + T^{-1} p^{1 + 2\delta_{r,1} -3\delta_{r,k_r}} q^{1 -\delta_{c,1} -\delta_{c,k_c}} + p^{1 +\delta_{r,1} -4\delta_{r,k_r}} q^{2 -\delta_{c,1} -\delta_{c,k_c}}  \\
        &\quad\quad\quad\quad\quad\quad\quad
        + T^{-1} q^{1 + 2\delta_{c,1} -3\delta_{c,k_c}} p^{1 -\delta_{r,1} -\delta_{r,k_r}} + q^{1 +\delta_{c,1} -4\delta_{c,k_c}} p^{2 -\delta_{r,1} -\delta_{r,k_r}}
        \big) .
\end{align*}
}
\end{theorem}

Theorem~\ref{thm: general_initial_consistency} details the behaviors of the estimators of factor loading matrices, core factors, and thus common components, under different factor strengths.
For a clear comparison with results in the literature,
Corollary~\ref{cor: consistency_strong_factor} (in the following) spells out the rates of convergence in Theorem~\ref{thm: general_initial_consistency} when all factors are strong.
From the corollary,
our factor loading estimators have comparable performance as the $\alpha$-PCA estimators in \cite{ChenFan2023}.
Also, our results align with, for instance, Theorem~4.1 of \cite{Heetal2024} and Theorem~4 of \cite{ChenFan2023}.

\begin{corollary}[Simplified Theorem~\ref{thm: general_initial_consistency} under pervasive factors]
\label{cor: consistency_strong_factor}
Let all assumptions in Theorem~\ref{thm: general_initial_consistency} hold, and further assume that $\delta_{r,i}= \delta_{c,j}=1$ for all $i\in[k_r]$, $j\in[k_c]$.
Define the re-normalized row and column loading estimators, and core factor estimator as
$\wh\A_r := \sqrt{p} \, \wh\Q_r$,
$\wh\A_c := \sqrt{q} \, \wh\Q_c$,
and $\wh\F_t := \wh\F_{Z,t} /\sqrt{pq}$,
respectively.
Then for any $t\in[T]$, $i\in[p]$, $j\in[q]$,
\begin{align*}
    \frac{1}{p} \big\| \wh\A_r -\A_r\H_r^\top \big\|_F^2 = O_P\bigg(\frac{1}{Tq} + \frac{1}{p^2} \bigg),  \quad
    \frac{1}{q} \big\| \wh\A_c -\A_c\H_c^\top \big\|_F^2 = O_P\bigg(\frac{1}{Tp} + \frac{1}{q^2} \bigg) &, \\
    \big\| \wh\F_{t} - (\H_r^{-1})^\top \F_{t} \H_c^{-1} \big\|_F^2 ,\quad 
    (\wh{C}_{t,ij} - C_{t,ij})^2 = O_P\bigg( \frac{1}{Tq} + \frac{1}{Tp} + \frac{1}{p^2} + \frac{1}{q^2}\bigg) &.
\end{align*}
\end{corollary}


\section{Sparsity Recovery by Regularized Estimation}\label{sec: sparsity_ME}

As a useful line of analysis, we study the sparsity of main effects in this section.
Throughout, the shift-varying function in Assumption~(IC1) takes the form in \eqref{eqn: set_quantile_function} with $u, c=0$,
i.e., $g(\cdot) = \min\{\cdot\}$ so that the main effects are identified by 
\begin{align}
    \min\{\balpha_t^\ast \} = 0 ,
    \quad
    \min\{\bbeta_t^\ast \} = 0 .
    \label{eqn: IC1_min0}
\end{align}
For each cross-sectional unit, the main effects can be sparse in certain periods.
Formally, consider the row effects $\{\balpha_t^\ast \}_{t\in[T]}$, and for any $i\in[p]$, we define the sparse block as $\cS_{\alpha,i} =\{t: \alpha_{t,i}^\ast=0\}$, and the dense block $\cB_{\alpha,i} = [T] \setminus \cS_{\alpha,i}$. The sparse and dense blocks for the column main effects are defined similarly, denoted as $\cS_{\beta,j}$ and $\cB_{\beta,j}$ respectively. Note that both $\cS_{\alpha,i}$ and $\cS_{\beta,j}$ can be potentially empty.
From a data generating point of view, the sparse blocks should be viewed as non-random sets of timestamps where the corresponding main effects vanish, and the dense blocks are the remaining periods.

For illustration, consider again the example of economic indices from Section~\ref{sec: introduction}, where each row corresponds to one country and hence the time series $\{\alpha_{t,i}^\ast\}$ represents the country-$i$'s effect which could disappear when the country's economy goes down for instance. To strengthen the idea that the sparsity can be piecewise, i.e., zeros in the main effects can be consecutive in timestamps, we address this by incorporating a total variation loss, which we defer to Section~\ref{subsec: estimation_regularize}. Note that there could be singular zeros which live in between non-zero main effects, so that our framework balances between interpretability and generality. To summarize, the benefit of such a sparsity framework is at least two-fold:
\begin{enumerate}
    \item With piecewise zeros, the main effects can be understood more easily by practitioners; on the other hand, singular zeros potentially imply influential events and merit further investigation.
    \item Allowing for general non-zeros retains adequate flexibility, compared to \text{e.g.} piecewise constants which is restricted and may not be realistic. Moreover, the dense blocks are able to feature the important patterns in the observed matrix time series, such as high volatility for financial return data.
\end{enumerate}

\subsection{Regularized estimation}
\label{subsec: estimation_regularize}

We address in this subsection the estimation problem under sparsity,
and our goal is to recover the sparse blocks for the row and column main effects.
By \eqref{eqn: quantile_tilde_balpha_t} and analogous steps in Section~\ref{subsec: estimation}, we first obtain the initial estimators for the main effects as
\begin{align}
    \wt\balpha_t &:= q^{-1}\X_t \1_q - q^{-1}\1_p \min\{ \X_t \1_q\} ,
    \label{eqn: tilde_balpha_t} \\
    \wt\bbeta_t &:= p^{-1}\X_t^\top \1_p - p^{-1}\1_q \min\{ \X_t^\top \1_p\} .
    \label{eqn: tilde_bbeta_t} 
\end{align}

Inspired by the Lasso \citep{Tibshirani1996}, we may employ an $L_1$ penalty to select the non-zero main effects. Nevertheless, such a regularization penalizes uniformly on each main effect value, potentially leading to over- or under-penalization which, in our scenario, fails block consistency (i.e.,  sparse blocks and dense blocks are estimated exactly as the true sets). To circumvent the restricted yet necessary irrepresentable condition in Lasso, \cite{Zou2006} proposes to adaptively penalize the magnitude of estimators by reweighing the $L_1$ loss by the inverse of some well-behaved initial estimators. This adaptive Lasso method enlightens us to consider a loss function such that, for $i\in[p]$,
\begin{align}
    L^\circ(\alpha_{1,i}, \dots, \alpha_{T,i}) := \frac{1}{2} \sum_{t=1}^T \big(\alpha_{t,i}^\ast - \alpha_{t,i}\big)^2 + \lambda_\alpha \sum_{t=1}^{T}\gamma_{\alpha, t, i}\big|\alpha_{t, i}\big| ,
    \label{eqn: loss_function_AL}
\end{align}
where $\gamma_{\alpha,t,i} =1/ \wt\alpha_{t,i}$ and $\lambda_\alpha$ is the tuning parameter. The estimators obtained by minimizing \eqref{eqn: loss_function_AL} are theoretically solid, but suffer from two flaws in practice. First, as a part of the latent MEFM representation, the row main effects cannot be directly observed, and hence all the $\alpha_{t,i}^\ast$'s in \eqref{eqn: loss_function_AL} are unavailable. Secondly, even if the true sparse block contains consecutive indices, the resulted main effect estimators from \eqref{eqn: loss_function_AL} do not favor piecewise sparsity, thus undermining the interpretation empirically.

To address the first concern above, we leverage the initial estimator $\wt\alpha_{t,i}$ which turns out to be an appropriate proxy to $\alpha_{t,i}^\ast$ under very general assumptions; see Assumption~(R2) in Section~\ref{subsec: thm_sparse_guarantee}. To promote smoothness in the estimator, we further include a total variation loss in the objective function, which is akin to the fused Lasso method \citep{Tibshiranietal2005, Rinaldo2009}. Different from the traditional use of the total variation loss, we borrow the idea from adaptive Lasso again and reweigh the penalty by $1/\max\{\wt\alpha_{t,i}, \wt\alpha_{t-1,i}\}$, so that the smoothness is mainly encouraged on the sparse blocks. To this end, we introduce a new regularized estimator called the \textit{Doubly Adaptive Fused Lasso (DAFL)} estimator, detailed as follows. For each $i\in[p]$, the DAFL estimator for the $i$-th row effect, $\{\wh\alpha_{t,i}\}_{t\in[T]}$, is obtained by minimizing the penalized loss
\begin{align}
    L(\alpha_{1,i}, \dots, \alpha_{T,i}) := \frac{1}{2} \sum_{t=1}^T \big(\wt\alpha_{t,i} - \alpha_{t,i}\big)^2 + \lambda_\alpha \sum_{t=2}^T u_{\alpha, t, i}\big|\alpha_{t, i} - \alpha_{t-1, i}\big| + \lambda_\alpha \sum_{t=1}^{T}\gamma_{\alpha, t, i}\big|\alpha_{t, i}\big| ,
    \label{eqn: DAFL}
\end{align}
where $u_{\alpha,t,i} = 1/\max\{\wt\alpha_{t,i}, \wt\alpha_{t-1,i}\}$, and $\lambda_\alpha$ and $\gamma_{\alpha,t,i}$ are defined in \eqref{eqn: loss_function_AL}.
Note that $L(\cdot)$ depends on $\{\wt\alpha_{t,i}\}_{t\in[T]}$, which is not implied from our notation for the ease of presentation. The DAFL estimators for the column effects can be obtained similarly, denoted by $\{\wh\beta_{t,j}\}_{t\in[T]}$ for $j\in[q]$. Then the sparse block estimators and their corresponding dense block estimators for the row and column main effects can be respectively defined as follows, for $i\in[p]$, $j\in[q]$,
\begin{equation}
\label{eqn: block_estimator}
\begin{split}
    & \wh\cS_{\alpha,i} := \{t: \wh\alpha_{t,i}\leq 0 \}, \quad
    \wh\cS_{\beta,j} := \{t: \wh\beta_{t,j}\leq 0 \}, \\
    & \wh\cB_{\alpha,i} :=[T] \setminus \wh\cS_{\alpha,i}, \quad
    \wh\cB_{\beta,j} := [T] \setminus \wh\cS_{\beta,j}.
\end{split}
\end{equation}
Note that we expect the DAFL estimators to be non-negative since the initial estimators are non-negative by how they are constructed. Hence, $\wh\cS_{\alpha,i}$ and $\wh\cS_{\beta,j}$ should technically correspond to the periods where $\wh\alpha_{t,i}$ and $\wh\beta_{t,j}$ are exactly zero, which might not hold empirically due to outliers or over-penalization. We thus define the sparse blocks as in \eqref{eqn: block_estimator} to address this; see more details in Remark~\ref{remark: final_estimator}.

Lastly, we construct the final estimators for the main effects by replacing all entries in $\wt\balpha_t$ and $\wt\bbeta_t$ by zero except for those according to the estimated dense blocks. That is, for each $t\in[T]$, $i\in[p]$, $j\in[q]$,
\begin{align*}
    \wt\alpha_{t,i}^{\cB} := \begin{cases}
    \wt\alpha_{t,i}, & t\in \wh\cB_{\alpha,i}; \\
    0, & \text{otherwise.}
    \end{cases}, \quad
    \wt\beta_{t,j}^{\cB} := \begin{cases}
    \wt\beta_{t,j}, & t\in \wh\cB_{\beta,j}; \\
    0, & \text{otherwise.}
    \end{cases}
\end{align*}
This unconventional step is to compensate for the use of \text{e.g.} $\wt\alpha_{t,i}$ in \eqref{eqn: DAFL}, which already induces an estimation error. Thus, even though the DAFL estimator $\{\wh\balpha_t \}_{t\in[T]}$ can consistently recover the sparse and dense blocks, it inevitably inherits the error from $\wt\alpha_{t,i}$ and the rate of convergence deteriorates in general, \text{cf.} the discussion in Remark~\ref{remark: final_estimator}. Intuitively, $\{\wt\balpha_{t}^{\cB} \}_{t\in[T]}$ can be regarded as a thresholding estimator based on $\{\wt\balpha_{t}\}_{t\in[T]}$, except that the thresholding is carried out by solving a constrained optimization problem with desirable features. Those final estimators are a compromise for the fact that main effects are latent, pinpointing again the difficulty of the sparsity problem in factor models.

\subsection{Theoretical guarantee on sparsity recovery}\label{subsec: thm_sparse_guarantee}

Before we show the properties of the DAFL estimators, we require an additional rate assumptions as follows.

\begin{itemize}
  \item[(R2)] (Further rate assumptions).
{\em
We assume that,
\begin{align*}
    \lambda_\alpha^{-1} q^{-1} \log\Big(p \sum_{i=1}^p |\cS_{\alpha,i}|\Big)
    , \;\;\;
    \lambda_\beta^{-1} p^{-1} \log\Big(q \sum_{j=1}^q |\cS_{\beta,j}|\Big) &= o(1) , \\
    \Big\{ q^{-1} \log\Big(p \sum_{i=1}^p |\cB_{\alpha,i}|\Big) + \lambda_\alpha \Big\} \Big( \min_{i\in[p]} \min_{t\in\cB_{\alpha,i}} \{\alpha_{t,i}^{\ast 2}\} \Big)^{-1} &= o_P(1) ,\\
    \Big\{ p^{-1} \log\Big(q \sum_{j=1}^q |\cB_{\beta,j}|\Big) + \lambda_\beta \Big\} \Big( \min_{j\in[q]} \min_{t\in\cB_{\beta,j}} \{\beta_{t,j}^{\ast 2}\} \Big)^{-1} &= o_P(1) .
\end{align*}
}
\end{itemize}

Assumption~(R2) is a set of rate assumptions required for block consistency to hold in general, restricting the sizes of the sparse and dense blocks, the tuning parameters, and the behavior of the non-zero main effects. Note that we presume
\begin{align*}
    \min_{i\in[p]} \min_{t\in\cB_{\alpha,i}} \{\alpha_{t,i}^{\ast}\}, \, \min_{j\in[q]} \min_{t\in\cB_{\beta,j}} \{\beta_{t,j}^{\ast}\} =O_P(1) ,
\end{align*}
since otherwise block consistency can be trivially obtained by any constant thresholding, which is unrealistic. As long as the initial estimators of the main effects are consistent with rates according to Theorem~\ref{thm: particular_initial_consistency_maineffect}~(ii), we can read the first line in Assumption~(R2) as requiring $\log(T)/(q\lambda_\alpha + p\lambda_\beta) \to 0$. On the other hand, the second and third lines in (R2) restricts the tuning parameters to grow slower than the squared minimum non-zero main effects. Those constraints involving the sizes of the sparse and dense blocks are in parallel to standard assumptions on the zero and non-zero coefficients in variable selection in linear regression, cf.\ Assumption~(A4) in \cite{huang2008adaptive}. In what follows, we present the block consistency of the DAFL estimators, which further induces the results for the final estimators $\{\wt\balpha_t^\cB\}_{t\in[T]}$ and $\{\wt\bbeta_t^\cB\}_{t\in[T]}$.

\begin{theorem}[Properties of DAFL estimators]
\label{thm: oracle}
Let all assumptions in Theorem~\ref{thm: particular_initial_consistency_maineffect} hold,
with Assumption~(IC1)~(ii)--(iii) as in \eqref{eqn: IC1_min0}.
Further given Assumption~(R2),
and consider $\min\{p,q,T\} \to\infty$.
Then the DAFL estimators in Section~\ref{subsec: estimation_regularize} are consistent, and block consistency holds with
\begin{align*}
    \b{P}(\wh\cS_{\alpha,1} = \cS_{\alpha,1}, \dots, \wh\cS_{\alpha,p} = \cS_{\alpha,p}) \to 1 , \quad
    \b{P}(\wh\cS_{\beta,1} = \cS_{\beta,1}, \dots, \wh\cS_{\beta,q} = \cS_{\beta,q}) \to 1 .
\end{align*}
\end{theorem}

\begin{corollary}\label{cor: final_estimator}
Under the assumptions in Theorem~\ref{thm: oracle}, the final estimators for the main effects have the following properties.
\begin{itemize}
    \item [(i)] (Block consistency). As $\min\{p,q,T\} \to\infty$, it holds with probability tending 1 that
    \begin{align*}
        \{t: \wt\alpha_{t,i}^\cB =0\} &= \cS_{\alpha,i} \quad \text{for all $i\in[p]$} , \\
        \{t: \wt\beta_{t,j}^\cB =0\} &= \cS_{\beta,j} \quad \text{for all $j\in[q]$} .
    \end{align*}
    \item [(ii)] (Uniform convergence). The final estimators in the dense blocks are consistent with
        \begin{align*}
            \max_{i\in[p]} \max_{t\in\cB_{\alpha,i}} \big|\wt\alpha_{t,i}^\cB -\alpha_{t,i}^\ast \big| &= O_P\Big\{ q^{-1/2} \log^{1/2}\Big(p \sum_{i=1}^p |\cB_{\alpha,i}|\Big) \Big\} ,\\
            \max_{j\in[q]} \max_{t\in\cB_{\beta,j}} \big|\wt\beta_{t,j}^\cB -\beta_{t,j}^\ast \big| &= O_P\Big\{ p^{-1/2} \log^{1/2}\Big(q \sum_{j=1}^q |\cB_{\beta,j}|\Big) \Big\} .
        \end{align*}
\end{itemize}
\end{corollary}

Theorem~\ref{thm: oracle} is key to the desired properties in the final estimators. Note that although the DAFL estimators are consistent, to derive a comparable rate of convergence as Theorem~\ref{thm: particular_initial_consistency_maineffect}~(ii) requires stricter rate conditions than Assumption~(R2). This is circumvented in the final estimators, as shown in Corollary~\ref{cor: final_estimator}~(ii). Intuitively, the final estimators for the main effects treat the DAFL estimation as a thresholding procedure. In terms of thresholding, by Assumption~(R2), we require in the probability sense that the squared row main effects to asymptotically dominate the rate $q^{-1} \log\big(p \sum_{i=1}^p |\cS_{\alpha,i}|\big) + q^{-1} \log\big(p \sum_{i=1}^p |\cB_{\alpha,i}|\big)$ which is effectively of order $q^{-1} \log(pT)$; similar arguments hold for the column main effects. Finally, with the block consistency result from Theorem~\ref{thm: oracle} and hence Corollary~\ref{cor: final_estimator}~(i), we present the behaviors of the main effect estimators in the dense blocks as in Corollary~\ref{cor: final_estimator}~(ii).

\begin{remark}\label{remark: final_estimator}
Theorem~\ref{thm: oracle} in Section~\ref{subsec: theorem} shows that the DAFL estimators are only consistent with arbitrary rates, compared to the rate from the initial estimators according to Theorem~\ref{thm: particular_initial_consistency_maineffect}~(ii), unless more restrictive conditions hold. A more practical concern in directly using the DAFL estimators is related to Assumption~(IC1). As explained below \eqref{eqn: block_estimator}, the DAFL estimators might be negative in practice and would not fulfill Assumption~(IC1) under $g(\cdot) = \min\{\cdot\} = 0$, requiring that all main effects are at least zero. On the other hand, it is easy to see that the initial estimators from \eqref{eqn: tilde_balpha_t} and \eqref{eqn: tilde_bbeta_t} satisfy the identification. This property is retained for the final estimators due to the way they are constructed.
\end{remark}


\subsection{Practical implementation}\label{subsec: practical_implementation}
We now discuss the practical optimization to compute our DAFL estimators. We only focus on estimating the row  main effects by minimizing \eqref{eqn: DAFL}, as the arguments for the column main effects follow similarly. It turns out that our DAFL estimators can be obtained by equivalently solving a generalized lasso problem \citep[][]{TibshiraniTaylor2011}. 
More specifically, for $i\in[p]$, if we define $\balpha_{\cdot,i} := (\alpha_{1,i}, \dots, \alpha_{T,i})^\top$, $\D_{\alpha,i} := \big(\D_{\alpha,i}^{(\textnormal{\text{F}}) \top}, \D_{\alpha,i}^{(\textnormal{\text{L}}) \top} \big)^{\top}$, where
{\small
\begin{align*}
    \D_{\alpha,i}^{(\textnormal{\text{L}})} &:= \diag\big( \{1/\wt\alpha_{1,i}, \dots, 1/\wt\alpha_{T,i}\} \big) , \\
    \D_{\alpha,i}^{(\textnormal{\text{F}})} &:= 
    \begin{pmatrix}
        -\left(\wt\alpha_{1,i}\lor\wt\alpha_{2,i}\right)^{-1} & \left(\wt\alpha_{1,i}\lor\wt\alpha_{2,i}\right)^{-1} & 0 & \cdots & 0 \\
         0 & -\left(\wt\alpha_{2,i}\lor\wt\alpha_{3,i}\right)^{-1} & \left(\wt\alpha_{2,i}\lor\wt\alpha_{3,i}\right)^{-1} & \cdots & 0 \\
         \vdots & \vdots & \vdots & \ddots & \vdots \\
         0 & 0 & \cdots & -\left(\wt\alpha_{T-1,i}\lor\wt\alpha_{T,i}\right)^{-1} & \left(\wt\alpha_{T-1,i}\lor\wt\alpha_{T,i}\right)^{-1}
    \end{pmatrix}
    ,
\end{align*}
}
then \eqref{eqn: DAFL} can be rewritten as
\begin{equation}
\label{eqn: DAFL_generalized_lasso}
L(\balpha_{\cdot,i}) = \frac{1}{2} \big\| \wt\balpha_{\cdot,i} - \balpha_{\cdot,i} \big\|_2^2 + \lambda_\alpha \big\| \D_{\alpha,i} \balpha_{\cdot,i} \big\|_1 ,
\end{equation}
which has the form of Equation~(2) in \cite{TibshiraniTaylor2011} and the solution path can be computed by algorithms therein in $O(T^3)$ times; see Algorithm~2 and the discussion in Section~8 in \cite{TibshiraniTaylor2011}. It might be worth pointing out that although one can stack all rows and directly compute the solution path on $(\balpha_{\cdot,1}^\top, \dots, \balpha_{\cdot,p}^\top)^\top$, this is not recommended since the computational complexity is of order $T^3p^3$, compared to $T^3p$ by solving the above problem for each $i\in[p]$.

For such an $\ell_1$ penalized regression problem as \eqref{eqn: DAFL_generalized_lasso}, it is often of interest to study its degrees of freedom which characterizes the effective number of parameters of a fitting procedure \citep{Efron1986, Zouetal2007}. In brief, for a data vector $\mathbf{y}\in\b{R}^T$ whose elements are uncorrelated with homoscedastic mean $\mu$ and variance $\sigma^2$, the degrees of freedom of a function $g:\b{R}^T \to \b{R}^T$ is defined as
\[
\text{df}\,(g) := \frac{1}{\sigma^2} \sum_{t=1}^T \text{Cov}\{ g_t(\mathbf{y}), y_t\} .
\]
As our DAFL estimator $\wh\balpha_{\cdot,i}$ is obtained by minimizing \eqref{eqn: DAFL_generalized_lasso}, we are interested in $\text{df}\,(\wh\balpha_{\cdot,i})$ in terms of approximating the vector $\wt\balpha_{\cdot,i}$. For a fixed tuning parameter $\lambda_\alpha$, we may apply Theorem~1 in \cite{TibshiraniTaylor2011} and leverage the relation between the primal and dual solutions to conclude
\[
\text{df}\,(\wh\balpha_{\cdot,i}) = \b{E}\{ \text{nullity}(\D_{\cA})\} ,
\]
where $\D_\cA$ denotes the matrix $\D_{\alpha,i}$ with rows restricted on the index set $\cA =\{t: (\D_{\alpha,i} \wh\balpha_{\cdot,i})_t =0\}$. Without the adaptive terms $u_{\alpha,t,i}$ and $\gamma_{\alpha,t,i}$ in \eqref{eqn: DAFL_generalized_lasso}, our problem boils down to the sparse fused lasso problem in Equation~(47) in \cite{TibshiraniTaylor2011} and the degrees of freedom can be reduced to the expected number of non-zero fused groups in $\wh\balpha_{\cdot,i}$. Although we do not have such interpretation in our complicated scenario, we may readily use the realized $\text{nullity}(\D_{\cA})$ and hence compute an estimate of the degrees of freedom. This allows us to modify the Mallows's $C_p$ statistic \citep{Mallows1973} as
\[
\wh{C}_p(\lambda_\alpha) := \big\| \wt\balpha_{\cdot,i} - \wh\balpha_{\cdot,i} \big\|_2^2 -T\wh\sigma^2 + 2\wh\sigma^2 \text{nullity}(\D_{\cA}) ,
\]
where $\wh\sigma^2$ is the sample variance of $\wt\balpha_{\cdot,i}$. Therefore, in terms of model selection, we may choose the tuning parameter $\lambda_\alpha$ to minimize $\wh{C}_p(\lambda_\alpha)$, among a grid of candidates. Note that we use the same $\lambda_\alpha$ over $i\in[p]$ in \eqref{eqn: DAFL}, so that the theoretical guarantee holds uniformly. In practice, we can either minimize the aggregated $\wh{C}_p(\lambda_\alpha)$ over all $i\in[p]$ or, more generally, in our numerical experiments, we select different $\lambda_\alpha$'s for each $i$.

\section{Numerical Analysis}\label{sec: numerical_results}
\subsection{Simulation}\label{subsec: simulation}
In this subsection, we showcase the numerical performance of the proposed estimators using Monte Carlo experiments. We first evaluate the accuracy of the main effect and grand mean estimators  under various linear and quantile identifications.
We then examine the accuracy of the estimated common components and demonstrate their invariance to the identification imposed on the main effects. 
Finally, under the minimum-type identification, we investigate the performance of the DAFL and final estimators in terms of sparsity recovery and estimation accuracy.

\paragraph{Data generating process.}
We first generate grand mean series $\{\mu_t^\circ\}$, main effect series $\{\balpha_t^\circ\}$ and $\{\bbeta_t^\circ\}$. The grand mean and each entry of the two main effect processes independently follow AR(2) processes with coefficients (0.5,-0.2) and i.i.d.\ $\cN(0,1)$ innovations. Based on these, we construct different versions of grand mean and main effects according to different shift-varying functions in Assumption~(IC1), as detailed in each experiment.

The noise and factor series $\E_t$ and $\F_t$ are formed as in Assumptions~(E1), (E2), (E3) and (F1). Specifically, elements in $\F_t$ are jointly independent and each follows a standardized AR(2) process with coefficients (0.5, -0.3). Elements in $\F_{e,t}$ and $\bepsilon_t$ are similarly constructed, except that the AR coefficients are (-0.4, 0.4) and (0.6, 0.2), respectively. Furthermore, elements in the standard deviation matrix $\bSigma_\epsilon$ are generated by \text{i.i.d.} $|\cN(0,1)|$. The innovation processes in generating $\F_t$, $\F_{e,t}$ and $\bepsilon_t$ are \text{i.i.d.} $\cN(0,1)$. 
To incorporate factor strengths, the row factor loading matrix is generated as $\A_r = \M_p\U_r\B_r$, where $\U_r\in\b{R}^{p\times k_r}$ consists of \text{i.i.d.} $\cN(0,1)$ elements, $\M_p$ is defined in Section~\ref{sec: introduction} so that item (i) of Assumption~(IC1) is satisfied, and $\B_r= \diag(p^{-\zeta_{r,1}}, \dots, p^{-\zeta_{r,k_r}})$. Note that $\zeta_{r,j}\in [0,0.5]$, with pervasive factors represented by $\zeta_{r,j}=0$ and weak factors otherwise. The column loading $\A_c$ is similarly generated. Each entry of $\A_{e,r}$ and $\A_{e,c}$ is \text{i.i.d.} standard normal and has probability of 0.95 being exact 0. Throughout the simulation, we fix $k_{e,r}=k_{e,c}=2$ and $k_r=k_c=2$.
All numerical results are based on 500 Monte Carlo replications, unless specified otherwise.

To evaluate estimation accuracy, for any parameter series $\btheta = \{ \btheta_t \}_{t\in[T]}$, where $\btheta_t$ can be a scalar, vector, or matrix, and its estimator $\wh\btheta =\{ \wh\btheta_t \}_{t\in[T]}$, we define the mean squared errors as
$$
\text{MSE}(\wh\btheta) := \frac{\sum_{t=1}^{T} \|\btheta_{t}- \wh\btheta_t\|_F^2}{Td_{\btheta}} ,
$$
where $d_{\btheta}$ denotes the number of elements in $\btheta_t$. Moreover, for any given $\Q$ and $\wh\Q$, we use the column space distance to measure their discrepancy:
$$
\c{D}(\Q, \wh\Q) := \Big\|\Q(\Q^\top\Q)^{-1}\Q^\top -\wh\Q(\wh\Q^\top\wh\Q)^{-1}\wh\Q^\top \Big\|,
$$
which is widely used in the literature such as \cite{chen2022factor}, among others. Unless otherwise stated, all reported performance measures are averaged over replications.

\paragraph{Accuracy of grand mean and main effects.}
Given shift-varying functions $g_\alpha$ and $g_\beta$ and their solution functions $\wt{g}_\alpha$ and $\wt{g}_\beta$ (cf.\ Assumptions~(IC1) and (IC2)), 
the identified parameters are constructed as
\begin{align*}
    \balpha_t^\ast:= \balpha_t^\circ -\wt{g}_\alpha(\balpha_t^\circ) \1_p, \qquad
    \bbeta_t^\ast:= \bbeta_t^\circ-\wt{g}_\beta(\bbeta_t^\circ) \1_q, \qquad
    \mu_t^\ast:=\mu_t^\circ + \wt{g}_\alpha(\balpha_t^\circ) + \wt{g}_\beta(\bbeta_t^\circ) .
\end{align*}
It follows that $g_\alpha(\balpha_t^\ast)=g_\beta(\bbeta_t^\ast)=0$, while 
\begin{align*}
    \mu_t^\ast\1_p\1_q^\top + \balpha_t^\ast\1_q^\top + \1_p\bbeta_t^{\ast\top}
    =\mu_t^\circ\1_p\1_q^\top + \balpha_t^\circ\1_q^\top + \1_p\bbeta_t^{\circ\top}.
\end{align*}
Hence, within each realization, the observed matrix time series $\{\X_t\}$ is exactly the same under any identification, while the identification only changes the values of the grand mean and main effects. This allows all identifications to be compared using the same set of generated data. Throughout this experiment, all factors are pervasive such that $\zeta_{r,j}=\zeta_{c,j}=0$.
We first experiment the linear identifications described in $\cL$. 
Given any positive integer $d$ and $s\in [d]$,
let $\v_{d,s}:=(\1_s^\top,\0_{d-s}^\top)^\top \in\b{R}^d$. We set
\begin{align*}
    g_\alpha(\aa)=\v_{p,s_p}^\top\aa
    \;\text{\ and \ }
    g_\beta(\bb)=\v_{q,s_q}^\top\bb.
\end{align*}
A smaller value of $s_p$ or $s_q$ corresponds to a more concentrated identification, with $s_p=s_q=1$ giving a reference-unit identification. 

\begin{figure}[t!]
\centering
\includegraphics[width=0.64\textwidth]{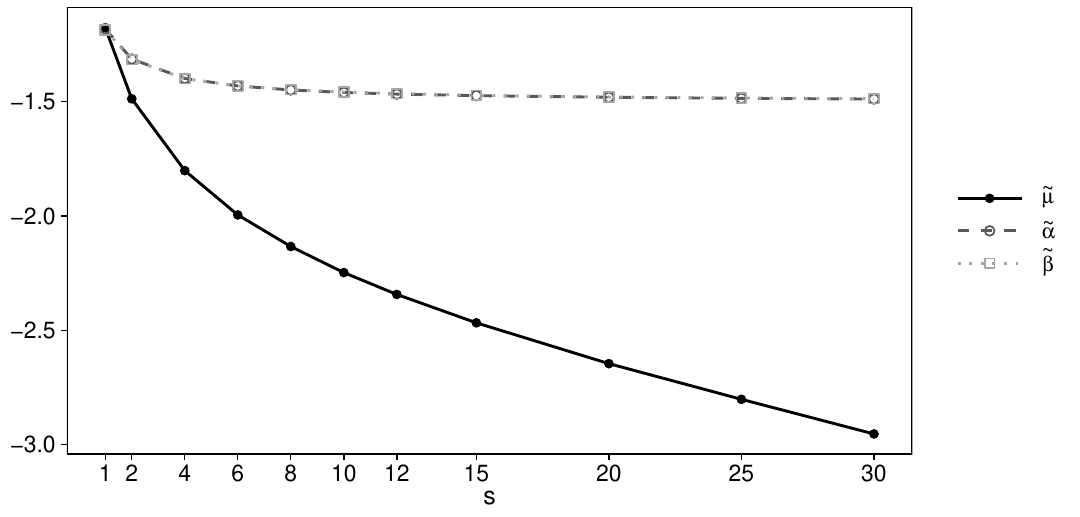}
\caption{MSEs under linear identification as the number $s$ varies (on the $\log_{10}$ scale).}
\label{fig:sim-linear-concentration}
\end{figure}

To examine the effect of weight concentration, we fix $T=100$ and $p=q=30$ and impose $s_p=s_q=s$. Figure~\ref{fig:sim-linear-concentration} shows how the accuracy of the grand mean and main effect estimators varies with $s$. The MSEs decrease as $s$ increases, showing that less concentrated identifications lead to more accurate estimation. The improvement is much stronger for the grand mean than for the main effects, consistent with Theorem~\ref{thm: particular_initial_consistency_maineffect}~(i).

We then examine the effects of varying $p$, $q$, and $T$ under the reference-unit and equal-weight identifications, corresponding to $(s_p,s_q)=(1,1)$ and $(p,q)$, respectively. The results reported in Table~\ref{tab:sim-linear-identification} are largely unchanged as $T$ increases. The equal-weight identification is substantially more accurate than the reference-unit identification, which is consistent with Figure~\ref{fig:sim-linear-concentration}. 
The grand mean becomes more accurate as either $p$ or $q$ increases. For the main effects, increasing $q$ improves the row main effect estimator but not the column main effect estimator, while increasing $p$ shows the reverse pattern. These patterns are also consistent with Theorem~\ref{thm: particular_initial_consistency_maineffect}~(i).

\begin{table}[t!]
\centering
\caption{Accuracy under linear identification. Entries are MSE means multiplied by $10^2$, with standard errors in parentheses.}
\label{tab:sim-linear-identification}
\resizebox{\textwidth}{!}{%
\begin{tabular}{ccccccccc}
\toprule
\multirow{2}{*}{$(s_p,s_q)$} & \multirow{2}{*}{$p$} & \multirow{2}{*}{$q$} & \multicolumn{2}{c}{$\MSE(\wt\mu)$} & \multicolumn{2}{c}{$\MSE(\wt\balpha)$} & \multicolumn{2}{c}{$\MSE(\wt\bbeta)$} \\
\cmidrule(lr){4-5}\cmidrule(lr){6-7}\cmidrule(lr){8-9}
& & & $T=100$ & $T=200$ & $T=100$ & $T=200$ & $T=100$ & $T=200$ \\
\midrule
\multirow{3}{*}{$(1,1)$} & 30 & 30 & $6.652(0.110)$ & $6.464(0.108)$ & $6.577(0.065)$ & $6.430(0.060)$ & $6.516(0.063)$ & $6.413(0.057)$ \\
 & 80 & 80 & $2.526(0.037)$ & $2.549(0.029)$ & $2.499(0.022)$ & $2.554(0.020)$ & $2.505(0.020)$ & $2.491(0.016)$ \\
 & 30 & 80 & $4.591(0.078)$ & $4.656(0.071)$ & $2.474(0.022)$ & $2.426(0.015)$ & $6.628(0.059)$ & $6.712(0.062)$ \\
\addlinespace[0.35em]
\multirow{3}{*}{$(p,q)$} & 30 & 30 & $0.111(0.002)$ & $0.112(0.001)$ & $3.245(0.011)$ & $3.237(0.010)$ & $3.256(0.012)$ & $3.249(0.010)$ \\
 & 80 & 80 & $0.016(0.000)$ & $0.016(0.000)$ & $1.250(0.002)$ & $1.247(0.002)$ & $1.251(0.002)$ & $1.249(0.002)$ \\
 & 30 & 80 & $0.042(0.001)$ & $0.042(0.000)$ & $1.223(0.004)$ & $1.223(0.003)$ & $3.332(0.007)$ & $3.339(0.007)$ \\
\bottomrule
\end{tabular}%
}
\end{table}

We next experiment the quantile identifications described in $\cQ$. Specifically, we set
\begin{align*}
    g_\alpha(\aa)=q_{u}(\aa)
    \;\text{\ and \ }
    g_\beta(\bb)=q_{u}(\bb),
\end{align*}
and vary the quantile level over $u\in\{0,0.25,0.5\}$. 
The results over different dimension combinations are reported in Table~\ref{tab:sim-quantile-identification}. 
The estimation accuracy improves as the identifying quantile moves toward the center. This is because tail quantiles are more sensitive to extreme cross-sectional noise, with $u=0$ being directly determined by the minimum. The results change little as $T$ increases, while the effects of the cross-sectional dimensions follow the same patterns as under linear identification.

\begin{table}[ht!]
\centering
\caption{Accuracy under quantile identification. Entries are MSE means multiplied by $10^2$, with standard errors in parentheses.}
\label{tab:sim-quantile-identification}
\resizebox{\textwidth}{!}{%
\begin{tabular}{ccccccccc}
\toprule
\multirow{2}{*}{$p$} & \multirow{2}{*}{$q$} & \multirow{2}{*}{$u$} & \multicolumn{2}{c}{$\MSE(\wt\mu)$} & \multicolumn{2}{c}{$\MSE(\wt\balpha)$} & \multicolumn{2}{c}{$\MSE(\wt\bbeta)$} \\
\cmidrule(lr){4-5}\cmidrule(lr){6-7}\cmidrule(lr){8-9}
& & & $T=100$ & $T=200$ & $T=100$ & $T=200$ & $T=100$ & $T=200$ \\ 
\midrule
\multirow{3}{*}{30} & \multirow{3}{*}{30} & 0 & $6.160(0.047)$ & $6.155(0.034)$ & $6.180(0.030)$ & $6.162(0.023)$ & $6.192(0.030)$ & $6.162(0.023)$ \\
 &  & 0.25 & $2.062(0.014)$ & $2.069(0.011)$ & $4.206(0.015)$ & $4.202(0.013)$ & $4.220(0.015)$ & $4.223(0.014)$ \\
 &  & 0.5 & $1.609(0.011)$ & $1.616(0.008)$ & $3.996(0.014)$ & $3.990(0.012)$ & $4.008(0.014)$ & $4.002(0.012)$ \\
\addlinespace[0.35em]
\multirow{3}{*}{80} & \multirow{3}{*}{80} & 0 & $2.362(0.015)$ & $2.358(0.011)$ & $2.400(0.009)$ & $2.401(0.007)$ & $2.402(0.009)$ & $2.400(0.007)$ \\
 &  & 0.25 & $0.579(0.004)$ & $0.574(0.003)$ & $1.529(0.003)$ & $1.523(0.003)$ & $1.534(0.003)$ & $1.530(0.003)$ \\
 &  & 0.5 & $0.458(0.003)$ & $0.459(0.002)$ & $1.471(0.003)$ & $1.469(0.002)$ & $1.474(0.003)$ & $1.471(0.002)$ \\
\addlinespace[0.35em]
\multirow{3}{*}{30} & \multirow{3}{*}{80} & 0 & $4.245(0.030)$ & $4.236(0.021)$ & $2.362(0.011)$ & $2.368(0.007)$ & $6.309(0.024)$ & $6.323(0.020)$ \\
 &  & 0.25 & $1.016(0.007)$ & $1.023(0.005)$ & $1.692(0.006)$ & $1.691(0.004)$ & $3.838(0.009)$ & $3.845(0.008)$ \\
 &  & 0.5 & $0.805(0.006)$ & $0.803(0.004)$ & $1.595(0.005)$ & $1.593(0.004)$ & $3.719(0.008)$ & $3.728(0.007)$ \\
\bottomrule
\end{tabular}%
}
\end{table}

\paragraph{Accuracy of common components.}
We next examine the estimation accuracy of $\wh\Q_r$, $\wh\Q_c$, and $\wh\C_t$. Since these estimators are invariant to the identification imposed on the main effects, we fix the median identification throughout this experiment, namely, $g_\alpha(\aa)=q_{0.5}(\aa)$ and $g_\beta(\bb)=q_{0.5}(\bb)$. The grand mean and main effects are constructed accordingly as in the preceding experiment. We vary $p$, $q$, $T$, and the factor strength to examine their effects on the estimation accuracy,
as reported in Table~\ref{tab:sim-common-component-accuracy}. 
In general, increasing $T$ or the cross-sectional dimensions improves estimation accuracy. 
Weak factors substantially reduce the accuracy of the loading estimators, 
while their effect on the common component is relatively mild. 
corroborating with Theorem~\ref{thm: general_initial_consistency}.

\begin{table}[ht!]
\centering
\caption{Loading-space and common component estimation accuracy, where $\zeta=\zeta_r=\zeta_c$. Entries are Monte Carlo means multiplied by $10^2$, with standard errors in parentheses.}
\label{tab:sim-common-component-accuracy}
\resizebox{\textwidth}{!}{%
\begin{tabular}{ccccccccc}
\toprule
\multirow{2}{*}{$p$} & \multirow{2}{*}{$q$} & \multirow{2}{*}{$\zeta$} & \multicolumn{2}{c}{$\c{D}(\Q_r,\wh \Q_r)$} & \multicolumn{2}{c}{$\c{D}(\Q_c,\wh \Q_c)$} & \multicolumn{2}{c}{$\MSE(\wh \C)$} \\
\cmidrule(lr){4-5}\cmidrule(lr){6-7}\cmidrule(lr){8-9}
& & & $T=100$ & $T=200$ & $T=100$ & $T=200$ & $T=100$ & $T=200$ \\ 
\midrule
\multirow{2}{*}{30} & \multirow{2}{*}{30} & $0$ & $2.016(0.022)$ & $1.445(0.017)$ & $1.976(0.022)$ & $1.455(0.020)$ & $0.637(0.004)$ & $0.550(0.003)$ \\
 &  & $0.2$ & $16.844(0.450)$ & $13.491(0.315)$ & $16.528(0.383)$ & $13.529(0.336)$ & $1.284(0.033)$ & $1.014(0.027)$ \\
\addlinespace[0.35em]
\multirow{2}{*}{80} & \multirow{2}{*}{80} & $0$ & $1.046(0.007)$ & $0.728(0.004)$ & $1.052(0.007)$ & $0.732(0.004)$ & $0.127(0.001)$ & $0.096(0.000)$ \\
 &  & $0.2$ & $10.064(0.158)$ & $8.601(0.304)$ & $10.020(0.135)$ & $8.567(0.314)$ & $0.235(0.004)$ & $0.238(0.033)$ \\
\addlinespace[0.35em]
\multirow{2}{*}{30} & \multirow{2}{*}{80} & $0$ & $1.179(0.012)$ & $0.873(0.009)$ & $1.751(0.015)$ & $1.227(0.009)$ & $0.288(0.001)$ & $0.231(0.001)$ \\
 &  & $0.2$ & $15.000(0.470)$ & $11.906(0.270)$ & $13.176(0.360)$ & $9.313(0.143)$ & $0.714(0.067)$ & $0.434(0.010)$ \\
\bottomrule
\end{tabular}%
}
\end{table}

\paragraph{Sparsity recovery under minimum identification.}
We investigate the sparsity recovery performance of the DAFL estimators and whether exploiting the recovered sparsity improves the final estimation of the main effects. Throughout, we impose the minimum identification, namely, $g_\alpha(\aa)=\min\{\aa\}$ and $g_\beta(\bb)=\min\{\bb\}$.
All factors are pervasive,
the grand mean, common component, and noise processes are generated as in the beginning of this section, except that the row and column main effects are generated from sparse processes described below.

We describe the generating mechanism for the row main effects, with the column main effects constructed analogously. For each $i\in[p]$, denote the stay-in probability for sparse blocks and dense blocks given some sparsity level by $\pi_{\alpha,i}^\cS$ and $\pi_{\alpha,i}^\cB$ respectively, such that
\begin{align}
    \b{P}(\alpha_{t+1,i}^\circ = 0 \mid \alpha_{t,i}^\circ = 0) = \pi_{\alpha,i}^\cS, \quad 
    \b{P}(\alpha_{t+1,i}^\circ > 0 \mid \alpha_{t,i}^\circ > 0) = \pi_{\alpha,i}^\cB,
    \label{eqn: stay-in-prob}
\end{align}
Essentially, larger stay-in probabilities imply longer piecewise blocks, while a larger $\pi_{\alpha,i}^\cS$ combined with a smaller $\pi_{\alpha,i}^\cB$ implies a larger $|\cS_{\alpha,i}|$, and vice versa. Define i.i.d.\ random variables $Z_{\alpha,t,i}\sim|\cN(0,1)|$ for $t\in[T]$ and $i\in[p]$. Then given constant $m_{\alpha}$, we generate $\alpha_{t,i}^\circ$ for $t\in[T]$ by the following steps:
\begin{enumerate}[itemsep=0pt, label = Step~\arabic*., left = 0pt]
    \item Define $p_{\alpha,i,\ast}:=(1-\pi_{\alpha,i}^{\cB})/(2-\pi_{\alpha,i}^{\cS}-\pi_{\alpha,i}^{\cB})$. Set $\alpha_{1,i}^\circ=0$ with probability $p_{\alpha,i,\ast}$, and set $\alpha_{1,i}^\circ=m_\alpha+Z_{\alpha,1,i}$ otherwise.
    \item For $t=2,\ldots,T$, conditional on the state at $t-1$, generate the state at $t$ according to \eqref{eqn: stay-in-prob}. Set $\alpha_{t,i}^\circ=0$ if the process is in the sparse block, and set $\alpha_{t,i}^\circ=m_\alpha+Z_{\alpha,t,i}$ if it is in the dense block.
\end{enumerate}

The following proposition presents the properties of the resulting process $\{\alpha_{t,i}^\circ\}_{t\in[T]}$.
\begin{proposition}\label{prop: dgp_sparsity}
For $\{\alpha_{t,i}^\circ\}_{t\in[T]}$ generated by Steps~1--2 above, it holds that:
(i) $\{\b1\{\alpha_{t,i}^{\circ} =0\}\}_{t\in[T]}$ is a stationary Markov process, with stationary probability
$\b{P}(\alpha_{t,i}^\circ=0) =p_{\alpha,i,\ast} := (1-\pi_{\alpha,i}^{\cB})/(2-\pi_{\alpha,i}^{\cS}-\pi_{\alpha,i}^{\cB})$;
(ii) with $\cS_{\alpha,i}^\circ := \{t:\alpha_{t,i}^\circ=0\}$, we have $\b{E}(|\cS_{\alpha,i}^\circ|) = Tp_{\alpha,i,\ast}$;
and (iii) the expected length of a consecutive sparse interval is $L_{\alpha,i}^{\cS}=(1-\pi_{\alpha,i}^\cS)^{-1}$.
\end{proposition}

Proposition~\ref{prop: dgp_sparsity} shows that the two stay-in probabilities provide a convenient way to control two distinct features of sparsity. In particular, $p_{\alpha,i,\ast}$ determines the overall proportion of sparse observations, whereas $L_{\alpha,i}^{\cS}$ determines their temporal persistence. In the experiments below, these parameters are taken to be homogeneous across units, and we write them as $p_\ast$ and $L_{\cS}$. For any admissible pair $(p_\ast,L_{\cS})$, the corresponding stay-in probabilities are
$\pi^{\cS}= 1- 1/L_{\cS}$
and $\pi^{\cB}= 1- \frac{p_\ast}{(1-p_\ast) L_{\cS}}$.
After generating $\balpha_t^\circ$ and $\bbeta_t^\circ$, we impose the minimum identification through the same shift transformation as in the preceding experiments. Specifically, we set
\begin{align*}
    \balpha_t^\ast=\balpha_t^\circ-\big(\min\{\balpha_{t}^{\circ}\}\big)\1_p, \;\;
    \bbeta_t^\ast=\bbeta_t^\circ-\big(\min\{\bbeta_{t}^{\circ}\}\big)\1_q, \;\;
    \mu_t^\ast=\mu_t^\circ+\min\{\balpha_{t}^{\circ}\}+\min\{\bbeta_{t}^{\circ}\}.
\end{align*}

For every replication, we obtain the DAFL estimators $\wh\balpha_t$ and $\wh\bbeta_t$ using the tuning parameters selected by the modified $C_p$ criterion in Section~\ref{subsec: practical_implementation}. We evaluate sparsity recovery using sensitivity and specificity. For the row main effects, define
\begin{align*}
    \text{sensitivity}_{\balpha} 
    := \frac{\sum_{i=1}^p|\wh\cB_{\alpha,i}\cap\cB_{\alpha,i}|} {\sum_{i=1}^p|\cB_{\alpha,i}|}, \quad
    \text{specificity}_{\balpha} 
    := \frac{\sum_{i=1}^p|\wh\cS_{\alpha,i}\cap\cS_{\alpha,i}|} {\sum_{i=1}^p|\cS_{\alpha,i}|}.
\end{align*}
Also, define a unit-level exact sparsity recovery rate as
$ \text{ESR}_{\balpha} := p^{-1} \sum_{i=1}^p \b1 \{\wh\cS_{\alpha,i} =\cS_{\alpha,i}\}$.
Lastly, to assess whether exploiting the estimated sparsity improves main-effect estimation, we report the relative MSE improvement
$\text{Impr}_\alpha := 1-\text{MSE}(\wt\balpha^{\cB}) / \text{MSE}(\wt\balpha)$.
Let the corresponding quantities $\text{sensitivity}_{\bbeta}$, $\text{specificity}_{\bbeta}$, $\text{ESR}_{\bbeta}$, and $\text{Impr}_\beta$ be analogously defined.

\begin{figure}[t!]
\centering
\includegraphics[width=0.98\textwidth]{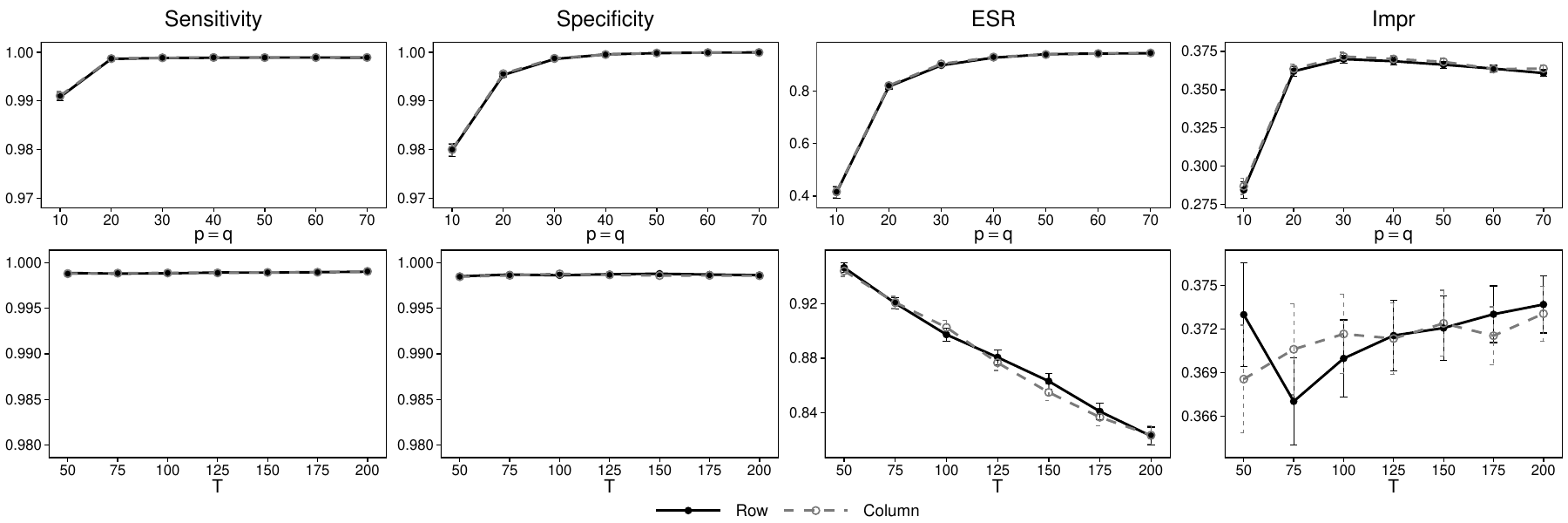}
\caption{Recovery performance. The columns report sensitivity, specificity, ESR, and Impr, respectively. The first row varies $p=q$ with $T=100$, while the second row varies $T$ with $p=q=30$. Curves show the row and column Monte Carlo means with 95\% Monte Carlo error bars.}
\label{fig:sim-dimension-time}
\end{figure}

We investigate sparsity recovery with respect to dimensions and series length, signal strength, and sparsity patterns. First, we vary the matrix dimensions and series length, while fixing a representative sparsity pattern with $p_\ast=0.4$ and $L_{\cS}=4$, and setting $m_\alpha=m_\beta=2$. 
Figure~\ref{fig:sim-dimension-time} presents the recovery performance. Except for $p=q=10$, the sensitivity and specificity both exceed $0.99$ and quickly approach one as $p=q$ increases. 
ESR shows a stronger dimensional effect and also approaches one as the cross-sectional dimensions increase. As $T$ increases, the sensitivity and specificity remain nearly unchanged near one, while the ESR gradually decreases but remains above $0.8$. 
This is natural because ESR requires the sparse block of a unit to be recovered exactly over its entire time trajectory. 
A longer series therefore makes exact recovery more demanding, even when the overall recovery accuracy remains nearly unchanged. Furthermore, $\text{Impr}_\alpha$ and $\text{Impr}_\beta$ remain well above zero across all settings. Hence, exploiting the recovered sparsity indeed improves the estimation accuracy of the main effects relative to the initial estimators, showing that sparsity recovery also leads to estimation gains.

\begin{figure}[t!]
\centering
\includegraphics[width=0.98\textwidth]{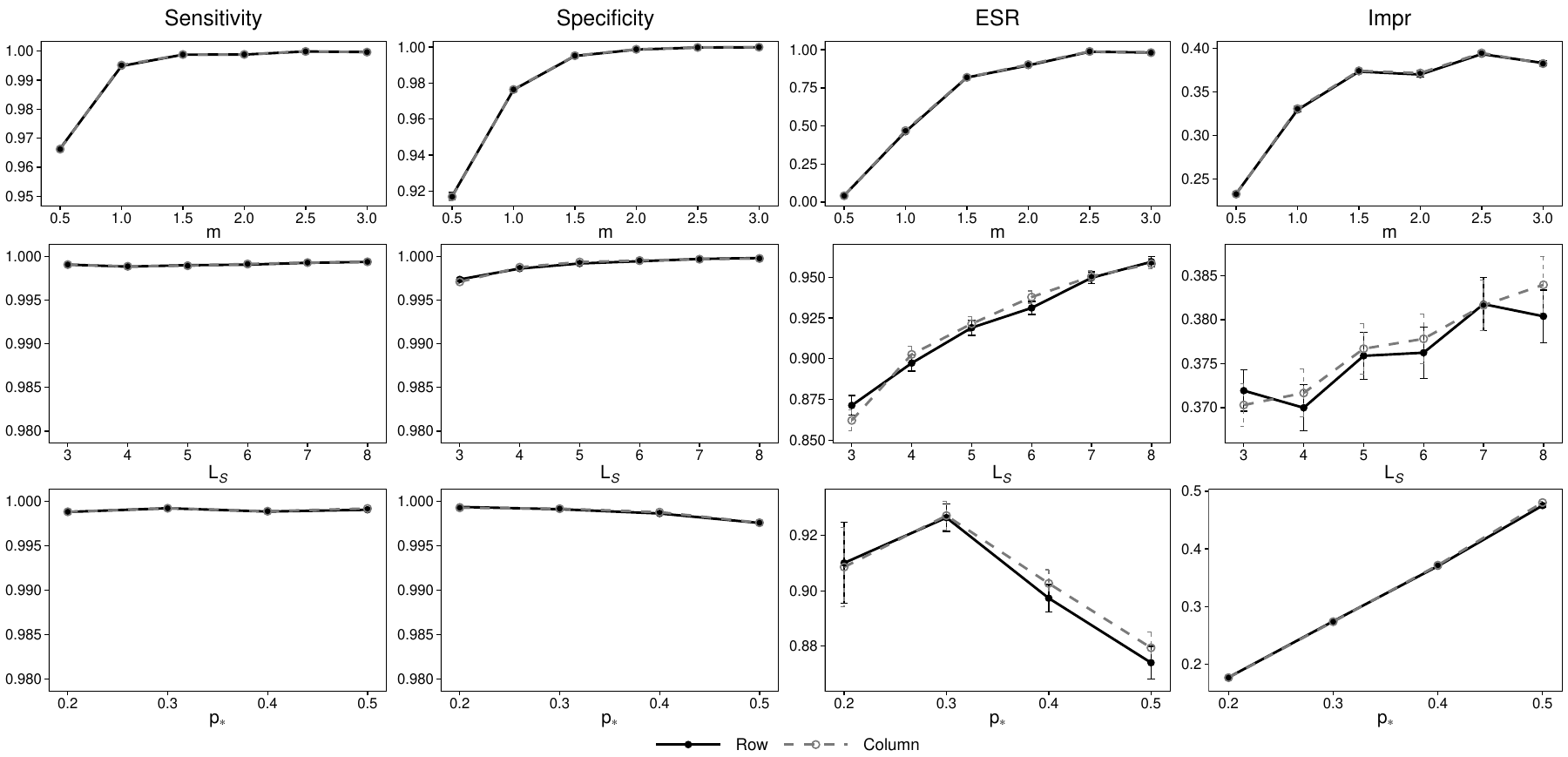}
\caption{Recovery performance. The rows vary signal strength, sparse-block persistence, and sparsity level, respectively. The columns report sensitivity, specificity, ESR, and Impr, respectively. Curves show the row and column Monte Carlo means with 95\% Monte Carlo error bars.}
\label{fig:sim-persistence-sparsity}
\end{figure}

Second, we examine the effect of signal strength under $p=q=30$ and $T=100$, while fixing $p_\ast=0.4$ and $L_{\cS}=4$. The first row of Figure~\ref{fig:sim-persistence-sparsity} shows that both sensitivity and specificity rapidly approach one as $m=m_\alpha=m_\beta$ increases. 
ESR increases from near zero under the weakest signal to nearly one. Under weak signals, specificity is lower than sensitivity, indicating that true sparse observations are more often classified as dense. When $m$ is small, sparse and dense effects are less well separated relative to estimation noise, so some sparse observations receive insufficient adaptive penalization and remain positive. As $m$ increases, sparse and dense effects become better separated, which is also consistent with Assumption~(R2).

We next examine the effects of sparse-block persistence and sparsity level under the same $p=q=30$, $T=100$, and $m_\alpha=m_\beta=2$ setting, as shown in the second and third rows of Figure~\ref{fig:sim-persistence-sparsity}. When $p_\ast=0.4$ is fixed, sensitivity and specificity remain close to one as $L_{\cS}$ varies, while ESR steadily approaches one as $L_{\cS}$ increases. This pattern is consistent with the role of the fused term in DAFL, which exploits longer piecewise sparse blocks for sparsity recovery.
When $L_{\cS}=4$ is fixed, increasing $p_\ast$ leaves $\pi^{\cS}=0.75$ unchanged but decreases $\pi^{\cB}$. Hence, the sparse blocks have the same expected length, while the dense blocks become shorter and sparse-dense transitions occur more frequently.

Across settings, $\text{Impr}_\alpha$ and $\text{Impr}_\beta$ remain positive. They generally increase as the signal strength or sparse-block persistence increases, following the improvement in sparsity recovery. The increase is particularly clear as $p_\ast$ increases, since a larger proportion of true zeros allows the final estimators to remove estimation error from more sparse entries.
Overall, under the minimum identification, DAFL accurately recovers the sparse blocks across a wide range of settings, with stronger signals and longer sparse blocks further improving recovery. Exploiting the recovered sparsity also consistently improves the final main effect estimators over their initial counterparts.

\subsection{A real application on labor market}\label{subsec: application_QCEW}

Consider the QCEW data introduced in Section~\ref{sec: introduction}, where the details on the dataset are included in the supplementary material.
To showcase our sparsity recovery procedure, we now in turn consider the identification
$\min\{\balpha_{t}^\ast\}=0$. The minimum is natural on the state categorization because its zero set identifies the lower boundary relevant for persistent downturn in state employment. 
This can be also of particular interest due to the fact that the raw minimum is highly unstable over time: among the 28 states, 25 attain the exact minimum at least once, and the minimum state changes in 41.5\% of adjacent months. We can therefore apply our procedure for sparsity recovery in Section~\ref{sec: sparsity_ME}, with the tuning parameter selected by the modified $C_p$ criterion in Section~\ref{subsec: practical_implementation}. DAFL exploits temporal persistence to recover a more stable zero set from the noisy lower-boundary estimates. 
Such the recovered zeros allow us to identify states with persistent lower-tail employment growth. 
For the sector main effect, we simply fix $\med\{\bbeta_{t}^\ast\}=0$ and omit its interpretation because of space limitation.

\begin{figure}[ht!]
    \centering
    \includegraphics[width=\textwidth]{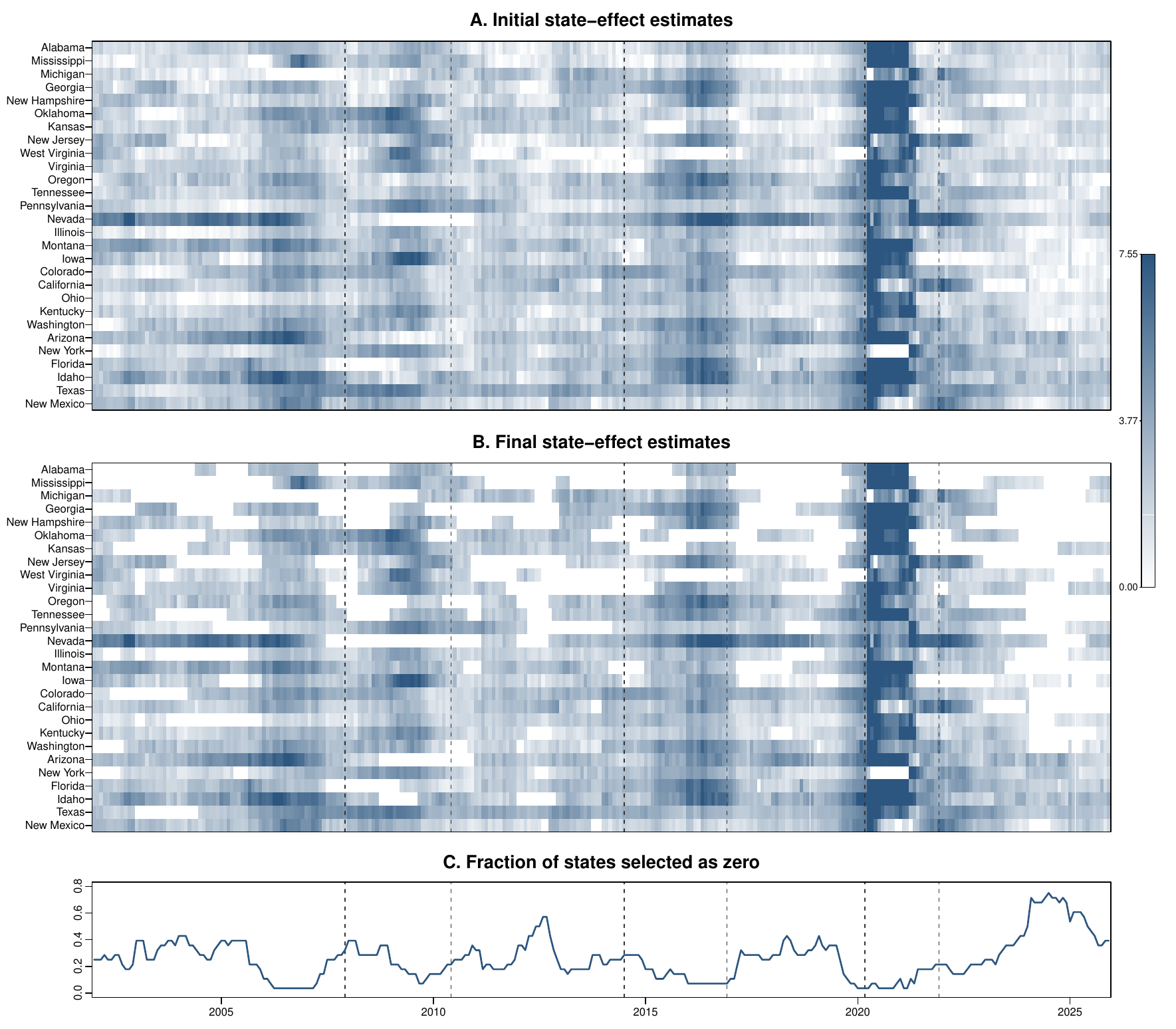}
    \caption{Initial and final state-effect estimates in the QCEW data. Panels~A and B show the initial estimates and the final estimates after DAFL zero-set recovery, respectively, using a common color scale. Panel~C reports the fraction of zero state effects over time. 
    Dashed vertical lines enclose the periods representing the GFC (line~1--2), the 2014--16 energy decline (line~3--4), and the pandemic (line~5--6).}
    \label{fig:qcew_persistent_regimes}
\end{figure}

Figure~\ref{fig:qcew_persistent_regimes} compares the initial and final state-effect estimates. 
Panel~A shows substantial month-to-month variation, whereas Panel~B reveals persistent zero blocks associated with recognizable events, such as Michigan before the GFC, Nevada and Florida during the housing bust, West Virginia during the 2014--16 energy decline, and New York during the pandemic around 2020.
We also report in Panel~C the fraction of estimated zero sets.
During the 2014--16 energy decline and the pandemic, such fractions are relatively narrow while much of the remaining cross section lies farther above the lower boundary. These large shocks widen cross-state differences and make the states with persistent low employment growth more distinct from the rest.

\begin{figure}[ht!]
    \centering
    \includegraphics[width=\textwidth]{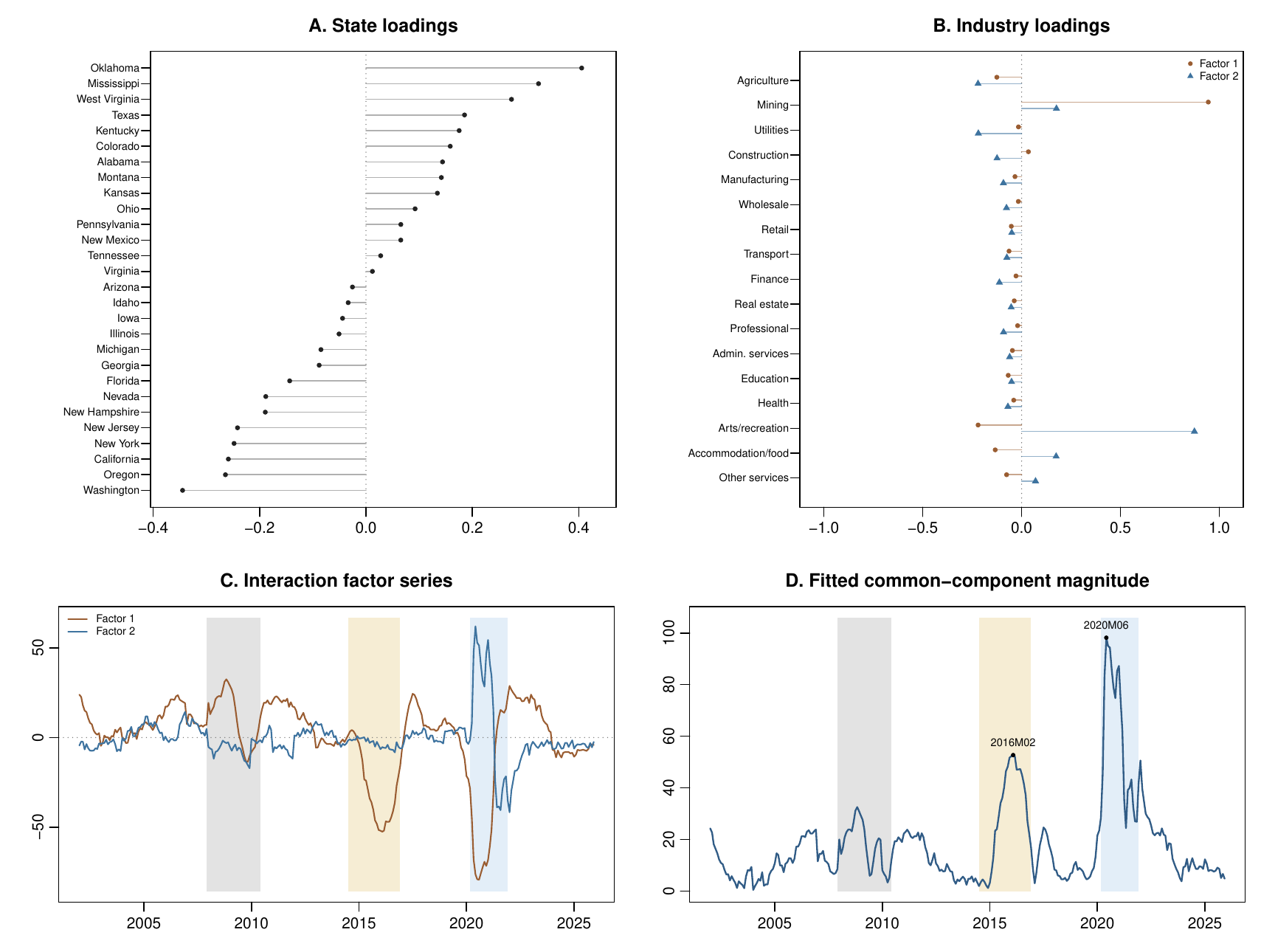}
    \caption{Estimated factor structure for the QCEW panel under $(\wh k_r,\wh k_c)=(1,2)$. Panels~A and B show the state and industry loadings, Panel~C the factor series, and Panel~D the Frobenius norm of the fitted common component. 
    Shaded areas mark the GFC (gray), the 2014--16 energy decline (orange), and the pandemic period (blue).}
    \label{fig:qcew_common_component}
\end{figure}

Next, we examine the estimated factor structure.
To select the factor numbers $(k_r, k_c)$,
we experiment the estimators in \cite{Yuetal2022} and Section~4.6 of \cite{lam2024matrix} (which can be viewed as \cite{Zhangetal2024} for matrix data) on $\wt\L_t$ in \eqref{eqn: tilde_L_t}, with both resulting in
$(\wh k_r,\wh k_c)=(1,2)$.
Figure~\ref{fig:qcew_common_component} summarizes the estimated common component. Panel~A shows a single state loading direction, with Oklahoma, Mississippi, West Virginia, and Texas toward one end and Washington, Oregon, California, and New York toward the other, indicating different patterns of industry-specific departures across states. In Panel~B, one estimated industry-loading direction is dominated by mining, and its factor series moves sharply during the 2014--16 energy decline shown in Panel~C. 
The second direction is dominated by arts and recreation, with a smaller contribution from accommodation and food services, and its factor shows significant movements during the pandemic, thus showing how the life styles in U.S.\ were affected by the pandemic.

We also display in Panel~D the magnitude of the estimated common components over time. 
From the diagram, the magnitude rises sharply during both episodes and reaches its sample maximum during the pandemic. The GFC is less prominent in the common component, compared to the peaks in the energy decline and the pandemic.

Taken together, 
this QCEW application highlights a different role of identification from that in Section~\ref{sec: introduction}.
The minimum identification makes the lower boundary of the state effects an object of interest, and DAFL turns its unstable monthly estimates into persistent slow employment growth that are more interpretable. 
The common component captures complementary state-industry interactions, most clearly during the energy decline and the pandemic.

\subsection{A real application on macroeconomic indices}
Motivated by \cite{chen2020constrained} and \cite{Yuetal2022},
we study the quarterly macroeconomic data from Organization for Economic Co-operation and Development (OECD). 
At each quarter $t$, $\X_t\in\b{R}^{8\times 10}$ contains 10 macroeconomic indicators for eight countries: the United States, the United Kingdom, Canada, France, Germany, Norway, Australia and New Zealand. The indicators consist of three price measures (total, energy and core CPI), two interest rates (3-month and long-term rates), three real-activity measures (GDP, industrial production and manufacturing production), and two trade measures (goods exports and imports).
The data contains $T=128$ quarters from 1988-Q3 to 2020-Q2.

For identification on the main effects,
we use the median country and the real-activity group as our primary choice, such that
\begin{align*}
    &\text{Identification I} :\quad
    \med\{\balpha_{t}^\ast\}=0, \quad
    (\beta_{t,\text{GDP}}^\ast + \beta_{t,\text{IP}}^\ast + \beta_{t,\text{Manuf}}^\ast)/3 = 0 . 
\end{align*}
Hence, the country effects are measured relative to a contemporaneous median country, while the indicator effects are measured relative to the average real-activity dynamics.  
Two alternative sets of identification are also considered: 
on the country side, we use the United States as a reference unit
while retaining the real-activity group as above;
and on the indicator side, we retain the median-country identification but use the three price indicators as the reference group.
Formally,
\begin{align*}
    &\text{Identification II} :\quad
    \alpha_{t,\text{USA}}^\ast=0 , \quad
    (\beta_{t,\text{GDP}}^\ast + \beta_{t,\text{IP}}^\ast + \beta_{t,\text{Manuf}}^\ast)/3 = 0 ; \\
    &\text{Identification III} :\quad
    \med\{\balpha_{t}^\ast\}=0, \quad
    (\beta_{t,\text{CPI}}^\ast+\beta_{t,\text{Energy}}^\ast+\beta_{t,\text{Core}}^\ast)/3=0.
\end{align*}

\begin{figure}[ht!]
    \centering
    \includegraphics[width=\textwidth]{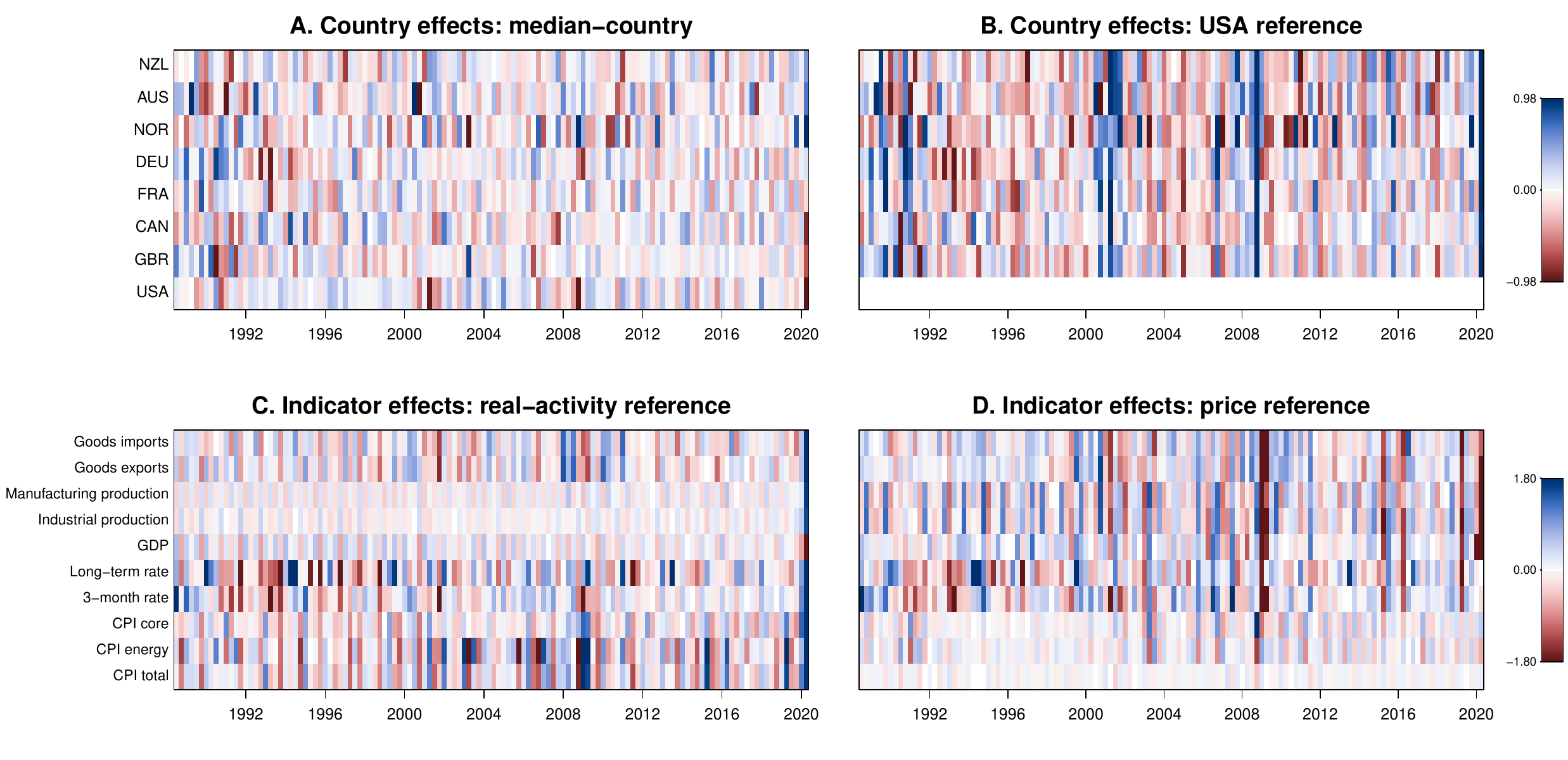}
    \caption{OECD main effects under different identifications. Panels~A and B show country effects under the median-country and U.S.\ reference identification, holding the real-activity reference fixed. Panels~C and D show indicator effects under the real-activity and price(-group) identification, holding the median-country (as in Panel~A) fixed. 
    Color scales are symmetrically clipped at the 98th percentile of absolute effect magnitudes for readability.}
    \label{fig:OECD_main_effect}
\end{figure}

Figure~\ref{fig:OECD_main_effect} shows the estimated main effects under these identifications. Panels~A and C report the estimated main effects under Identification I.
In Panel~A, 
around the global financial crisis (GFC) in 2008, the country main effect for the U.S.\ stepped below the median (of value zero) for several quarters, whereas that of Norway often moved in the opposite direction.
In Panel~C, indicator effects are measured relative to the average of GDP, industrial production and manufacturing production.
Around GFC in 2008, 
it is clear to see that 
Goods imports, Goods exports, and all CPI indicators stood prominently above our real-activity benchmark, indicating that, for instance, imports and exports were less influenced by GFC.
Also, the effects of Long-term rate and 3-month rate appeared more volatile in the early 1990s, while CPI energy stood out during the 2014--16 oil-price decline.

Panels~B and D provide alternative perspectives under Identifications~II and III. 
With the U.S.\ as the country reference in Panel~B, all other countries' main effects appeared above the U.S.\ baseline,
both during GFC around 2008 and the pandemic in 2020-Q2.
With the price group as the indicator reference in Panel~D, the influence by GFC around 2008 appeared in the opposite directions compared to that in Panel~C, which is due to real-activity reference moving well below the price reference.
Thus, the choice of identification could change the interpretations of main effects, while preserving the relative differences across countries and indicators. 

\begin{figure}[ht!]
    \centering
    \includegraphics[width=\textwidth]{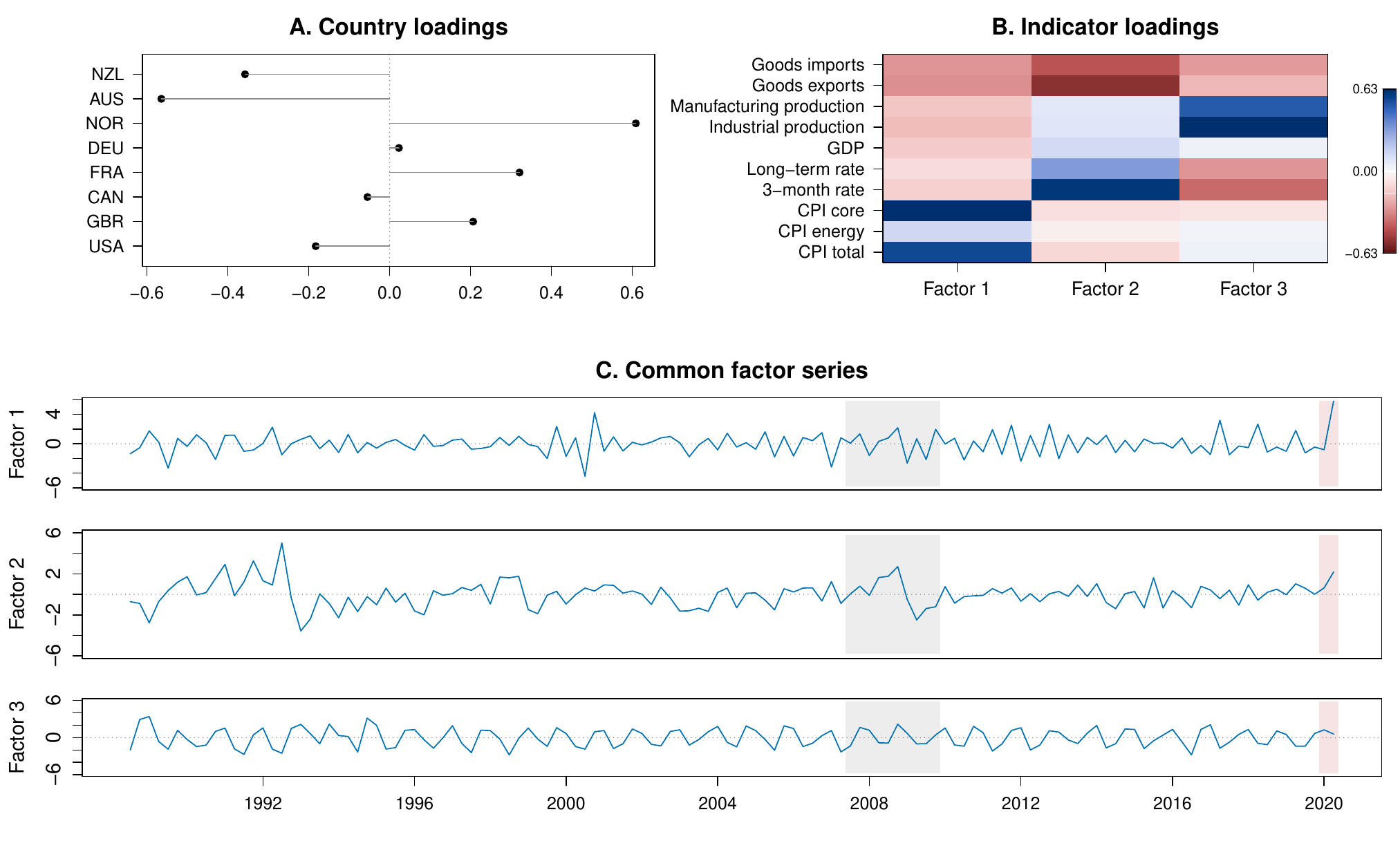}
    \caption{Estimated common factor structure for the OECD data. Panels~A and B show the country and indicator loadings, and Panel~C shows the three common factor series, with the GFC and pandemic periods shaded. 
    }
    \label{fig:oecd_factor_structure}
\end{figure}

We next examine the factor structure,
which, markedly, is estimated the same regardless of our identification on the main effects.
Similarly to Section~\ref{subsec: application_QCEW}, 
we apply the factor number estimators in \cite{Yuetal2022} and select $(\wh{k}_r,\wh{k}_c)=(1,3)$.
Figure~\ref{fig:oecd_factor_structure} reports the estimated loadings and factor series. 
Panel~A shows the estimated country loading matrix,
where a geographical grouping can be seen from the signs of the elements.
In Panel~B,
the first column can be regarded as CPI-related indices; the second is dominated by interest rates against trade; 
and the third represents industrial and manufacturing production with interest rates. 
Panel~C displays when these factors become more active. For instance, The second factor is more variable in the early 1990s and around the GFC, while the first shows a large movement in 2020-Q2.
Overall, the OECD application shows how general identification separates interpretable additive effects from a low-dimensional common component.

\section{Concluding Remarks}\label{sec: conclusion}
In this paper, we address a critical gap in the literature of matrix factor models by developing a flexible framework for MEFM identification based on shift‑varying functions, 
and by designing methodologies to recover sparse structures in the main effects. 
In particular, we reconcile the dual demands of interpretability for applied econometricians and structural parsimony inherent in matrix‑valued time series. 
Our general identification theory allows practitioners to choose domain‑appropriate constraints, from sum-to-zero to quantile-based or reference‑unit anchors,
while maintaining identifiability. Building on this, we propose a Doubly Adaptive Fused Lasso estimator that consistently recovers the sparse and dense blocks of the main effects, enabling a clear separation of periods with and without row‑ or column‑specific deviations.

Although we only consider sparsity in the main effects, it is possible to study approximately sparse loading matrices as in Assumption~2 in \cite{Freyaldenhoven2022}, or exactly sparse loadings as in Assumption~2 in \cite{UematsuYamagata2023}, among others. 
Their methods can be directly applied to the vectorization of the matrix time series $\wt{\L}_t$ ($t\in[T]$) defined in \eqref{eqn: tilde_L_t}, at the cost of overlooking the matrix nature in each observed data matrix and inflating number of parameters. 
These observations demand further development of sparsity in matrix factor models. In fact, methodologies designed for the classical factor models are all applicable to $\wt{\L}_t$ ($t\in[T]$) in practice, and hence the approach proposed in this paper can be used as a powerful pre‑processing of the observed data set. 
We hope that this article sheds light on sparsity structures in data analysis for matrix‑valued or even higher‑order tensor time series, and provides new insights in addressing dependency in complicated data structures.


\bibliographystyle{apalike}
\bibliography{ref}	


\clearpage
\appendix
\renewcommand{\thesection}{\Alph{section}}
\renewcommand{\theHsection}{\Alph{section}} 
\renewcommand{\thesubsection}{\Alph{section}.\arabic{subsection}}
\renewcommand{\theHsubsection}{\Alph{section}.\arabic{subsection}}

\numberwithin{equation}{section}
\numberwithin{figure}{section}
\numberwithin{table}{section}
\numberwithin{theorem}{section}
\numberwithin{lemma}{section}
\numberwithin{proposition}{section}
\numberwithin{definition}{section}

\begin{center}
{\Large \textbf{Supplementary materials for the paper 
``Main Effect Factor Models in High-Dimensional Matrix Time Series: Identification and Sparsity''}}
\end{center}

\section{Additional Experiments and Details}\label{sec: additional_num}

\subsection{Additional simulation experiments}\label{subsec: additional_sim}

\paragraph{Accuracy of grand mean and main effects}
Table~\ref{tab:sim-main-effect-heavy-tail} examines the accuracy of the grand mean and main effect estimators when the Gaussian noise innovations are replaced by standardized $t_5$ and $t_3$ innovations. We consider representative identifications from both $\cL$ and $\cQ$. Overall, the patterns under Gaussian noise are well preserved as the tails become heavier. The equal-weight identification remains the most accurate among the linear identifications, while the reference-unit identification has errors similar to those under the minimum identification.

The central and tail quantile identifications respond differently to heavier tails. The MSEs under the median identification remain nearly unchanged from Gaussian to $t_3$ noise, including their upper replication-level quantiles. In contrast, the minimum identification shows some deterioration under $t_3$ noise, most clearly for the grand mean. This reflects the greater sensitivity of an extreme quantile to large noise realizations. The results therefore also illustrate the relative robustness of using a central quantile as the identifying benchmark.

\begin{table}[ht!]
\centering
\caption{Accuracy under heavy-tailed noise. Entries report MSE means multiplied by $10^2$, with the 95th percentile of replication-level MSEs in brackets.}
\label{tab:sim-main-effect-heavy-tail}
{%
\begin{tabular}{ccccc}
\toprule
Noise & Identification & $\MSE(\wt\mu)$ & $\MSE(\wt\balpha)$ & $\MSE(\wt\bbeta)$ \\
\midrule
\multirow{4}{*}{Gaussian} & $g(\aa)=a_1$ & $6.439 [10.634]$ & $6.454 [8.899]$ & $6.504 [9.220]$ \\
 & $g(\aa)=\1_d^\top\aa$ & $0.111 [0.172]$ & $3.246 [3.728]$ & $3.257 [3.667]$ \\
 & $g(\aa)=\min\{\aa\}$ & $6.149 [8.145]$ & $6.151 [7.358]$ & $6.185 [7.410]$ \\
 & $g(\aa)=\med\{\aa\}$ & $1.621 [2.043]$ & $4.003 [4.549]$ & $4.008 [4.522]$ \\
\addlinespace[0.35em]
\multirow{4}{*}{$t_5$} & $g(\aa)=a_1$ & $6.574 [11.191]$ & $6.477 [9.494]$ & $6.505 [9.321]$ \\
 & $g(\aa)=\1_d^\top\aa$ & $0.109 [0.173]$ & $3.262 [3.710]$ & $3.262 [3.703]$ \\
 & $g(\aa)=\min\{\aa\}$ & $6.219 [8.238]$ & $6.199 [7.476]$ & $6.190 [7.420]$ \\
 & $g(\aa)=\med\{\aa\}$ & $1.630 [2.075]$ & $4.021 [4.571]$ & $4.018 [4.614]$ \\
\addlinespace[0.35em]
\multirow{4}{*}{$t_3$} & $g(\aa)=a_1$ & $6.558 [11.406]$ & $6.361 [9.164]$ & $6.527 [9.644]$ \\
 & $g(\aa)=\1_d^\top\aa$ & $0.112 [0.180]$ & $3.217 [3.791]$ & $3.218 [3.780]$ \\
 & $g(\aa)=\min\{\aa\}$ & $6.971 [8.435]$ & $6.341 [7.604]$ & $6.336 [7.784]$ \\
 & $g(\aa)=\med\{\aa\}$ & $1.552 [1.980]$ & $3.936 [4.601]$ & $3.946 [4.630]$ \\
\bottomrule
\end{tabular}%
}
\end{table}

\paragraph{Accuracy of common components}
Table~\ref{tab:sim-common-component-accuracy-t3} repeats the common component experiment using standardized $t_3$ noise innovations. The main qualitative patterns are preserved. Increasing $T$ and the cross-sectional dimensions generally improves estimation accuracy, while weak factors lead to larger loading-space errors.

Compared with the Gaussian setting, the deterioration is mild under pervasive factors in most cases. The loading-space errors remain of similar magnitude, and the common component is also estimated accurately in most settings. The effect of heavy tails is more pronounced under weak factors. Both loading-space errors and common-component MSEs increase substantially relative to the Gaussian setting. Increasing $T$ still improves the results in general.

\begin{table}[ht!]
\centering
\caption{Loading-space and common component estimation accuracy under standardized $t_3$ noise innovations, where $\zeta=\zeta_r=\zeta_c$. Entries are Monte Carlo means multiplied by $10^2$, with standard errors in parentheses.}
\label{tab:sim-common-component-accuracy-t3}
\resizebox{\textwidth}{!}{%
\begin{tabular}{ccccccccc}
\toprule
\multirow{2}{*}{$p$} & \multirow{2}{*}{$q$} & \multirow{2}{*}{$\zeta$} & \multicolumn{2}{c}{$\c{D}(\Q_r,\wh \Q_r)$} & \multicolumn{2}{c}{$\c{D}(\Q_c,\wh \Q_c)$} & \multicolumn{2}{c}{$\MSE(\wh \C)$} \\
\cmidrule(lr){4-5}\cmidrule(lr){6-7}\cmidrule(lr){8-9}
& & & $T=100$ & $T=200$ & $T=100$ & $T=200$ & $T=100$ & $T=200$ \\
\midrule
\multirow{2}{*}{30} & \multirow{2}{*}{30} & $0$ & $2.346(0.200)$ & $1.652(0.042)$ & $2.217(0.060)$ & $1.610(0.029)$ & $0.764(0.085)$ & $0.575(0.008)$ \\
 &  & $0.2$ & $26.713(1.052)$ & $21.815(0.928)$ & $26.205(1.049)$ & $21.693(0.958)$ & $5.540(1.021)$ & $3.632(0.468)$ \\
\addlinespace[0.35em]
\multirow{2}{*}{80} & \multirow{2}{*}{80} & $0$ & $1.277(0.196)$ & $0.761(0.007)$ & $1.284(0.195)$ & $0.754(0.007)$ & $1.653(1.521)$ & $0.097(0.001)$ \\
 &  & $0.2$ & $21.170(1.074)$ & $15.785(0.934)$ & $21.196(1.101)$ & $15.648(0.927)$ & $1.782(0.235)$ & $1.466(0.329)$ \\
\addlinespace[0.35em]
\multirow{2}{*}{30} & \multirow{2}{*}{80} & $0$ & $1.406(0.046)$ & $1.290(0.152)$ & $1.831(0.020)$ & $1.379(0.076)$ & $0.316(0.009)$ & $0.417(0.106)$ \\
 &  & $0.2$ & $24.381(1.005)$ & $20.455(0.928)$ & $21.308(0.981)$ & $16.764(0.932)$ & $2.898(0.635)$ & $2.205(0.477)$ \\
\bottomrule
\end{tabular}%
}
\end{table}

\subsection{Additional details on real data applications}

\paragraph{QCEW.}
For our QCEW data studied in the main text, the industries cover agriculture, mining, utilities, construction, manufacturing, wholesale and retail trade, transportation, finance, real estate, professional and administrative services, education, health care, arts and recreation, accommodation and food services, and other services,
all classified according to the North American Industry Classification System. 
The observation for each state--industry pair is computed as the year-over-year log employment growth, i.e.,
\begin{align*}
    X_{t,ij}=100 \, (\log N_{t,ij}-\log N_{t-12,ij}),
\end{align*}
where $N_{t,ij}$ denotes monthly employment.

\paragraph{OECD.}
We process the OECD raw data as follows.
Following the pre-processing adopted in \cite{chen2020constrained} and \cite{Yuetal2022}, logarithmic and differencing transformations are applied: the CPI series are second log-differenced, the interest rates are first-differenced, the real-activity series are first log-differenced, and exports and imports are transformed into growth rates. Each transformed country--indicator series is further standardized over time to remove scale differences across heterogeneous macroeconomic variables.

\subsection{Additional real data application on taxi trips}\label{subsec: additional_real_data}

We illustrate the proposed method for MEFM using the publicly available New York City (NYC) Yellow Taxi Trip Records; see the details in \url{https://www.nyc.gov/site/tlc/about/tlc-trip-record-data.page}. 
Each raw trip record includes pick-up/drop-off timestamps, pick-up/drop-off location codes, trip distance, fare components, payment information, etc. The Taxi and Limousine Commission's map partitions NYC into 265 zones, and we focus on the 69 zones on Manhattan Island which accounts for the vast majority of trip activity.
Let $\X_t\in\b{R}^{69\times 24}$ denote the matrix for calendar day $t$. Its $(i,j)$-th entry counts all trips with drop-off zone $i$ and pick-up time in the $j$-th hourly slot. We analyze the period from 1 January 2019 to 31 December 2022, which spans the dramatic mobility collapse during the 2020 pandemic lockdown and the subsequent recovery. 
To respect the distinct spatial/temporal rhythms of the city, we further split the series into weekday and weekend subsets, containing 1044 and 417 days, respectively, and estimate the model on each subset separately. In this application, we impose the minimum identification
\begin{align*}
    g_\alpha(\aa)=\min\{\aa\}
    \;\text{\ and \ }
    g_\beta(\bb)=\min\{\bb\},
\end{align*}
so that the location and hour main effects are nonnegative and, at each time point, their minimum is normalized to zero. Accordingly, a zero main effect indicates the absence of an additional location- or hour-specific additive effect relative to the estimated baseline, rather than the absence of taxi trips.
The tuning parameters are chosen via the modified $C_p$ criterion described in Section~\ref{subsec: practical_implementation}.

\begin{figure}[htp!]
\centering
\includegraphics[width=.9\columnwidth]{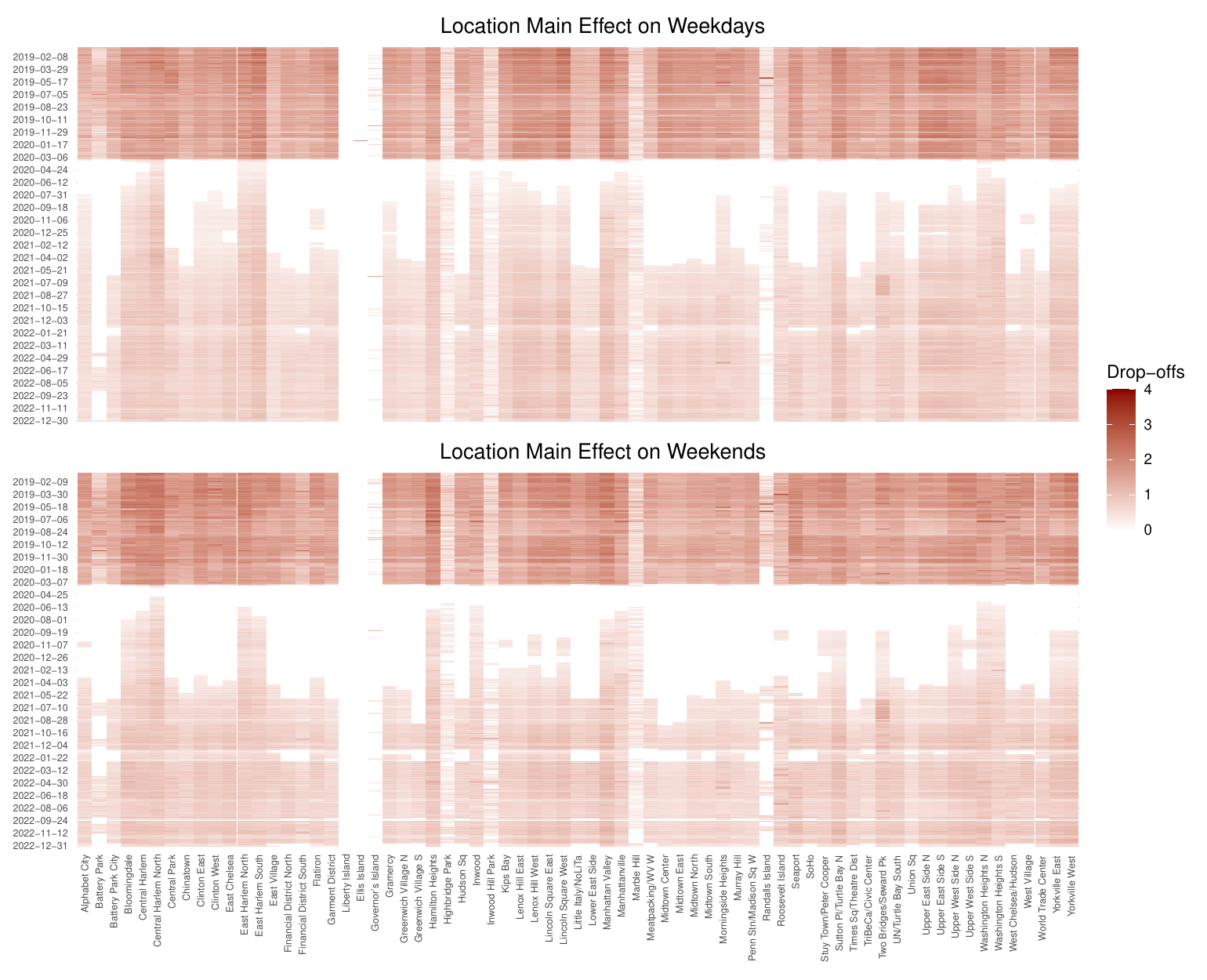}
\caption{Heatmap of the location main effect $\wt\balpha_t^{\cB}$.}
\label{fig: location_ME}
\end{figure}

In Figure~\ref{fig: location_ME}, the location main effects for weekdays and weekends are demonstrated. Firstly, an abrupt shift from dark red to pale white in mid-March 2020 signals the pandemic lockdown, during which the location main effect estimators for almost all Manhattan zones collapse to zero for several months. Secondly, the rebound is asymmetric: weekday activity re-emerges more quickly than weekend activity in most zones, reflecting a faster return to work routines while leisure remains depressed. Thirdly, several locations---Liberty, Ellis, Governors, and Randalls Islands---stay faint almost throughout, yet the model still detects sporadic taxi traffic. On the other hand, zones such as Harlem, Hamilton Heights, Manhattan Valley, and Washington Heights preserve moderate intensities even during the pandemic, underscoring local demand. Finally, neither weekdays nor weekends regain their 2019 pre-pandemic intensity by late 2021--2022, demonstrating the impact of the pandemic on Manhattan mobility and, by extension, urban economic and lifestyle patterns. Notably, the fused penalty is what allows the model to delineate the prolonged ``silent'' pandemic interval, yet it remains sensitive enough to uncover the faint residual taxi activity that persisted throughout the lockdown.

\begin{figure}[htp!]
\centering
\includegraphics[width=.95\columnwidth]{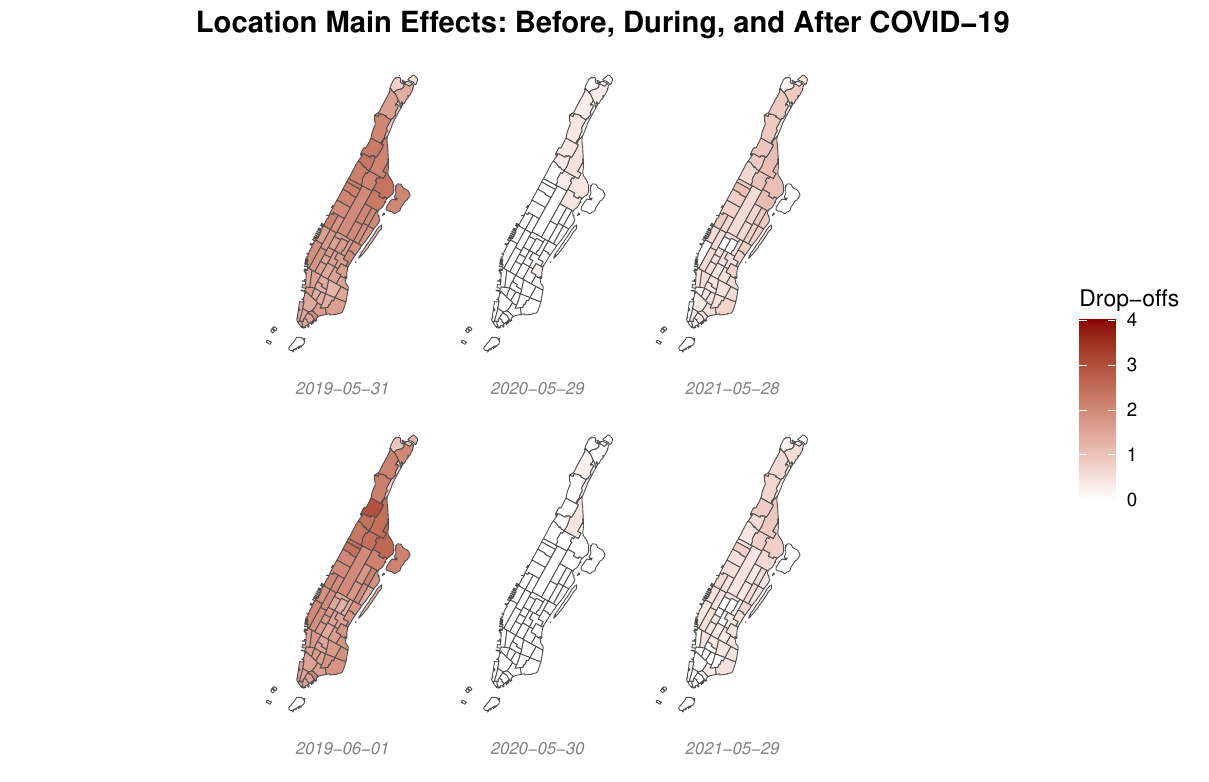}
\caption{Snapshot heatmaps of location main effect. The top rows are weekdays, and the bottom rows are weekends. The color scale matches Figure~\ref{fig: location_ME}.}
\label{fig: location_ME_single_day}
\end{figure}

To further demonstrate the spatial pattern, Figure~\ref{fig: location_ME_single_day} offers a snapshot of how the estimated location main effects evolved across Manhattan. Six representative days are selected to demonstrate the before, during, and after the pandemic period for weekdays and weekends. Before the pandemic, both rows are covered by dark reds, but their centers of taxi rides differ: weekdays peak in business core in Midtown whereas weekends tilt toward leisure areas in the West Village, SoHo and the Lower East Side. However, during lockdown, the island converts almost uniformly pale, and only a few residential and hospital-adjacent blocks in Upper Manhattan retain faint activity during weekday, while the corresponding weekend map is blank, underscoring the sharper collapse of leisure travel. After a year of rehabilitation, demand partially returns, but the asymmetry remains: weekend drop-offs stay muted across several entertainment districts downtown, and weekday rebounds more strongly, especially in residential uptown zones.

\begin{figure}[htp!]
\centering
\includegraphics[width=.9\columnwidth]{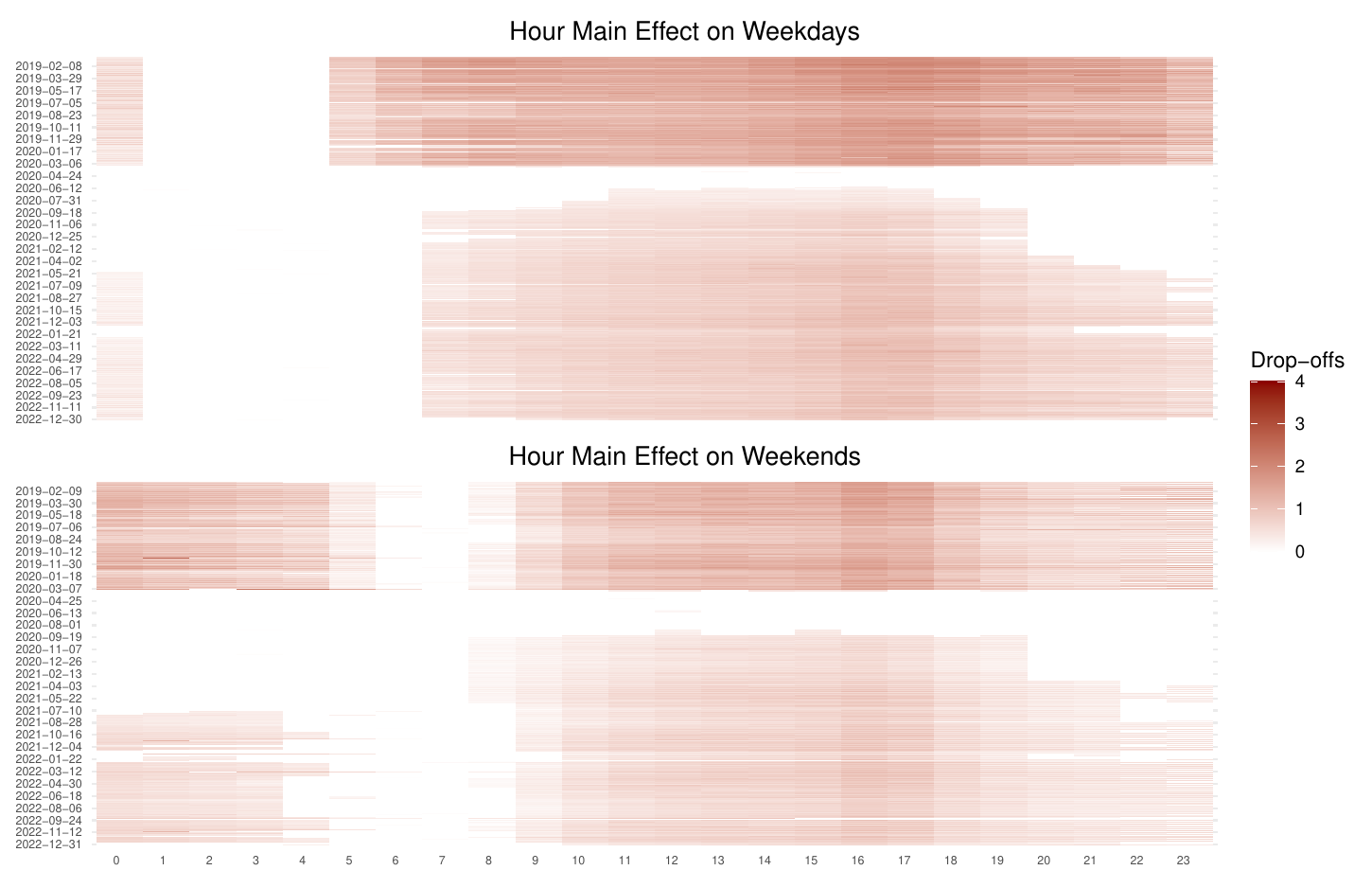}
\caption{Heatmap of the hour main effect $\wt\bbeta_t^{\cB}$.}
\label{fig: hour_ME}
\end{figure}

Furthermore, the hour main effect is displayed in Figure~\ref{fig: hour_ME}. We may observe that weekday mobility peaks in the late-afternoon commute (3pm--6pm) and tapers off by 1am, whereas weekends sustain a pronounced nightlife pulse until about 4am and begin later in the morning. In addition, the lockdown starts in mid-March 2020 and erases nearly all hourly traffic for several weeks. The weekend panel stays blank longer, demonstrating a sharper hit to entertainment travel. Moreover, recovery unfolds asymmetrically, where weekday evening activity re-emerges first (from mid-June 2020), with the signal band sliding from 5pm down to midnight over the next year. Weekend nightlife restarts only in the early autumn of 2020, and does not regain its pre-pandemic 4am endpoint by 2022, indicating a lasting contraction of late-night leisure demand.

\newpage

\section{Technical Proofs}\label{sec: proofs}

\subsection{Proof of theorems}

\begin{proof}[Proof of Theorem~\ref{thm: identification}]
First, consider part~(i) in the statement.
Suppose we have another set of parameters, $\big(\ddot\mu_t, \ddot\balpha_t, \ddot\bbeta_t, \ddot\C_t \big)$ for $t\in[T]$, also satisfying \eqref{eqn: MEFM}. For each $t\in[T]$, right-multiply $\1_q$ on both sides of \eqref{eqn: MEFM}, by Assumption~(IC1)~(i), it holds that
\begin{align*}
    \X_t \1_q &= q\balpha_t^\ast + \1_p (q\mu_t + \1_q^\top \bbeta_t^\ast) + \E_t \1_q \\
    &=
    q\ddot\balpha_t + \1_p (q\ddot\mu_t + \1_q^\top \ddot\bbeta_t) + \E_t \1_q .
\end{align*}
Note that both $(q\mu_t + \1_q^\top \bbeta_t^\ast)$ and $(q\ddot\mu_t + \1_q^\top \ddot\bbeta_t)$ are scalars, we obtain $\balpha_t^\ast - \ddot\balpha_t = c_t\1_p$ for some scalar $c_t$. This implies each row main effect is at least identified up to a shift.
Then by Assumption~(IC1)~(ii), we have
\begin{align*}
    0= g_\alpha(\balpha_t^\ast) = g_\alpha(\ddot\balpha_t + c_t\1_p) ,
\end{align*}
which immediately implies $c_t=0$ by (IC1)~(iii), and thus $\balpha_t^\ast = \ddot\balpha_t$.
For the column main effects, all the above arguments follow analogously by left-multiplying $\1_p^\top$ on both sides of \eqref{eqn: MEFM}, and hence $\bbeta_t^\ast = \ddot\bbeta_t$.
For each $t\in [T]$, by respectively left- and right-multiplying $\1_p^\top$ and $\1_q$ on both sides of \eqref{eqn: MEFM}, we obtain
\begin{align*}
    \1_p^\top \X_t \1_q &= pq\mu_t + q\1_p^\top \balpha_t^\ast + p\1_q^\top \bbeta_t^\ast + \1_p^\top \E_t \1_q \\
    &=
    pq\ddot\mu_t + q\1_p^\top \ddot\balpha_t + p\1_q^\top \ddot\bbeta_t + \1_p^\top \E_t \1_q ,
\end{align*}
so that $\mu_t =\ddot\mu_t$ given all main effects are identified previously. Lastly, the remaining components $\C_t$ and $\ddot\C_t$ trivially coincide for all $t\in[T]$. This completes the ``if'' part in the statement of this theorem.

To show the ``only if'' part, it suffices to show that Assumption~(IC1)~(iii) is necessary to identify the row main effects.
Suppose Assumption~(IC1)~(iii) does not hold. Namely, there exists some vector $\bf{a}$ with $g_\alpha(\bf{a}) =0$ and non-zero constant $d$ such that $g_\alpha(\bf{a} + d\1_p) =0$. This implies both row main effects $\bf{a}$ and $\bf{a} + d\1_p$ fulfill Assumption~(IC1)~(ii), but $d\neq 0$ and hence those two vectors are not the same. The row main effect is thus unidentifiable in this case, indicating that Assumption~(IC1)~(iii) is necessary in identifying the row main effects.
This completes the proof of part~(i) in the theorem.

For part~(ii), let $\X_t = \Acute\A_r \F_t \Acute\A_c^\top +\E_t =: \Acute\C_t + \E_t$. For each $t\in[T]$, define $c_{\alpha,t}$ as the scalar such that 
\[
g_\alpha(q^{-1}\M_p \Acute\C_t \1_q +c_{\alpha,t}\1_p)= 0 ,
\]
with $c_{\alpha,t}$ guaranteed to exist according to the condition in the statement of Theorem~\ref{thm: identification}~(ii). Define analogously $c_{\beta,t}$.
Then we can construct
\begin{align*}
    \A_r &= \M_p \Acute\A_r = (\I_p - p^{-1}\1_p\1_p^\top) \Acute\A_r , 
    \;\;\;
    \balpha_t^\ast = q^{-1}\M_p \Acute\C_t \1_q + \1_p c_{\alpha,t} , \\
    \A_c &= \M_q \Acute\A_c = (\I_q - q^{-1}\1_q\1_q^\top) \Acute\A_c , 
    \;\;\;
    \bbeta_t^\ast = p^{-1}\M_q \Acute\C_t^\top \1_p + \1_q c_{\beta,t} , \\
    \mu_t &= (pq)^{-1}\1_p^\top \Acute\C_t \1_q - c_{\alpha,t} - c_{\beta,t} .
\end{align*}
By Remark~1 in \cite{lam2024matrix}, it holds immediately that
\begin{align*}
    \Acute\C_t &= (pq)^{-1}(\1_p^\top \Acute\C_t \1_q) \1_p \1_q^\top + q^{-1}(\M_p \Acute\C_t \1_q) \1_q^\top + p^{-1}\1_p (\M_q \Acute\C_t^\top \1_p)^\top + (\M_p \Acute\A_r) \F_t (\M_q \Acute\A_c)^\top \\
    &=
    \mu_t \1_p \1_q^\top + \balpha_t^\ast \1_q^\top + \1_p \bbeta_t^{\ast \top} + \A_r\F_t\A_c^\top ,
\end{align*}
as in \eqref{eqn: MEFM}. By the definition of $c_{\alpha,t}$ and $c_{\beta,t}$, it is easy to see $g_\alpha(\balpha_t^\ast)=0$ and $g_\beta(\bbeta_t^\ast)=0$, thus satisfying Assumption~(IC1)~(ii). Assumption~(IC1)~(i) is indeed satisfied by $\M_a\1_a = \0$ and hence
\[
\1_p^\top (\M_p \Acute\A_r) = \0 , \;\;\;
\1_q^\top (\M_q \Acute\A_c) = \0 .
\]
For completeness, Assumption~(IC1)~(iii) holds automatically. This completes the proof of the theorem.
\end{proof}

\begin{proof}[Proof of Theorem~\ref{thm: general_initial_consistency_maineffect}]
First, consider the estimator of row main effects.
Fix any $t\in[T]$.
To begin with, by Lemma~\ref{lemma: correlation_Et_Ft}~(i), we have
\begin{align}\label{eqn: sum_of_rows_bound}
    \b{E}\big[\|q^{-1}\E_t\1_q\|^2\big] = q^{-2}\sum_{i=1}^p\var\big(\sum_{j=1}^q E_{t,ij}\big) = O(pq^{-1}).
\end{align}
As discussed in Section~\ref{subsec: estimation}, using Assumptions~(IC1) and (IC2), it holds that
\begin{align*}
    \balpha_t^\ast 
    &= q^{-1} \X_t \1_q - \1_p (\mu_t + q^{-1} \1_q^\top \bbeta_t^\ast) - q^{-1} \E_t \1_q \\
    &=
    q^{-1} \X_t \1_q - \1_p \cdot \wt{g}_\alpha( q^{-1} \X_t \1_q  -q^{-1} \E_t \1_q ) - q^{-1} \E_t \1_q .
\end{align*}
Together with \eqref{eqn: general_tilde_balpha_t}, we have
\begin{align}
    \wt\balpha_t - \balpha_t^\ast
    =
    \1_p \cdot \left\{ \wt{g}_\alpha( q^{-1} \X_t \1_q  -q^{-1} \E_t \1_q ) - \wt{g}_\alpha( q^{-1} \X_t \1_q) \right\} + q^{-1} \E_t \1_q .
    \label{eqn: alpha_tilde_ast_diff}
\end{align}
This gives
\begin{align*}
    \frac{1}{p}\|\wt\balpha_t - \balpha_t^\ast\|^2 &=
    \frac{1}{p}\left\| \1_p \cdot \left\{ \wt{g}_\alpha( q^{-1} \X_t \1_q  -q^{-1} \E_t \1_q ) - \wt{g}_\alpha( q^{-1} \X_t \1_q) \right\} + q^{-1} \E_t \1_q \right\|^2 \\
    &\leq
    2 \left\{ \wt{g}_\alpha( q^{-1} \X_t \1_q  -q^{-1} \E_t \1_q ) - \wt{g}_\alpha( q^{-1} \X_t \1_q) \right\}^2
    + \frac{2}{p} \left\| q^{-1} \E_t \1_q \right\|^2 \\
    &\leq
    2 L_\alpha^2 \left\| q^{-1} \E_t \1_q \right\|_\ast^2
    + \frac{2}{p} \left\| q^{-1} \E_t \1_q \right\|^2 
    =
    O_P\left( \frac{L_\alpha^2}{q^2} \left\| \E_t \1_q \right\|_\ast^2 
    + \frac{1}{q} \right) .
\end{align*}
where the second inequality used Assumption~(S1) and the notations therein,
and the last equality used \eqref{eqn: sum_of_rows_bound}.
Analogous arguments follow for the estimator of the column main effects.

Next, consider the grand mean estimator at $t\in[T]$.
For this, note that by \eqref{eqn: alpha_tilde_ast_diff} and Assumption~(S1) again, we have
\begin{align}
    \frac{1}{p^2} \big\{ \1_p^\top (\wt\balpha_t - \balpha_t^\ast) \big\}^2 &=
    \frac{1}{p^2} \left[ p \left\{ \wt{g}_\alpha( q^{-1} \X_t \1_q  -q^{-1} \E_t \1_q ) - \wt{g}_\alpha( q^{-1} \X_t \1_q) \right\} + q^{-1} \1_p^\top \E_t \1_q \right]^2 \notag \\
    &\leq
    2 \left\{ \wt{g}_\alpha( q^{-1} \X_t \1_q  -q^{-1} \E_t \1_q ) - \wt{g}_\alpha( q^{-1} \X_t \1_q) \right\}^2
    + \frac{2}{p^2 q^2} (\1_p^\top \E_t\1_q)^2 \notag \\
    &\leq
    2 L_\alpha^2 \left\| q^{-1} \E_t \1_q \right\|_\ast^2
    + \frac{2}{p^2 q^2} (\1_p^\top \E_t\1_q)^2 .
    \label{eqn: mu_tilde_diff_step1}
\end{align}
Similarly, it holds that
\begin{align}
    \frac{1}{q^2} \big\{ \1_q^\top (\wt\bbeta_t - \bbeta_t^\ast) \big\}^2 &\leq
    2 L_\beta^2 \left\| p^{-1} \E_t^\top \1_p \right\|_\ast^2
    + \frac{2}{p^2 q^2} (\1_p^\top \E_t\1_q)^2 .
    \label{eqn: mu_tilde_diff_step2}
\end{align}

Recalling $\mu_t = (pq)^{-1} \1_p^\top \X_t \1_q - p^{-1} \1_p^\top \balpha_t^\ast - q^{-1} \1_q^\top \bbeta_t^\ast - (pq)^{-1} \1_p^\top \E_t \1_q$
and the construction in \eqref{eqn: general_tilde_mu_t}, we have
\begin{align*}
    (\wt\mu_t - \mu_t)^2 &= \big\{ p^{-1} \1_p^\top (\wt\balpha_t -\balpha_t^\ast) + q^{-1} \1_q^\top (\wt\bbeta_t -\bbeta_t^\ast) - (pq)^{-1}\1_p^\top \E_t\1_q \big\}^2 \\
    &\lesssim
    \frac{1}{p^2} \big\{ \1_p^\top (\wt\balpha_t -\balpha_t^\ast) \big\}^2
    + \frac{1}{q^2} \big\{ \1_q^\top (\wt\bbeta_t -\bbeta_t^\ast) \big\}^2
    + \frac{1}{p^2 q^2} (\1_p^\top \E_t\1_q)^2 \\
    &\lesssim
    L_\alpha^2 \left\| q^{-1} \E_t \1_q \right\|_\ast^2
    + L_\beta^2 \left\| p^{-1} \E_t^\top \1_p \right\|_\ast^2
    + \frac{1}{p^2 q^2} (\1_p^\top \E_t\1_q)^2 \\
    &=
    O_P\left( \frac{L_\alpha^2}{q^2} \left\| \E_t \1_q \right\|_\ast^2 
    + \frac{L_\beta^2}{p^2} \left\| \E_t^\top \1_p \right\|_\ast^2 
    + \frac{1}{pq} \right) .
\end{align*}
where the third line used \eqref{eqn: mu_tilde_diff_step1} and \eqref{eqn: mu_tilde_diff_step2},
and the last equality used the fact that 
\[
\big((pq)^{-1}\1_p^\top\E_t\1_q\big)^2 = O_P(p^{-1}q^{-1}) 
\]
from the proof of Theorem~2 in \cite{lam2024matrix}.
This completes the proof for the grand mean estimator and hence finishes this proof.
\end{proof}

\begin{proof}[Proof of Theorem~\ref{thm: particular_initial_consistency_maineffect}]
Before any part, first consider a series of sub-Gaussian random variables $\{X_i\}_{i=1}^{p}$ with zero mean and variance proxy $\sigma^2$. For any $\lambda \geq 0$, by the Jensen's inequality, we have
\begin{align*}
    \exp\{\lambda\b{E} (\max_{i\in[p]} \{X_i\})\} &\leq 
    \b{E}(\exp\{\lambda \max_{i\in[p]} \{X_i\}\}) 
    \leq 
    \b{E}(\max_{i\in[p]} \{\exp\{\lambda X_i\} \}) \\
    &\leq 
    \sum_{i=1}^p \b{E}(e^{\lambda X_i}) 
    \leq pe^{\lambda^2\sigma^2/2} ,
\end{align*}
implying that $\b{E} (\max_{i\in[p]} \{X_i\}) \leq \log(p)/\lambda + \lambda\sigma^2/2 $. 
Choosing $\lambda=\sqrt{2\log(p)}/\sigma$ yields $\b{E} (\max_{i\in[p]} \{X_i\})\leq \sqrt{2\sigma^2\log(p)}$.
Therefore, we can conclude that $\b{E} (\max_{i\in[p]} |X_i|) \leq \sqrt{2\sigma^2\log(2p)}$.
With the above argument, note that for each $i\in[p]$, $X_i = \sum_{j=1}^q E_{t,ij}$ is a sub-Gaussian random variables with variance proxy $\sigma^2 \asymp q$, according to Assumptions~(E1), (E2), and (E3). Subsequently, we have $\b{E}(\max_{i\in[p]}\{|\sum_{j=1}^q E_{t,ij}|\}) = O(\sqrt{q\log(p)})$, which in turn yields
\begin{align}
    \|\E_t\1_q \|_{\infty} 
    =  \max_{i\in[p]}\Bigg( \bigg| \sum_{j=1}^q E_{t,ij} \bigg| \Bigg)
    = O_P\left( \sqrt{q\log(p)} \right) .
    \label{eqn: Et1q_max_rate}
\end{align}

Now consider part~(i).
In this case, note that Assumption~(S1) can be satisfied when $\|\cdot\|_\ast$ is set as the Euclidean norm $\|\cdot\|$.
Then for any shift-varying function in the set $\cL$, by \eqref{eqn: linear_tilde_balpha_t}, we have
\[
|\wt{g}_\alpha(\x_1) - \wt{g}_\alpha(\x_2)|
= \left| (\v^\top\1_p)^{-1} \v^\top(\x_1-\x_2) \right|
\leq 
\frac{\| \v\|}{|\v^\top\1_p|} \cdot \|\x_1-\x_2\| ,
\]
thus indicating the Lipschitz constant as $L_\alpha = \|\v\| /|\v^\top\1_p|$.
Also, by \eqref{eqn: sum_of_rows_bound}, $\| \E_t\1_q \|^2 = O_P(pq)$. Then we may apply Theorem~\ref{thm: general_initial_consistency_maineffect} and obtain
\begin{align}
    \frac{1}{p}\|\wt\balpha_t - \balpha_t^\ast\|^2 &=
    O_P\bigg( \frac{L_\alpha^2}{q^2} \left\| \E_t \1_q \right\|^2 
    + \frac{1}{q} \bigg) 
    =
    O_P\bigg( \frac{p \|\v\|^2}{q (\v^\top\1_p)^2}
    + \frac{1}{q} \bigg) .
    \label{eqn: alpha_rate_cL_Euclidean}
\end{align}
On the other hand, Assumption~(S1) can also be satisfied with $\|\cdot\|_\ast$ set as the $L_{\infty}$ norm $\|\cdot\|_{\infty}$, in which case we have $L_\alpha = \|\v\|_1 /|\v^\top\1_p|$.
Together with \eqref{eqn: Et1q_max_rate} and Theorem~\ref{thm: general_initial_consistency_maineffect}, we obtain
\begin{align}
    \frac{1}{p}\|\wt\balpha_t - \balpha_t^\ast\|^2 &=
    O_P\bigg( \frac{L_\alpha^2}{q^2} \left\| \E_t \1_q \right\|_{\infty}^2 
    + \frac{1}{q} \bigg) 
    =
    O_P\bigg( \frac{\log(p) \|\v\|_1^2}{q (\v^\top\1_p)^2}
    + \frac{1}{q} \bigg) .
    \label{eqn: alpha_rate_cL_max}
\end{align}
Combining \eqref{eqn: alpha_rate_cL_Euclidean} and \eqref{eqn: alpha_rate_cL_max}, we have the desired for the row main effect estimators.
The results for the column main effect estimators and the grand mean estimators follow similarly.
Now when $\v=\1$, it is direct that the Lipschitz constant (for the row main effect estimators) has value $L_\alpha = \|\v\| /|\v^\top\1_p| = p^{-1/2}$, and similarly $L_\beta = q^{-1/2}$. It remains to show that the grand mean estimator in this case would not carry these two rates of $1/p$ and $1/q$.
To see this, note that by \eqref{eqn: linear_tilde_balpha_t},
\begin{align*}
    \1_p^\top \wt\balpha_t &= 
    q^{-1} \1_p^\top \X_t \1_q - \frac{p}{\v^\top \1_p} (q^{-1} \v^\top \X_t \1_q - c) ,
\end{align*}
which has value $c$ if $\v=\1_p$. Also, recall that Assumption~(IC1) requires $\1_p^\top \balpha_t^\ast = c$ in such case.
Then from the decomposition of $(\wt\mu_t - \mu_t)^2$ in the proof of Theorem~\ref{thm: general_initial_consistency_maineffect}, it is direct that $(\wt\mu_t - \mu_t)^2 = O_P\{ 1/(pq)\}$.
This completes the proof for part~(i) of the theorem.

For part~(ii), recall that
\begin{align*}
    |\wt{g}_\alpha(\x_1) - \wt{g}_\alpha(\x_2)|
    = |q_u(\x_1) - q_u(\x_2)|
    \leq 
    \|\x_1-\x_2\|_{\infty} 
    \leq
    \|\x_1-\x_2\| .
\end{align*}
Thus in Assumption~(S1), the Lipschitz constant $L_\alpha = 1$ when $\|\cdot\|_\ast$ is set as either the Euclidean norm or the $L_{\infty}$ norm.
Since $\left\| \E_t \1_q \right\|^2$ has order larger than that of $\left\| \E_t \1_q \right\|_{\infty}^2$, we use the $L_{\infty}$ norm and apply Theorem~\ref{thm: general_initial_consistency_maineffect}. This gives
\begin{align*}
    \frac{1}{p}\|\wt\balpha_t - \balpha_t^\ast\|^2 &=
    O_P\bigg( \frac{L_\alpha^2}{q^2} \left\| \E_t \1_q \right\|_{\infty}^2 
    + \frac{1}{q} \bigg) 
    =
    O_P\bigg( \frac{\log(p)}{q}
    + \frac{1}{q} \bigg) 
    =
    O_P\bigg( \frac{\log(p)}{q} \bigg) ,
\end{align*}
as desired, which concludes the proof of this theorem.
\end{proof}

\begin{proof}[Proof of Theorem~\ref{thm: general_initial_consistency}]
First, we show the consistency of the factor loading estimators. 
It suffices to consider the row loading estimator, as the arguments for the column loading estimator are analogous.
Recall the notation $\M_a=\I_a - a^{-1}\1_a\1_a^\top$, then from \eqref{eqn: tilde_L_t}, we have
\begin{align*}
    \wt\L_t &= \X_t - (pq)^{-1}\1_p^\top\X_t\1_q\1_p\1_q^\top + p^{-1}\1_p^\top\wt\balpha_t\1_p\1_q^\top + q^{-1}\1_q^\top\wt\bbeta_t\1_p\1_q^\top - \wt\balpha_t\1_q^\top - \1_p\wt\bbeta_t^\top \\
    &=
    \M_p\X_t(\1_q\1_q^\top/q) + \X_t\M_q - \M_p \wt\balpha_t \1_q^\top - \1_p\wt\bbeta_t^\top \M_q \\
    &=
    \M_p\X_t(\1_q\1_q^\top/q) + \X_t\M_q - \M_p\X_t(\1_q\1_q^\top/q) - (\1_p\1_p^\top/p) \X_t\M_q
    = \M_p \X_t \M_q ,
\end{align*}
where the second last equality used the definitions of $\wt\balpha_t$ and $\wt\bbeta_t$ in \eqref{eqn: general_tilde_balpha_t} and \eqref{eqn: general_tilde_bbeta_t} respectively, together with the key that $\M_a \1_a = \0$. Hence, 
\begin{align}
    \wt\L_t\wt\L_t^\top  &= \M_p\X_t\M_q\M_q^\top\X_t^\top\M_p^\top 
    = \M_p\X_t\M_q\X_t^\top\M_p \notag \\
    &= \M_p(\C_t + \E_t)\M_q(\C_t^\top + \E_t^\top)\M_p
    = \C_t\C_t^\top + \C_t\E_t^\top\M_p + \M_p\E_t\C_t^\top + \M_p\E_t\M_q\E_t^\top\M_p \notag \\
    &= \C_t\C_t^\top + \C_t\E_t^\top + \E_t\C_t^\top + \E_t\E_t^\top + (pq)^{-1}\E_t\1_q\1_q^\top\E_t^\top\1_p\1_p^\top + (pq)^{-1}\1_p\1_p^\top\E_t\1_q\1_q^\top\E_t^\top \notag \\ 
    &\hspace{12pt}
    - p^{-1}\C_t\E_t^\top\1_p\1_p^\top - p^{-1}\1_p\1_p^\top\E_t\C_t^\top - p^{-1}\E_t\E_t^\top\1_p\1_p^\top - p^{-1}\1_p\1_p^\top\E_t\E_t^\top - q^{-1}\E_t\1_q\1_q^\top\E_t^\top \notag \\
    &\hspace{12pt}
    + p^{-2}\1_p\1_p^\top\E_t\E_t^\top\1_p\1_p^\top - (pq)^{-1}p^{-1}\1_p\1_p^\top\E_t\1_q\1_q^\top\E_t^\top\1_p\1_p^\top .
    \label{eqn: sample_covariance_row}
\end{align}
To ease notations, define 
\begin{align*}
    \R_{r,t} &:= \C_t\E_t^\top + \E_t\C_t^\top + \E_t\E_t^\top + (pq)^{-1}\E_t\1_q\1_q^\top\E_t^\top\1_p\1_p^\top + (pq)^{-1}\1_p\1_p^\top\E_t\1_q\1_q^\top\E_t^\top \notag \\ 
    &\quad
    - p^{-1}\C_t\E_t^\top\1_p\1_p^\top - p^{-1}\1_p\1_p^\top\E_t\C_t^\top - p^{-1}\E_t\E_t^\top\1_p\1_p^\top - p^{-1}\1_p\1_p^\top\E_t\E_t^\top - q^{-1}\E_t\1_q\1_q^\top\E_t^\top \notag \\
    &\quad
    + p^{-2}\1_p\1_p^\top\E_t\E_t^\top\1_p\1_p^\top - (pq)^{-1}p^{-1}\1_p\1_p^\top\E_t\1_q\1_q^\top\E_t^\top\1_p\1_p^\top , 
\end{align*}
so that from \eqref{eqn: sample_covariance_row},
we can write 
$\wt\L_t\wt\L_t^{\top} = \C_t\C_t^\top + \R_{r,t}$.
In particular, from Lemma~2 of \cite{lam2024matrix}, we have
\begin{align}
    \bigg\| \sum_{t=1}^T \R_{r,t} \bigg\|_F^2 = O_P(Tp^2q + T^2pq^2) .
    \label{eqn: R_rt_rate}
\end{align}
Next, we show
\begin{align}
\big\| \wh\Q_r - \Q_r\H_r^\top \big\|_F^2 &= 
O_P\Big( T^{-1}p^{2(1-\delta_{r,k_r})}q^{1-2\delta_{c,1}} + p^{1-2\delta_{r,k_r}}q^{2(1-\delta_{c,1})} \Big) , 
\label{eqn: row_loading_initial_rate}
\end{align}
where $\H_r$ satisfies that 
$\norm{\H_r\Q_r^\top\Q_r\H_r^\top-\I_{k_r}}_F^2=o_P(1)$, which, since $\Q_r^\top\Q_r\to\bSigma_{A,r}$ and $\H_r=O_P(1)$, implies $\|\H_r\bSigma_{A,r}\H_r^\top-\I_{k_r}\|_F^2=o_P(1)$ and the asymptotic invertibility of $\H_r$.
This would follow similar arguments in the proof of Theorem~2 in \cite{lam2024matrix}, which we delineate as follows for completeness.
By construction, it holds that 
$\wh\Q_r \wh\D_r = T^{-1} \sum_{t=1}^T \wt\L_t \wt\L_t^\top \wh\Q_r$,
and hence
for any $j\in[p]$:
\begin{align}
    \wh\Q_{r,j\cdot} - \H_r\Q_{r,j\cdot} =
    \frac{1}{T} \wh\D_r^{-1} \sum_{i=1}^p
    \wh\Q_{r,i\cdot} \sum_{t=1}^T (\R_{r,t})_{ij} .
    \label{eqn: row_loading_est_decomp}
\end{align}
Hence by sub-multiplicativity of the Frobenius norm, we have
\begin{align*}
    \big\| \wh\Q_r - \Q_r\H_r^\top \big\|_F^2 
    &= 
    \sum_{j=1}^p \big\|\wh\Q_{r,j\cdot} - \H_r\Q_{r,j\cdot} \big\|^2 
    = \sum_{j=1}^p \bigg\| \frac{1}{T} \wh\D_r^{-1} \wh\Q_r^\top \Big(\sum_{t=1}^T \R_{r,t}\Big)_{\cdot j} \bigg\|^2 \\
    &\leq
    \frac{1}{T^2} \big\|\wh\D_r^{-1} \big\|_F^2 \cdot
    \big\|\wh\Q_r \big\|_F^2 \cdot \bigg\| \sum_{t=1}^T \R_{r,t} \bigg\|_F^2 .
\end{align*}
Combining the above with \eqref{eqn: R_rt_rate} and that $\|\wh\D_r^{-1} \|_F = O_P\big( p^{-\delta_{r,k_r}} q^{-\delta_{c,1}} \big)$ from Lemma~3 of \cite{lam2024matrix},
we may conclude \eqref{eqn: row_loading_initial_rate}.
We further improve this rate of convergence in the following.
To this end, note that for any $j\in[p]$, \eqref{eqn: row_loading_est_decomp} gives 
\begin{align*}
    \wh\Q_{r,j\cdot} - \H_r\Q_{r,j\cdot} &=
    \frac{1}{T} \wh\D_r^{-1} \sum_{i=1}^p
    (\wh\Q_{r,i\cdot} - \H_r \Q_{r,i\cdot}) \sum_{t=1}^T (\R_{r,t})_{ij}
    + \frac{1}{T} \wh\D_r^{-1} \sum_{i=1}^p
    \H_r\Q_{r,i\cdot} \sum_{t=1}^T (\R_{r,t})_{ij} .
\end{align*}
Thus similar to the above, we can write
\begin{align*}
    &\quad
    \big\| \wh\Q_r - \Q_r\H_r^\top \big\|_F^2 
    = 
    \sum_{j=1}^p \big\|\wh\Q_{r,j\cdot} - \H_r\Q_{r,j\cdot} \big\|^2  \\
    &\leq
    \sum_{j=1}^p \bigg\| \frac{1}{T} \wh\D_r^{-1} (\wh\Q_r - \Q_r\H_r^\top)^\top \Big(\sum_{t=1}^T \R_{r,t}\Big)_{\cdot j} \bigg\|^2
    +
    \sum_{j=1}^p \bigg\| \frac{1}{T} \wh\D_r^{-1} \H_r \Q_r^\top \Big(\sum_{t=1}^T \R_{r,t}\Big)_{\cdot j} \bigg\|^2 \\
    &\leq
    \frac{1}{T^2} \big\|\wh\D_r^{-1} \big\|_F^2 \cdot
    \big\| \wh\Q_r - \Q_r\H_r^\top \big\|_F^2 \cdot \bigg\| \sum_{t=1}^T \R_{r,t} \bigg\|_F^2
    + \frac{1}{T^2} \big\|\wh\D_r^{-1} \big\|_F^2 \cdot \|\H_r\|_F^2 \cdot
    \bigg\| \sum_{t=1}^T \Q_r^\top \R_{r,t} \bigg\|_F^2 \\
    &=
    O_P\Big(
    T^{-1}p^{2(1-\delta_{r,k_r})}q^{1-2\delta_{c,1}} 
    + p^{1-3\delta_{r,k_r}} q^{2-2\delta_{c,1}}
    \Big) ,
\end{align*}
where the last equality used the rate on $\|\wh\D_r^{-1} \|_F$, \eqref{eqn: R_rt_rate}, \eqref{eqn: row_loading_initial_rate}, and the arguments in the proof of Theorem~5 of \cite{lam2024matrix}.
This concludes the proof for the factor loading estimators.
Finally, the proofs for the rates of the error of estimated factor series and individual common components are similar to that of Theorem~3 in \cite{lam2024matrix}, which completes the proof of this theorem.
\end{proof}

\begin{proof}[Proof of Theorem~\ref{thm: oracle}]
It is sufficient to show the results for the row main effects, as the proof for the column main effects follows similarly. We first show the block consistency or equivalently, sign consistency of the DAFL estimators according to \eqref{eqn: DAFL} in estimating the row main effects.
Define $\wt\balpha_{\cdot,i} = (\wt\alpha_{1,i}, \dots, \wt\alpha_{T,i})^\top$ and $\wh\balpha_{\cdot,i} = (\wh\alpha_{1,i}, \dots, \wh\alpha_{T,i})^\top$, for the ease of notation. By the KKT condition, $\wh\balpha_{\cdot,i}$ is a solution minimizing the loss function in \eqref{eqn: DAFL} if and only if there exists a subgradient
\begin{equation}
\label{eqn: subgradient_h}
\begin{split}
        \h &= \partial \Big(\sum_{t=2}^T u_{\alpha, t, i} \big|\wh\alpha_{t, i} - \wh\alpha_{t-1, i}\big| \Big)
    = \left\{
    \h \in \b{R}^{T} \text{ such that} \right. \\
    &\hspace{15pt}
    \left. h_1: \left\{
    \begin{array}{ll}
	h_1 = -u_{\alpha, 2, i} \, \text{sign}(\wh\alpha_{2, i} - \wh\alpha_{1, i}), & \hbox{$\wh\alpha_{1, i} \neq \wh\alpha_{2, i}$;} \\
	|h_1| \leq u_{\alpha, 2, i}, & \hbox{otherwise.}
    \end{array}
    \right.; \right. \\
    &\hspace{15pt}
    \left. h_T: \left\{
    \begin{array}{ll}
	h_T = u_{\alpha, T, i} \, \text{sign}(\wh\alpha_{T, i} - \wh\alpha_{T-1, i}), & \hbox{$\wh\alpha_{T, i} \neq \wh\alpha_{T-1, i}$;} \\
	|h_T| \leq u_{\alpha, T, i}, & \hbox{otherwise.}
    \end{array}
    \right.; \right. \\
    &\hspace{15pt}
    \text{for } t=2,\dots,T-1, \\
    &\hspace{15pt}
    \left. h_t: \left\{
    \begin{array}{ll}
	h_t = u_{\alpha, t, i} \, \text{sign}(\wh\alpha_{t, i} - \wh\alpha_{t-1, i}) - u_{\alpha, t+1, i} \, \text{sign}(\wh\alpha_{t+1, i} - \wh\alpha_{t, i}), & \hbox{$\wh\alpha_{t-1, i} \neq \wh\alpha_{t, i} \neq \wh\alpha_{t+1, i}$;} \\
    \big|h_t + u_{\alpha, t+1, i} \, \text{sign}(\wh\alpha_{t+1, i} - \wh\alpha_{t, i}) \big| \leq u_{\alpha, t, i} , & \hbox{$\wh\alpha_{t-1, i} = \wh\alpha_{t, i} \neq \wh\alpha_{t+1, i}$;} \\
    \big|h_t - u_{\alpha, t, i} \, \text{sign}(\wh\alpha_{t, i} - \wh\alpha_{t-1, i}) \big| \leq u_{\alpha, t+1, i} , & \hbox{$\wh\alpha_{t-1, i} \neq \wh\alpha_{t, i} = \wh\alpha_{t+1, i}$;} \\
	\big|h_t\big| \leq u_{\alpha, t, i} + u_{\alpha, t+1, i}, & \hbox{otherwise.}
    \end{array}
    \right.
    \right\} ,
\end{split}
\end{equation}
and also a subgradient
\begin{equation}
\label{eqn: subgradient_g}
\begin{split}
    \g = \partial \Big(\sum_{t=1}^{T} \gamma_{\alpha, t, i} \big|\wh\alpha_{t, i}\big| \Big)
    = \left\{
	\g \in \b{R}^{T} : \left\{
	\begin{array}{ll}
		g_t = \gamma_{\alpha, t, i} \, \text{sign}(\wh\alpha_{t,i}), & \hbox{$\wh\alpha_{t,i} \neq 0$;} \\
		|g_t| \leq \gamma_{\alpha, t, i}, & \hbox{otherwise.}
	\end{array}
	\right.
	\right\} ,
\end{split}
\end{equation}
such that differentiating $L(\wh\balpha_{\cdot,i}) \equiv L(\wh\alpha_{1,i}, \dots, \wh\alpha_{T,i})$ with respect to $\wh\balpha_{\cdot,i}$, we have
\begin{equation}
\label{eqn: oracle_KKT}
\wh\balpha_{\cdot,i} -\wt\balpha_{\cdot,i} = -\partial \Big(\lambda_\alpha \sum_{t=2}^T u_{\alpha, t, i}\big|\wh\alpha_{t, i} - \wh\alpha_{t-1, i}\big| + \lambda_\alpha \sum_{t=1}^{T}\gamma_{\alpha, t, i}\big|\wh\alpha_{t, i}\big| \Big)
= -\lambda_\alpha \h - \lambda_\alpha \g .
\end{equation}

Without loss of generality, we can always partition each sparse block into $\cS_{\alpha,i} = \overline\cS_{\alpha,i} \cup \cS_{\alpha,i}^\circ$ , where $\overline\cS_{\alpha,i}$ is the periods when the main effect is piecewise zero and $\cS_{\alpha,i}^\circ$ is the periods when the main effect is distinct zero. Formally, $\overline\cS_{\alpha,i}$ is the largest set such that for all $t\in \overline\cS_{\alpha,i}$: (i) $\alpha_{t,i}^\ast=0$; (ii) $\{t-1, t+1\} \cap \overline\cS_{\alpha,i} \neq \emptyset$.

Due to the intricacy in \eqref{eqn: subgradient_h}, we define the interior of $\overline\cS_{\alpha,i}$ as
\[
\cS_{\alpha,i}^\ast = \big\{t\in \overline\cS_{\alpha,i}: \text{ either }  t\in\{1,T\} \text{ or } \{t-1,t+1\} \subseteq \overline\cS_{\alpha,i} \big\} ,
\]
that is, both (or only one neighbor when $t=1,T$) neighboring timestamps are still in the sparse block.

To have sign consistency, we first show consistency in recovering the sparse block. \eqref{eqn: oracle_KKT} implies that uniformly over $i\in[p]$ and $t\in \cS_{\alpha,i}$,
\begin{align}
    \wt\alpha_{t,i} = \lambda_\alpha h_t + \lambda_\alpha g_t ,
    \label{eqn: zero_consistency_KKT}
\end{align}
where $|g_t| \leq \gamma_{\alpha, t, i}$ according to \eqref{eqn: subgradient_g}. Firstly consider $t\in \overline\cS_{\alpha,i}$, we discuss this in three cases to ease the discussion, where Cases~(a) and (c) correspond to consistency in recovering the interior of the sparse block, Case~(b) corresponds to the margin, $\overline\cS_{\alpha,i} \setminus \cS_{\alpha,i}^*$, of the sparse block:
\begin{itemize}
    \item [(a)] $t\in \cS_{\alpha,i}^\ast$ and $t\notin \{1,T\}$, i.e., both the previous and the next timestamps are also in the sparse block;
    \item [(b)] $t\in \overline\cS_{\alpha,i} \setminus \cS_{\alpha,i}^*$, recall that it means either the previous or the next timestamp is in the sparse block (but not both), i.e., $\{t-1, t\}\subseteq \cS_{\alpha,i}$ or $\{t, t+1\}\subseteq \cS_{\alpha,i}$; the definition also implies $t\neq 1$ and $t\neq T$, which facilitates us to better discuss the subgradient $\h$;
    \item [(c)] $t\in \cS_{\alpha,i}^\ast$ and $t\in \{1,T\}$.
\end{itemize}

Consider Case (a) first, then the subgradient $\h$ in \eqref{eqn: subgradient_h} can only take values $\big|h_t\big| \leq u_{\alpha, t, i} + u_{\alpha, t+1, i}$. Hence for \eqref{eqn: zero_consistency_KKT}, we need to show that (uniformly over $i\in[p]$ and $t\in \overline\cS_{\alpha,i}$),
\begin{align*}
    |\wt\alpha_{t,i}| &= \big| \lambda_\alpha h_t + \lambda_\alpha g_t \big|
    \leq
    \lambda_\alpha |h_t| + \lambda_\alpha |g_t|
    \leq
    \lambda_\alpha (u_{\alpha, t, i} + u_{\alpha, t+1, i}) + \lambda_\alpha \gamma_{\alpha, t, i} .
\end{align*}
The right hand side above is lower bounded by $\lambda_\alpha /\wt\alpha_{t,i}$, according to the definition of $u_{\alpha, t, i}$ and $\gamma_{\alpha, t, i}$. Hence it suffices to show
\begin{align*}
    \max_{i\in[p]} \max_{t\in\overline\cS_{\alpha,i}: t-1\in \overline\cS_{\alpha,i}} \big\{ \wt\alpha_{t,i}^2 \big\} =
    \max_{i\in[p]} \|(\wt\balpha_{\cdot,i})_{\cS_{\alpha,i}}\|_{\infty}^2
    = o_P( \lambda_\alpha ),
\end{align*}
where the second equality indeed holds true from Assumption~(R2) and the first result of Lemma~\ref{lemma: alpha_tilde_max}.

For Case~(b), the subgradient value for $h_t$ satisfies
\begin{align*}
    \text{either} \;\;\;
    \big|h_t + u_{\alpha, t+1, i} \big| \leq u_{\alpha, t, i} \;\;\; \text{with} \;\;\;
    0= \wh\alpha_{t-1, i} = \wh\alpha_{t,i} \neq \wh\alpha_{t+1, i} & \\
    \text{or} \;\;\;
    \big|h_t + u_{\alpha, t, i} \big| \leq u_{\alpha, t+1, i} \;\;\; \text{with} \;\;\;
    \wh\alpha_{t-1, i} \neq \wh\alpha_{t,i} = \wh\alpha_{t+1, i} =0 & .
\end{align*}
We only consider the first scenario above as the second follows a similar proof. Notice that $\wt\alpha_{t,i}\geq 0$, so with \eqref{eqn: zero_consistency_KKT}, it suffices to show that (uniformly over $i\in[p]$ and $t\in\overline\cS_{\alpha,i}$ with $t+1\in \cB_{\alpha,i}$),
\begin{equation}
\label{eqn: zero_consistency_(b)}
\left\{
\begin{array}{ll}
	& -\lambda_\alpha u_{\alpha, t+1, i} - \lambda_\alpha u_{\alpha, t, i} - \lambda_\alpha \gamma_{\alpha,t,i} \leq 0 ,\\
	& \wt\alpha_{t,i} \leq -\lambda_\alpha u_{\alpha, t+1, i} + \lambda_\alpha u_{\alpha, t, i} + \lambda_\alpha \gamma_{\alpha,t,i} ,
\end{array}
\right.
\end{equation}
where the first inequality is immediate by noticing all $u_{\alpha, t+1, i}$, $u_{\alpha, t, i}$, and $\gamma_{\alpha,t,i}$ are non-negative. For the second inequality, since $\wt\alpha_{t,i}$ is dominated by $\lambda_\alpha u_{\alpha, t, i} +\lambda_\alpha \gamma_{\alpha,t,i}$ for $t\in \overline\cS_{\alpha,i}$ using the argument in Case~(a), it remains to show $u_{\alpha, t+1, i}$ is stochastically dominated by $\gamma_{\alpha,t,i}$, for which, noting that $\overline\cS_{\alpha,i} \subseteq \cS_{\alpha,i}$, it suffices to show
\begin{equation}
\label{eqn: sign_consistency_case(b)}
\max_{i\in[p]} \max_{t\in \cS_{\alpha,i}: t+1\in \cB_{\alpha,i}} \{ \wt\alpha_{t,i} \} = o_P\big( \min_{i\in[p]} \min_{t\in \cS_{\alpha,i}: t+1\in \cB_{\alpha,i}} \{\wt\alpha_{t+1,i}\} \big).
\end{equation}
In detail, we have
\begin{equation*}
\begin{split}
\left\{
\begin{array}{ll}
	&\hspace{11pt}
    \max_{i\in[p]} \max_{t\in \cS_{\alpha,i}: t+1\in \cB_{\alpha,i}} \{ \wt\alpha_{t,i} \} \leq \max_{i\in[p]} \|(\wt\balpha_{\cdot,i})_{\cS_{\alpha,i}}\|_{\infty} ,\\
	&\hspace{11pt}
    \min_{i\in[p]} \min_{t\in \cS_{\alpha,i}: t+1\in \cB_{\alpha,i}} \{\wt\alpha_{t+1,i}\} \\
    & \geq
    \min_{i\in[p]} \min_{t\in \cS_{\alpha,i}: t+1\in \cB_{\alpha,i}} \{ \alpha_{t+1,i}^\ast\} - \max_{i\in[p]} \max_{t\in \cS_{\alpha,i}: t+1\in \cB_{\alpha,i}} |\wt\alpha_{t+1,i} -\alpha_{t+1,i}^\ast| ,
\end{array}
\right.
\end{split}
\end{equation*}
which, together with Assumption~(R2) and Lemma~\ref{lemma: alpha_tilde_max}, shows \eqref{eqn: sign_consistency_case(b)}.

Lastly, consider Case~(c). If $t=1$, then we have $|h_1|\leq u_{\alpha,2,i}$; otherwise if $t=T$, then $|h_T|\leq u_{\alpha,T,i}$. Suppose $t=1$, with \eqref{eqn: zero_consistency_KKT}, it is required that
\begin{align*}
    |\wt\alpha_{1,i}| &= \big| \lambda_\alpha h_1 + \lambda_\alpha g_1 \big|
    \leq
    \lambda_\alpha |h_1| + \lambda_\alpha |g_1|
    \leq
    \lambda_\alpha u_{\alpha, 2, i} + \lambda_\alpha \gamma_{\alpha, 1, i} .
\end{align*}
Similar to Case~(a), it is sufficient to show $\max_{i\in[p]} \{\wt\alpha_{1,i}^2\} \leq \max_{i\in[p]} \|(\wt\balpha_{\cdot,i})_{\cS_{\alpha,i}} \|_{\infty}^2 = o_P( \lambda_\alpha)$, which is true by exactly the same argument in Case~(a). The argument for the scenario $t=T$ is similar and hence omitted here. This completes the proof of consistency in recovering the sparse blocks.

It remains to consider $t\in \cS_{\alpha,i}^\circ$ for consistency in the sparse block. Using \eqref{eqn: subgradient_h}, we require uniformly over $i\in[p]$ and $t\in\cS_{\alpha,i}^\circ$ that
\[
\big| \wt\alpha_{t,i} + u_{\alpha, t, i} + u_{\alpha, t+1, i} \big| = \wt\alpha_{t,i} + \lambda_\alpha u_{\alpha, t, i} + \lambda_\alpha u_{\alpha, t+1, i} \leq \lambda_\alpha \gamma_{\alpha,t,i} ,
\]
which holds true if we can show
\[
\left\{
\begin{array}{ll}
	& \max_{i\in[p]} \|(\wt\balpha_{\cdot,i})_{\cS_{\alpha,i}^\circ}\|_{\infty}^2 \leq \max_{i\in[p]} \|(\wt\balpha_{\cdot,i})_{\cS_{\alpha,i}}\|_{\infty}^2 = o_P(\lambda_\alpha) ,\\
	& \max_{i\in[p]} \max_{t\in \cS_{\alpha,i}^\circ :t+1\in \cB_{\alpha,i}} \{ \wt\alpha_{t,i}\} = o_P\big( \min_{i\in[p]} \min_{t\in \cS_{\alpha,i}^\circ :t+1\in \cB_{\alpha,i}} \{\wt\alpha_{t+1,i}\} \big) .
\end{array}
\right.
\]
The first equality is direct from Lemma~\ref{lemma: alpha_tilde_max} and Assumption~(R2), while the second equality from \eqref{eqn: sign_consistency_case(b)}.

For sign consistency in $\cB_{\alpha,i}$, from \eqref{eqn: oracle_KKT}, we require with probability approaching 1 that
\[
\0< (\wh\balpha_{\cdot,i})_{\cB_{\alpha,i}} =(\wt\balpha_{\cdot,i})_{\cB_{\alpha,i}} -\lambda_\alpha (\h)_{\cB_{\alpha,i}} - \lambda_\alpha (\g)_{\cB_{\alpha,i}} ,
\]
which holds true if both $\lambda_\alpha (\h)_{\cB_{\alpha,i}}$ and $\lambda_\alpha (\g)_{\cB_{\alpha,i}}$ are stochastically dominated by $(\wt\balpha_{\cdot,i})_{\cB_{\alpha,i}}$. By \eqref{eqn: subgradient_h} and \eqref{eqn: subgradient_g}, it suffices to show $\lambda_\alpha$ is stochastically dominated by $\min_{i\in[p]} \min_{t\in \cB_{\alpha,i}} \{\wt\alpha_{t,i}^2\}$. This is direct by combining Assumption~(R2), Lemma~\ref{lemma: alpha_tilde_max}, and
\begin{equation}
\label{eqn: minmin_tilde_alpha}
\min_{i\in[p]} \min_{t\in \cB_{\alpha,i}} \{\wt\alpha_{t,i}\} \geq
\min_{i\in[p]} \min_{t\in \cB_{\alpha,i}} \{\alpha_{t,i}^\ast\} - \max_{i\in[p]} \max_{t\in \cB_{\alpha,i}} |\wt\alpha_{t,i} - \alpha_{t,i}^\ast|.
\end{equation}
This ends the proof of block consistency. It remains to show the consistency for the DAFL estimators. 

To this end, note that for any $i\in[p]$, $t\in \cB_{\alpha,i}$, 
\[
\wh\alpha_{t,i} - \alpha_{t,i}^\ast = 
(\wh\alpha_{t,i} - \wt\alpha_{t,i}) + 
(\wt\alpha_{t,i} - \alpha_{t,i}^\ast)
=:\cI_1 + \cI_2.
\]
For $\cI_1$, combining the KKT condition \eqref{eqn: oracle_KKT} and the subgradients from \eqref{eqn: subgradient_h} and \eqref{eqn: subgradient_g}, we have 
\begin{equation}
\label{eqn: DAFL_rate_cI_1}
\begin{split}
    |\cI_1| &\leq \lambda_\alpha |u_{\alpha,t,i} + u_{\alpha,t+1,i} + \gamma_{\alpha,t,i}|
    \leq \frac{\lambda_\alpha}{\min_{i\in[p]} \min_{t\in\cB_{\alpha,i}} \{\wt\alpha_{t,i}\}} \\
    &\leq
    \frac{\lambda_\alpha}{\min_{i\in[p]} \min_{t\in\cB_{\alpha,i}} \{\alpha_{t,i}^\ast\} - \max_{i\in[p]} \max_{t\in\cB_{\alpha,i}} |\wt\alpha_{t,i} - \alpha_{t,i}^\ast|} ,
\end{split}
\end{equation}
where the last inequality used \eqref{eqn: minmin_tilde_alpha}. Using the fact that $\min_{i\in[p]} \min_{t\in\cB_{\alpha,i}} \{\alpha_{t,i}^{\ast}\} =O_P(1)$ and Assumption~(R2), we have
\[
\lambda_\alpha = o_P\Big( \min_{i\in[p]} \min_{t\in\cB_{\alpha,i}} \{\alpha_{t,i}^{\ast 2}\} \Big) = o_P\Big( \min_{i\in[p]} \min_{t\in\cB_{\alpha,i}} \{\alpha_{t,i}^{\ast}\} \Big) ,
\]
so that in \eqref{eqn: DAFL_rate_cI_1}, $|\cI_1| =o_P(1)$. On the other hand, consider $\cI_2$. For any $t\in[T]$ and $i\in\cB_{\alpha,i}$, from the proof of Theorem~\ref{thm: general_initial_consistency_maineffect}, we can decompose
\begin{equation}
\label{eqn: DAFL_rate_cI_2}
\begin{split}
    \wt\alpha_{t,i} - \alpha_{t,i}^\ast &= q^{-1}\E_{t,i\cdot}^{\top}\1_q - q^{-1}\min\{\E_t\1_q\} \\
    &=
    q^{-1}(\A_{e,r})_{i\cdot}^{\top} \F_{e,t} \A_{e,c}^{\top} \1_q + q^{-1}\sum_{j=1}^q \Sigma_{\varepsilon,ij} \varepsilon_{t,ij} - q^{-1}\min\{\E_t\1_q\}.
\end{split}
\end{equation}
Since $\|\A_{e,r}\|_1, \|\A_{e,c}\|_1 = O(1)$ from Assumption~(E1), we have $q^{-1}(\A_{e,r})_{i\cdot}^{\top}\F_{e,t}\A_{e,c}^{\top}\1_q = O_P(q^{-1})$. Moreover, $\b{E}(q^{-1}(\bSigma_\epsilon \circ \bepsilon_t)\1_q) = \0$ and $\var(q^{-1}\sum_{j=1}^q\Sigma_{\epsilon,ij}\epsilon_{t,ij}) = q^{-2}\sum_{j=1}^q\Sigma_{\epsilon, ij}^2 = O(q^{-1})$, implying that $|q^{-1}((\bSigma_\epsilon \circ \bepsilon_t)\1_q)_i| =O_P(q^{-1/2})$. We also have $q^{-1}\min\{\E_t\1_q\} = O_P(q^{-1/2}\sqrt{\log(p)})$ from the proof of Theorem~\ref{thm: particular_initial_consistency_maineffect}. Finally, let $\gamma_{\alpha,i}^2 := \lim_{q\to\infty} q^{-1}\sum_{i=1}^q \Sigma_{\epsilon, ij}^2$, then using Theorem~1 in \cite{AyvazyanUlyanov2023}, we have $q^{-1/2}\sum_{j=1}^q \Sigma_{\epsilon,ij} \epsilon_{t,ij} \xrightarrow{\c{D}} \c{N}(0, \gamma_{\alpha,i}^2)$ and hence $|q^{-1} \sum_{j=1}^q \Sigma_{\varepsilon,ij} \varepsilon_{t,ij}| =O_P(q^{-1/2})$. We may now conclude from \eqref{eqn: DAFL_rate_cI_2} that $|\cI_2|=O_P(q^{-1/2} \sqrt{\log(p)})$. Hence finally,
\[
|\wh\alpha_{t,i} - \alpha_{t,i}^\ast| = |\cI_1 + \cI_2| =o_P(1),
\]
which completes the proof of the theorem.
\end{proof}

\subsection{Auxiliary results and proofs}

\begin{proof}[Proof of Corollary~\ref{cor: consistency_strong_factor}]
The results are direct from Theorem~\ref{thm: general_initial_consistency}.
\end{proof}

\begin{proof}[Proof of Corollary~\ref{cor: final_estimator}]
Result~1 is direct from Theorem~\ref{thm: oracle}, while result~2 is immediate by Lemma~\ref{lemma: alpha_tilde_max}, given the way we construct the final estimators.
\end{proof}

\begin{proof}[Proof of Proposition~\ref{prop: dgp_sparsity}]
For ease of notation, consider a time series $\{x_t\}_{t=1}^T$ with stay-in probabilities $\pi^{\cS}$ and $\pi^{\cB}$ and an initial probability $p_1=\b{P}(x_1=0)$. Let $p_t:=\b{P}(x_t=0)=\b{E}(\b1{\{x_t=0\}})$. Then for any $t=2,\dots,T$, we have
$$
p_{t+1} = p_t\pi^{\cS} + (1-p_t)(1-\pi^{\cB}) = (1-\pi^{\cB}) + (\pi^{\cS} + \pi^{\cB} - 1) p_t, 
$$
which can be solved recursively as
$$
p_t = p_\ast + (p_1-p_\ast)(\pi^{\cS} + \pi^{\cB} - 1)^{t-1},
$$
where $p_\ast=(1-\pi^{\cB})/(2-\pi^{\cS}-\pi^{\cB})$. We hence choose $p_1$ to be the same as $p_\ast$, and $p_t=p_\ast$ is satisfied. Furthermore, we shall compute the expected number of zeros over the whole series, showcased by
$$
\b{E}(\#\{t: x_t=0\}) = \sum_{t=1}^T p_t = Tp_\ast.
$$
Consequently, we can also compute the expected length of each sparse sub-block $\{t_\ell+1,\dots,t_\ell+m_{\ell}\}$. Every run of consecutive zeros is a geometric string, i.e., for $k=1,2,\dots$, $\b{P}(\text{block length} = k) = (\pi^{\cS})^{k-1}(1-\pi^{\cS})$. Hence the unconditional expected length of a sparse sub-block is 
$$
\b{E}(\text{length of a sparse sub-block}) = \sum_{k=1}^\infty k(1-\pi^{\cS})(\pi^{\cS})^{k-1} = \frac{1}{1-\pi^{\cS}}.
$$
This completes the proof of Proposition~\ref{prop: dgp_sparsity}.
\end{proof}

Our setting of factor structure is similar to that in \cite{CenLam2025} with tensor order $K=2$ therein. For convenience, we state their Proposition~1 (in their supplement) as the following lemma. 

\begin{lemma}\label{lemma: correlation_Et_Ft}
Let Assumptions (F1), (E1) and (E2) hold. Then
\begin{itemize}
    \item [(i)]
    (Weak correlation of noise $\E_t$ across different rows, columns and times). We only require Assumptions~(E1) and (E2) in this part.
    There exists some positive constant $C<\infty$ so that for any $t\in[T], i,j\in[p], h\in[q]$,
    \begin{align*}
    & \sum_{k=1}^{p}\sum_{l=1}^{q}
    \Big|\b{E}[E_{t,ih}
    E_{t,kl}]\Big| \leq C, \\
    & \sum_{l=1}^{q} \sum_{s=1}^T \Big| \textnormal{cov}(E_{t,ih} E_{t,jh}, E_{s,il} E_{s,jl} ) \Big| \leq C.
    \end{align*}
    \item [(ii)]
    (Weak dependence between factor $\F_t$ and noise $\E_t$). There exists some positive constant $C < \infty$ so that for any $j\in[p], i\in[q]$, and any deterministic vectors $\bf{u}\in\b{R}^{k_r}$ and $\bf{v}\in\b{R}^{k_c}$ with constant magnitudes,
    \begin{equation*}
    \b{E}\Bigg(
        \frac{1}{(qT)^{1/2}}\sum_{h=1}^{q}
        \sum_{t=1}^T
        E_{t,jh}
        \bf{u}^\top \F_t \bf{v}
    \Bigg)^2 \leq C,
    \;\;\;
    \b{E}\Bigg(
        \frac{1}{(pT)^{1/2}}\sum_{h=1}^{p}
        \sum_{t=1}^T
        E_{t,hi}
        \bf{v}^\top \F_t^\top \bf{u}
    \Bigg)^2 \leq C .
    \end{equation*}
    \item [(iii)]
    (Further results on factor $\F_t$).
    We only require Assumption~(F1) in this part.
    For any $t\in[T]$, all elements in $\F_t$ are independent of each other, with mean $0$ and unit variance. Moreover,
    \begin{equation*}
    \frac{1}{T}\sum_{t=1}^T \F_t\F_t^\top
    \xrightarrow{p} \bSigma_r := k_c\I_{k_r}
    , \;\;\;
    \frac{1}{T}\sum_{t=1}^T \F_t^\top\F_t
    \xrightarrow{p} \bSigma_c := k_r\I_{k_c} ,
    \end{equation*}
    with the number of factors $k_r$ and $k_c$ fixed as $\min\{T,p,q\}\to \infty$.
\end{itemize}
\end{lemma}

\begin{lemma}\label{lemma: alpha_tilde_max}
Let Assumptions~(IC1)
(with parts~(ii)--(iii) as in \eqref{eqn: IC1_min0}), 
(IC2), (E1), (E2), and (E3) hold. 
For each $i\in[p]$, define the notation $\balpha_{\cdot,i}^\ast = (\alpha_{1,i}^\ast, \dots, \alpha_{T,i}^\ast)^\top$ and $\wt\balpha_{\cdot,i} = (\wt\alpha_{1,i}, \dots, \wt\alpha_{T,i})^\top$, then we have
\[
\max_{i\in[p]} \big\|(\wt\balpha_{\cdot,i})_{\cS_{\alpha,i}} \big\|_{\infty} =
O_P\Big\{ q^{-1/2} \log^{1/2}\Big(p \sum_{i=1}^p |\cS_{\alpha,i}|\Big) \Big\},
\]
where $(\wt\balpha_{\cdot,i})_{\cS_{\alpha,i}}$ denotes the vector of $\wt\balpha_{\cdot,i}$ with indices restricted on $\cS_{\alpha,i}$, i.e., the vector consisting of $\{\wt\alpha_{t,i}\}_{t\in \cS_{\alpha,i}}$. Let $(\wt\balpha_{\cdot,i} -\balpha_{\cdot,i}^\ast)_{\cB_{\alpha,i}}$ be similarly defined by restricting indices on the set $\cB_{\alpha,i}$. Then we have
\begin{align*}
    \max_{i\in[p]} \big\|(\wt\balpha_{\cdot,i} -\balpha_{\cdot,i}^\ast)_{\cB_{\alpha,i}} \big\|_{\infty} =
    O_P\Big\{ q^{-1/2} \log^{1/2}\Big(p \sum_{i=1}^p |\cB_{\alpha,i}|\Big) \Big\} &.
\end{align*}
\end{lemma}

\begin{proof}[Proof of Lemma~\ref{lemma: alpha_tilde_max}]
To show the first result, from \eqref{eqn: MEFM} and Assumption~(IC1), we have for any $i\in[p]$,
\begin{align*}
    (\X_t \1_q)_i &= \big(\1_p (q\mu_t +\1_q^\top \bbeta_t^\ast) + q\balpha_t^\ast + \E_t \1_q \big)_i 
    =
    q\mu_t +\1_q^\top \bbeta_t^\ast + q \alpha_{t,i}^\ast + \1_q^\top \E_{t,i\cdot} ,
\end{align*}
so that using \eqref{eqn: tilde_balpha_t}, it holds for any $t\in[T]$,
\begin{equation}
\label{eqn: tilde_alpha_decomp}
\begin{split}
    &\hspace{13pt}
    (\wt\balpha_{\cdot,i})_t \equiv \wt\alpha_{t,i} = q^{-1} (\X_t \1_q)_i - q^{-1} \min\{ \X_t \1_q\} \\
    &=
    \mu_t + q^{-1} \1_q^\top \bbeta_t^\ast + \alpha_{t,i}^\ast + q^{-1} \1_q^\top \E_{t,i\cdot} - \min_{j\in[p]} \big\{ \mu_t + q^{-1} \1_q^\top \bbeta_t^\ast + \alpha_{t,j}^\ast + q^{-1} \1_q^\top \E_{t,j\cdot} \big\} \\
    &=
    \alpha_{t,i}^\ast + q^{-1} \1_q^\top \E_{t,i\cdot} -\min_{j\in[p]} \big\{ \alpha_{t,j}^\ast + q^{-1} \1_q^\top \E_{t,j\cdot} \big\} .
\end{split}
\end{equation}
Then using the min-zero identification condition in the main effects, and with \eqref{eqn: tilde_alpha_decomp}, we have
\begin{equation}
\label{eqn: tilde_alpha_max_step}
\begin{split}
    &\hspace{13pt}
    \big\| (\wt\alpha_{1,i} -\alpha_{1,i}^\ast, \dots, \wt\alpha_{T,i} -\alpha_{T,i}^\ast )_{\cS_{\alpha,i}}^\top \big\|_{\infty} 
    = \max_{t\in \cS_{\alpha,i}} \big| q^{-1} \1_q^\top \E_{t,i\cdot} -\min_{j\in[p]} \big\{ \alpha_{t,j}^\ast + q^{-1} \1_q^\top \E_{t,j\cdot} \big\}\big| \\
    &\leq
    \max_{t\in \cS_{\alpha,i}} \big| q^{-1} \1_q^\top \E_{t,i\cdot} \big| 
    + \max_{t\in \cS_{\alpha,i}} \max_{j\in[p]} \big| q^{-1} \1_q^\top \E_{t,j\cdot} \big|
    + 0
    \leq 2\max_{t\in \cS_{\alpha,i}} \max_{j\in[p]} \big| q^{-1} \1_q^\top \E_{t,j\cdot} \big|.
\end{split}
\end{equation}
For any $j\in[p]$, $t\in[T]$, by Assumption~(E1) and (E2),
\[
E_{t,jh} = \A_{e,r, j\cdot}^\top \Big( \sum_{w\geq 0}a_{e,w}\X_{e,t-w} \Big) \A_{e,c, h\cdot} + \Sigma_{\epsilon, jh} \Big( \sum_{g\geq 0} a_{\epsilon,g} X_{\epsilon,t-g,jh} \Big),
\]
so that we have
\begin{align*}
    &\hspace{13pt}
    \1_q^\top \E_{t,j\cdot} = \sum_{h=1}^q E_{t,jh} 
    =
    \A_{e,r, j\cdot}^\top \Big( \sum_{w\geq 0}a_{e,w}\X_{e,t-w} \Big) \sum_{h=1}^q \A_{e,c, h\cdot} + \sum_{h=1}^q \Sigma_{\epsilon, jh} \Big( \sum_{g\geq 0} a_{\epsilon,g} X_{\epsilon,t-g,jh} \Big) .
\end{align*}
By the sparsity of $\A_{e,c}$ from Assumption (E1), and Assumption (E2) and (E3), we conclude that the first term above is a zero-mean sub-Gaussian random variable with variance proxy of constant order. Similarly, the second term above is also zero-mean sub-Gaussian except that the variance proxy is of order $q$, and independent of the first term. Thus, $q^{-1}\1_q^\top \E_{t,j\cdot} \sim \text{subG}(C/q)$ for some arbitrary constant $C>0$, and hence for any $\varepsilon >0$ we have
\begin{align*}
    \b{P}\Big( \max_{i\in[p]} \max_{t\in \cS_{\alpha,i}} \max_{j\in[p]} \big| q^{-1} \1_q^\top \E_{t,j\cdot} \big| > \varepsilon\Big) \leq 2 \exp\Big\{ \log\Big(p \sum_{i=1}^p |\cS_{\alpha,i}|\Big) - q \varepsilon^2 /2C\Big\} ,
\end{align*}
implying that
\begin{align*}
    \max_{i\in[p]} \max_{t\in \cS_{\alpha,i}} \max_{j\in[p]} \big| q^{-1} \1_q^\top \E_{t,j\cdot} \big| =O_P\big\{ q^{-1/2} \log^{1/2}\big(p \sum_{i=1}^p |\cS_{\alpha,i}|\big) \big\} .
\end{align*}
Together with the equation \eqref{eqn: tilde_alpha_max_step} and the fact that $\alpha_{t,i}^\ast =0$ for any $i\in[p]$, $t\in \cS_{\alpha,i}$, we conclude the first result of Lemma~\ref{lemma: alpha_tilde_max}. The remaining result of the lemma can be shown by repeating all previous arguments, except that $\alpha_{t,i}^\ast$ is non-zero in \eqref{eqn: tilde_alpha_decomp}. This completes the proof of Lemma~\ref{lemma: alpha_tilde_max}.
\end{proof}

\begin{lemma}\label{lemma: loading_est_fulfill_iden}
Under Assumptions~(IC1) and (IC2),
both loading estimators $\wh\Q_r$ and $\wh\Q_c$ satisfy Assumption~(IC1)~(i) such that $\A_r^\top \1_p=\0$ and $\A_c^\top \1_q=\0$.
\end{lemma}

\begin{proof}[Proof of Lemma~\ref{lemma: loading_est_fulfill_iden}]
From the proof of Theorem~\ref{thm: general_initial_consistency},
we have
\[
\wt\L_t \wt\L_t^\top = \M_p\X_t\M_q\X_t^\top\M_p .
\]
Together with $\M_a \1_a = \0$ for any positive integer $a$, it is direct that the any eigenvector of $\wt\L_t\wt\L_t^\top$ is orthogonal to $\1_p$. 
This shows $\A_r^\top \1_p=\0$.
Similar arguments hold for $\wt\L_t^\top \wt\L_t$ and thus $\A_c^\top \1_q=\0$. 
This completes the proof.
\end{proof}


\end{document}